\documentclass[a4paper, reqno, 12pt]{amsart}
\usepackage[english]{babel}
\usepackage{amsmath,amssymb,amsthm, verbatim,bbm, float}
\usepackage{calrsfs}
\usepackage{dsfont}
\usepackage[pdftex]{xcolor}
\usepackage{graphicx}
\usepackage{tikz}
\usepackage{tabularx}
\usepackage{graphicx}
\usepackage[bookmarks=true, hyperindex, pdftex, colorlinks, allcolors=blue]{hyperref}
\usepackage[capitalize]{cleveref}
\usepackage{mathtools}
\usepackage{enumerate}
\usepackage{pgfmath} 

\usetikzlibrary{calc}

\numberwithin{equation}{section}

\usepackage{xspace}

\newtheorem{theorem}{Theorem}[section]
\newtheorem{corollary}[theorem]{Corollary}
\newtheorem{lemma}[theorem]{Lemma}
\newtheorem{remark}[theorem]{Remark}
\newtheorem{definition}[theorem]{Definition}
\newtheorem{proposition}[theorem]{Proposition}
\newtheorem{example}[theorem]{Example}
\newtheorem{question}[theorem]{Question}

\newtheorem{conjecture}[theorem]{Conjecture}

\AddToHook{env/corollary/begin}{\crefalias{theorem}{corollary}}
\AddToHook{env/lemma/begin}{\crefalias{theorem}{lemma}} 
\AddToHook{env/remark/begin}{\crefalias{theorem}{remark}}
\AddToHook{env/definition/begin}{\crefalias{theorem}{definition}}
\AddToHook{env/proposition/begin}{\crefalias{theorem}{proposition}}
\AddToHook{env/example/begin}{\crefalias{theorem}{example}}
\AddToHook{env/question/begin}{\crefalias{theorem}{question}}
\AddToHook{env/problem/begin}{\crefalias{theorem}{problem}}
\AddToHook{env/conjecture/begin}{\crefalias{theorem}{conjecture}}

\newcommand{\alp}{\mu_i}
\newcommand{\bet}{\nu_{i+n}}
\newcommand{\beti}{\nu_i}

\newcommand{\term}[1]{\textbf{#1}}
\newcommand{\tsym}{$\rho$-symmetric\xspace}
\newcommand{\Tsym}{Symmetric\xspace}
\newcommand{\tred}{minimal \tsym\xspace}

\newcommand{\tpairal}{$\rho$-paired alignment\xspace} 
\newcommand{\tpairals}{$\rho$-paired alignments\xspace}
\newcommand{\tcenteral}{$\rho$-central alignment\xspace} 
\newcommand{\tcenterals}{$\rho$-central alignments\xspace}

\newcommand{\cell}[1]{\Pi^{>0}_{#1}}
\newcommand{\symbridge}[1]{\beta_{#1}}
\newcommand{\Imin}[1]{I_{#1}^{min}}
\newcommand{\Imax}[1]{I_{#1}^{max}}

\DeclareMathOperator{\mat}{Mat}
\DeclareMathOperator{\antisymmat}{Persym}
\DeclareMathOperator{\tnnmat}{TNN}
\DeclareMathOperator{\rot}{\rho}
\DeclareMathOperator{\rotaff}{\widetilde{\rho}}

\DeclareMathOperator{\alt}{alt}

\DeclareMathOperator{\pairs}{pairs}
\DeclareMathOperator{\sing}{single}
\DeclareMathOperator{\rowspan}{rowspan}
\DeclareMathOperator{\spn}{span}
\DeclareMathOperator{\swap}{swap}
\DeclareMathOperator{\symmdim}{d_{\form}}

\DeclareMathOperator{\symmell}{\ell_{\form}}
\DeclareMathOperator{\mainop}{\mathcal{R}}
\DeclareMathOperator{\codim}{codim}
\DeclareMathOperator{\absell}{a\ell}

\DeclareMathOperator{\interior}{Int}

\newcommand{\rotcomp}[1]{\mainop(#1)}

\newcommand{\sless}{<_{R}}

\newcommand{\sleq}{\leq_{R}}

\newcommand{\M}{\mathcal{M}}
\newcommand{\B}{\mathcal{B}}
\newcommand{\Bkn}{\B(k,n)}
\newcommand{\Bnn}{\B(n,2n)}
\newcommand{\sbn}{\B^{\form}(n,2n)}
\newcommand{\sj}{\invol^{\form}}
\newcommand{\sbb}{\B^{\form}}
\newcommand{\embed}{\sigma}
\newcommand{\form}{R}

\newcommand{\symcross}[2]{\beta_{#1}(#2)}

\newcommand{\lollipop}[1]{L_{#1}}

\newcommand{\wf}{\widetilde{f}}

\usepackage{stackengine}

\newcommand{\celllgn}[1]{\Pi^{>0,{\form}}_{#1}}

\newcommand{\trsp}[1]{t_{#1} }
\newcommand{\plusn}[1]{#1 + n}

\newcommand{\twist}[1]{\rtimes{s_{#1}}}
\newcommand{\hultmantwist}[1]{\rtimes{s_{#1}}}
\newcommand{\affine}{\widetilde{S}_{2n}}
\newcommand{\naffine}{\widetilde{S}_{n}}
\newcommand{\affinen}{\widetilde{S}_{2n}^n}
\newcommand{\affinekn}{\widetilde{S}_{n}^k}
\newcommand{\affinezeron}{\widetilde{S}_{n}^0}
\newcommand{\invol}{\mathcal{J}}

\newcommand{\rvline}{\hspace*{-\arraycolsep}\vline\hspace*{-\arraycolsep}}

\DeclareMathOperator{\grrrr}{Gr}
\DeclareMathOperator{\grrnn}{Gr_{\geq 0}}
\DeclareMathOperator{\og}{OG}
\DeclareMathOperator{\ognnn}{OG_{\geq 0}}
\DeclareMathOperator{\lgrr}{LG^{\form}}
\DeclareMathOperator{\lgrnn}{LG^{\form}_{\geq 0}}

\newcommand{\ognn}{\ognnn(n,2n)}
\newcommand{\grknnn}{\grrnn(k,n)}
\newcommand{\grkn}{\grrrr(k,n)}
\newcommand{\gr}{\grrrr(n,2n)}
\newcommand{\grnn}{\grrnn(n,2n)}
\newcommand{\RR}{\mathbb{R}}

\newcommand{\ZZ}{\mathbb{Z}}

\newcommand{\lgr}{\lgrr (n,2n)}
\newcommand{\lgnn}{\lgrnn (n,2n)}
\newcommand{\lgnp}{LG^{\form}_{> 0}(n,2n)}

\newcommand{\Spn}{Sp(2n)}
\newcommand{\Spalgn}{\mathfrak{sp}(2n)}

\usepackage{subcaption}

\newcommand{\bdot}[2]{
\draw[line width=0.3 mm,blue] [black, fill = black] (#1, #2) circle [radius = 0.1];
}

\newcommand{\wdot}[2]{
\draw[line width=0.3 mm,blue] [black, fill = white] (#1, #2) circle [radius = 0.1];
}
\newcommand{\cntr}{\draw[line width=0.3 mm,red] [red, fill = red] (0, 0) circle [radius = 0.2];}

\newcommand{\posetscale}{0.6}
\newcommand{\posetline}{0.8}
\newcommand{\posetcircle}{0.8}

\newcommand{\constructionscale}{0.4}
\newcommand{\constructionline}{1}
\newcommand{\constructioncircle}{1.4}

\newcommand{\graphline}{1}
\newcommand{\graphcircle}{1}
\newcommand{\graphscale}{0.5}

\newcommand{\bridgeline}{1}
\newcommand{\bridgecircle}{1.7}
\newcommand{\bridgescale}{0.6}

\newcommand{\topcellscale}{0.4}

\newcommand{\coordx}[1]{cos(90*(2-#1}
\newcommand{\coordy}[1]{sin(90*(2-#1}

\newcommand{\coordxsix}[1]{{3*(cos(60*(2-#1)))}}
\newcommand{\coordysix}[1]{{3*(sin(60*(2-#1)))}}

\newcommand{\loosns}{7}

\def\levl{2.5}
\def\levll{5}
\def\levlll{7.5}
\newcommand\edge[4]{
\draw[line width=0.5] (#1#3.south)--(#2#4.north);
}
\newcommand{\circlenodes}{
        \node[above] at (\coordx{1},\coordy{1}) (1) {};
        \node[right] at (\coordx{2},\coordy{2}) (2) {};
        \node[below] at (\coordx{3},\coordy{3}) (3) {};
        \node[left] at (\coordx{4},\coordy{4}) (4) {};
}

\newcommand{\circlenodessix}{
    \node [above right] at ({3*cos(60)},{3*sin(60)}) {1};
    \node [right] at ({3*cos(0)},{3*sin(0)}) {2};
    \node [below right] at ({3*cos(-60)},{3*sin(-60)}) {3};
    \node [below left] at ({3*cos(-120)},{3*sin(-120)}) {4};
    \node [left] at ({3*cos(180)},{3*sin(180)}) {5};
    \node [above left] at ({3*cos(120)},{3*sin(120)}) { 6};
}

\newcommand{\symmlineblue}[5]{
     \draw[line width=#5, blue] (#1,#2) -- (#3,#4);
     \draw[line width=#5, blue] (-#1,-#2) -- (-#3,-#4);
}

\newcommand{\symmbdot}[2]{
    \bdot{#1}{#2};
    \wdot{{-#1}}{{-#2}};
}
\newcommand{\symmwdot}[2]{
    \wdot{#1}{#2};
    \bdot{{-#1}}{{-#2}};
}

\newcommand{\symmbdotsix}[1]{
    \symmbdot{\coordxsix{#1}}{\coordysix{#1}}
}
\newcommand{\symmwdotsix}[1]{
    \symmwdot{\coordxsix{#1}}{\coordysix{#1}}
}

\newcommand{\twonchoosen}{{{[2n]}\choose n}}

\usepackage{environ} 

\DeclareMathAlphabet{\mathcal}{OMS}{cmsy}{m}{n}

\newcommand{\st}{\,|\,}

\title{Rotationally symmetric plabic graphs and the Lagrangian Grassmannian}
\author{Olha Shevchenko}

\begin{document}

\maketitle

\begin{abstract}
    We introduce the totally nonnegative Lagrangian Grassmannian $\lgnn$, a new subset of the totally nonnegative Grassmannian consisting of subspaces isotropic with respect to a certain bilinear form $R$. We describe its cell structure and show that each cell admits a representation by a rotationally symmetric (not necessarily reduced) plabic graph. Along the way, we develop new techniques for working with non-reduced plabic graphs.
\end{abstract}

\setcounter{tocdepth}{1}
\tableofcontents

\section{Introduction}



The \term{totally nonnegative (TNN) Grassmannian} \( \grknnn \), introduced in \cite{postnikov, lusztig}, is the subset of the real Grassmannian \( \grkn \) consisting of \( k \)-dimensional subspaces of \( \mathbb{R}^n \) that can be represented by matrices with all nonnegative Plücker coordinates. Since its introduction, the TNN Grassmannian has gathered significant interest due to its rich combinatorial structure and connections to other mathematical areas, particularly through its relationship with \term{plabic graphs}---planar bicolored graphs that parametrize each cell in \( \grknnn \). We introduce the relevant background in \Cref{background-section}.

Within the totally nonnegative Grassmannian, one can define the totally nonnegative orthogonal and Lagrangian Grassmannians, the sets of subspaces isotropic with respect to a bilinear form. For particular choices of the bilinear form, such subspaces of $\grknnn$ have been shown to be related to objects in statistical mechanics such as the Ising model \cite{pashas} and electrical resistor networks \cite{chepuri-george-speyer, lam-electroids}. Furthermore, these subspaces can often be represented by plabic graphs that satisfy particular symmetry conditions, as demonstrated in \cite{chepuri-george-speyer, pashas, karpman-su-combinatorics}.

We compile these previous results in a table similar to the one presented in \cite{george-twist}.

\begin{table}[h!]
    \centering
    \begin{tabularx}{\textwidth}{|X|X|X|}
    \hline
        \textbf{Subspace of $\boldsymbol{\grknnn}$} & \textbf{Plabic graphs} & \textbf{Application}\\ \hline
        TNN Grassmannian \cite{postnikov,lusztig} &  All graphs & Bipartite dimer model \\  \hline
        TNN orthogonal Grassmannian \cite{pashas} &  Pseudo-symmetric with respect to the color swap & Ising model \\ \hline
        TNN Lagrangian Grassmannian \cite{chepuri-george-speyer, lam-electroids} & Temperley's graphs & Electrical networks\\ \hline
        TNN Lagrangian Grassmannian \cite{karpman-su-combinatorics} & Mirror-symmetric &  Lusztig's TNN partial flag varieties of type C  \cite{lusztig}\\ \hline
    \end{tabularx}
    \label{applications-table}
\end{table}

A natural question arises.
\begin{question}
    \label{question: intro}
    Which other choices of bilinear form yield interesting combinatorial structures? Are there other forms that would correspond to plabic graphs exhibiting some symmetry?
\end{question}

It turns out that such well-behaved structures are quite rare. For many quadratic forms, the associated totally nonnegative orthogonal or Lagrangian Grassmannian either fails to be full dimensional or even has no real points at all. For example, \cite[Example~1.2]{yelena-yassine-orthogonal-grassmannian} provides explicit instances with no real points. Moreover, \cite[Proposition~4.2.4]{yelena-dissertation} (originally due to Yassine El Maazouz and Yelena Mandelshtam, as explained there) shows that the totally nonnegative orthogonal Grassmannian considered in \cite{pashas} is the only full-dimensional TNN orthogonal Grassmannian corresponding to a diagonal form --- namely, the only diagonal case for which $\dim \ognn = \dim \og(n,2n)$.

In this paper, we make progress on \Cref{question: intro} by introducing a new Lagrangian Grassmannian $\lgnn$ and showing that it corresponds to the (non-reduced) rotationally-symmetric plabic graphs.

\subsection{Main results}
The primary object of study in this paper is a \term{Lagrangian Grassmannian}, corresponding to a bilinear form $\form$ not previously explored. This Lagrangian Grassmannian is defined as 
\begin{equation}
    \label{lgr-def-equation}
    \lgr = \{X \in \gr\st \form(x,y) = 0 \text{ for all } x,y \in X\},
\end{equation}
where $\form$ is the skew-symmetric bilinear form

\begin{equation}
    \label{Q-def-equation}
    \form(x,y) = \sum_{i=1}^n (-1)^{i-1} (x_i y_{n+i} - x_{n+i}y_i).
\end{equation}

We are interested in the \term{totally nonnegative Lagrangian Grassmannian}
\begin{equation}
    \label{lgnn-def-equation}
    \lgnn = \lgr \cap \grnn.
\end{equation}

We define the Lagrangian positroid cells as $\celllgn{f}= \cell{f} \cap \lgnn$, indexed by bounded affine permutations $f$. We call a permutation $f$ \term{\tsym} if the map $i \mapsto f(i) - n$ is an involution, and denote the set of such permutations as $\sbn$. For each $f \in \sbn$, we denote $\symmell(f)$ to be the symmetric length of $f$ (see \Cref{def: symmell-symmdim}). We call a plabic graph $G$ \term{\tsym} if it is equal to the color swap of its $180^\circ$ rotation. See \Cref{symm-perm-section} and \Cref{symm-graphs-section} for definitions.

We are now ready to state our main results.
\begin{theorem}
\label{main-theorem}
    \begin{enumerate}
        \item $\lgnn = \bigsqcup\limits_{f \in \sbn} \celllgn{f}$, where the union is taken over all \tsym affine permutations $f$.
        \item For $f \in \sbn$, 
        the cell $\celllgn{f}$ is homeomorphic to $\RR^{{{n+1} \choose 2}-\symmell(f)},$ 
        and can be parametrized by a \tsym, non-reduced, plabic graph.

        \item For two \tsym permutations $f$ and $g$, 
        $$\cell{f} \subset \overline{\cell{g}} \quad \text{if and only if } \quad f \leq g \text{ in the affine Bruhat order}.$$
        Furthermore,
        \[
        \overline{\cell{g}} = \bigsqcup_{\substack{f \leq g\\
        f \in \sbn}} \cell{f}.
        \]

        \item $\lgnn \simeq [0,1]^{{n+1 \choose 2}}$, i.e. $\lgnn$ is homeomorphic to a closed ball. 
    \end{enumerate}
\end{theorem}

This result is parallel to results by \cite{postnikov} (see \Cref{postnikov-theorem}) and \cite[Theorem~1.1]{galashin-karp-lam-ball}, which prove similar results for $\grknnn$. We will prove parts \textit{(1)} and \textit{(2)} in \Cref{symm-graphs-section} and parts \textit{(3)} and \textit{(4)} in \Cref{cell-structure-section}.

We further observe that by part \textit{(1)} of \Cref{main-theorem}, points of $\lgnn$ correspond to affine permutations which, after a shift by $n$, are involutions. Consequently, $\lgnn$ provides a geometric realization of an order ideal in the poset of affine involutions. While this poset has been investigated in depth in previous work --- see \cite{hultman, hansson-hultman-involutions, hamaker-marberg-pawlowski-involutions, marberg-involutions} --- a geometric interpretation did not exist before. Additional details appear in \Cref{subsection: symm-perm-involutions}.

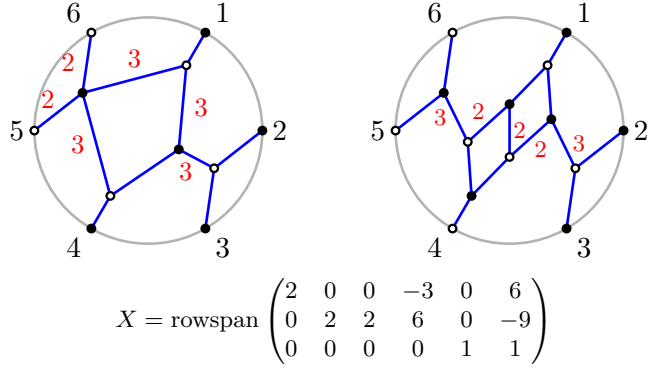
\begin{figure}
    \centering
    \begin{center}

\begin{tabular}{ccc}

    {\begin{tikzpicture}[scale = \graphscale, baseline={(0,0)}]
    \footnotesize
    \draw[line width=\graphcircle,opacity=0.3] (0,0) circle (3);
    \draw[line width=\graphline, blue] ({3*cos(60)},{3*sin(60)}) --
    ({2*cos(60)},{2*sin(60)})--node[above,red,inner sep=4pt]{3}({-2*cos(-30)},{-2*sin(-30)})--node[left,red,inner sep=4pt]{3}({-2*cos(60)},{-2*sin(60)});
    \draw[line width=\graphline, blue] ({-3*cos(60)},{-3*sin(60)}) --
    ({-2*cos(60)},{-2*sin(60)})--
    (.8,-.5)--node[right,red,inner sep=4pt]{3}({2*cos(60)},{2*sin(60)});
    \draw[line width=\graphline, blue] (.8,-.5) --node[below left,red,inner sep=1pt]{3} ({2*cos(-30)},{2*sin(-30)});
    \draw[line width=\graphline, blue] ({3*cos(0)},{3*sin(0)}) --
    ({2*cos(-30)},{2*sin(-30)}) --
    ({3*cos(-60)},{3*sin(-60)});
    \draw[line width=\graphline, blue] ({-3*cos(0)},{-3*sin(0)}) --node[above left,red,inner sep=1pt]{2} ({-2*cos(-30)},{-2*sin(-30)}) --node[left,red,inner sep=4pt]{2} ({-3*cos(-60)},{-3*sin(-60)});
    
    \bdot{{3*cos(60)}}{{3*sin(60)}}; \bdot{{3*cos(0)}}{{3*sin(0)}}; \bdot{{3*cos(-60)}}{{3*sin(-60)}}; \bdot{{3*cos(-120)}}{{3*sin(-120)}}; \wdot{{3*cos(180)}}{{3*sin(180)}}; \wdot{{3*cos(120)}}{{3*sin(120)}};\wdot{{2*cos(60)}}{{2*sin(60)}};\wdot{{2*cos(-120)}}{{2*sin(-120)}};\wdot{{2*cos(-30)}}{{2*sin(-30)}};\bdot{{-2*cos(-30)}}{{-2*sin(-30)}}; \bdot{.8}{-.5};
    \small
    \node [above right] at ({3*cos(60)},{3*sin(60)}) {1};
    \node [right] at ({3*cos(0)},{3*sin(0)}) {2};
    \node [below right] at ({3*cos(-60)},{3*sin(-60)}) {3};
    \node [below left] at ({3*cos(-120)},{3*sin(-120)}) {4};
    \node [left] at ({3*cos(180)},{3*sin(180)}) {5};
    \node [above left] at ({3*cos(120)},{3*sin(120)}) { 6};
    \normalsize
\end{tikzpicture}}
    &
    {
    }
    &
    {\begin{tikzpicture}[scale = \graphscale,baseline={(0,0)}]
    \scriptsize
    \draw[line width=\graphcircle,opacity = 0.3] (0,0) circle (3);
    \draw[line width=\graphcircle, blue] ({3*cos(60)},{3*sin(60)}) -- ({2*cos(60)},{2*sin(60)})--
    (0,.7)--node[above left,red,inner sep=1pt]{2}(-1.1,-.3)--
    ({-2*cos(60)},{-2*sin(60)});
    \draw[line width=\graphcircle, blue] ({-3*cos(60)},{-3*sin(60)}) -- ({-2*cos(60)},{-2*sin(60)})--
    (0,-.7)--node[below right,red,inner sep=1pt]{2}(1.1,.3)--
    ({2*cos(60)},{2*sin(60)});
    \draw[line width=\graphcircle, blue] (0,.7)--node[right,red,inner sep=1pt]{2}(0,-.7);
    \draw[line width=\graphcircle, blue] (1.1,.3) --node[right,red,inner sep=3pt]{3} ({2*cos(-30)},{2*sin(-30)});
    \draw[line width=\graphcircle, blue] (-1.1,-.3) --node[left,red,inner sep=3pt]{3} ({-2*cos(-30)},{-2*sin(-30)});
    \draw[line width=\graphcircle, blue] ({3*cos(0)},{3*sin(0)}) -- ({2*cos(-30)},{2*sin(-30)}) -- ({3*cos(-60)},{3*sin(-60)});
    \draw[line width=\graphcircle, blue] ({-3*cos(0)},{-3*sin(0)}) -- ({-2*cos(-30)},{-2*sin(-30)}) -- ({-3*cos(-60)},{-3*sin(-60)});
    
    \bdot{{3*cos(60)}}{{3*sin(60)}}; \bdot{{3*cos(0)}}{{3*sin(0)}}; \bdot{{3*cos(-60)}}{{3*sin(-60)}}; \wdot{{3*cos(-120)}}{{3*sin(-120)}}; \wdot{{3*cos(180)}}{{3*sin(180)}}; \wdot{{3*cos(120)}}{{3*sin(120)}};\wdot{{2*cos(60)}}{{2*sin(60)}};\bdot{0}{.7}; \wdot{-1.1}{-.3};\bdot{{2*cos(-120)}}{{2*sin(-120)}};\wdot{{2*cos(-30)}}{{2*sin(-30)}};\bdot{{-2*cos(-30)}}{{-2*sin(-30)}};\wdot{0}{-.7}; \bdot{1.1}{.3};
    \small
    \node [above right] at ({3*cos(60)},{3*sin(60)}) {1};
    \node [right] at ({3*cos(0)},{3*sin(0)}) {2};
    \node [below right] at ({3*cos(-60)},{3*sin(-60)}) {3};
    \node [below left] at ({3*cos(-120)},{3*sin(-120)}) {4};
    \node [left] at ({3*cos(180)},{3*sin(180)}) {5};
    \node [above left] at ({3*cos(120)},{3*sin(120)}) { 6};
    \normalsize
\end{tikzpicture}}
    \end{tabular}\\[1mm]
    \scriptsize
    $X = \rowspan \begin{pmatrix}
        2 & 0 & 0 & -3 & 0 & 6\\
        0 & 2 & 2 & 6 & 0 & -9\\
        0 & 0 & 0 & 0 & 1 & 1
    \end{pmatrix}$
    \normalsize
    \end{center}
    \normalsize
    \caption{The two networks represent the same point $X$. The one on the left is reduced, but not \tsym, while the one on the right is \tsym, but not reduced. Here, unmarked edges have weight $1$.}
    \label{two-graphs-figure}
\end{figure}

\begin{figure}
    \centering
    \begin{center}

\begin{tabular}{ccc}

    {
\begin{tikzpicture}[scale = \topcellscale, baseline={(0,0)}]
        \draw[line width=\graphcircle,opacity=0.3] (0,0) circle (3);

        \draw[line width=\graphline,blue] (2.12132034356,2.12132034356) -- (1,1) -- (1,-1) -- (2.12132034356,-2.12132034356);
        \draw[line width=\graphline,blue] (-2.12132034356,2.12132034356) -- (-1,1) -- (-1,-1) -- (-2.12132034356,-2.12132034356);
        \draw[line width=\graphline,blue] (1,1) -- (-1,1);
        \draw[line width=\graphline,blue] (1,0) -- (-1,0);
        \draw[line width=\graphline,blue] (1,-1) -- (-1,-1);

        \wdot{2.12132034356}{2.12132034356}; \wdot{2.12132034356}{-2.12132034356};\bdot{-2.12132034356}{-2.12132034356};\bdot{-2.12132034356}{2.12132034356};\wdot{-1}{1};\bdot{1}{1};\bdot{-1}{0};\wdot{1}{0};\wdot{-1}{-1};\bdot{1}{-1};

        \scriptsize
        \node [above right] at (2.12132034356,2.12132034356) {1};
        \node [below right] at (2.12132034356,-2.12132034356) {2};
        \node [below left] at (-2.12132034356,-2.12132034356) {3};
        \node [above left] at (-2.12132034356,2.12132034356) {4};
        
\end{tikzpicture}

}
    &
    {
    }
    &
    {
\newcommand{\xcoord}{0.776457135308}

\begin{tikzpicture}[scale = \topcellscale, baseline={(0,0)}]
        \draw[line width=\graphcircle,opacity=0.3] (0,0) circle (3);

        \draw[line width=\graphline,blue] (3,0) -- (2,0) -- (\xcoord,1) -- (-\xcoord,1) -- (-2,0) -- (-3,0);

        \draw[line width=\graphline,blue] (2,0) -- (\xcoord,-1) -- (-\xcoord,-1) -- (-2,0);

        \draw[line width=\graphline,blue] (\xcoord,-2.89777747887) -- (\xcoord,2.89777747887);

        \draw[line width=\graphline,blue] (-\xcoord,-2.89777747887) -- (-\xcoord,2.89777747887);

        \draw[line width=\graphline,blue] (\xcoord,2) -- (-\xcoord,2);
        \draw[line width=\graphline,blue] (\xcoord,-2) -- (-\xcoord,-2);
        \draw[line width=\graphline,blue] (\xcoord,0) -- (-\xcoord,0);
        
        \bdot{3}{0}; \wdot{-3}{0};\bdot{{3*cos(285)}}{{3*sin(285)}};\wdot{{3*cos(255)}}{{3*sin(255)}};\bdot{{3*cos(75)}}{{3*sin(75)}};\wdot{{3*cos(105)}}{{3*sin(105)}};\bdot{-2}{0};\wdot{2}{0};\bdot{-\xcoord}{0};\wdot{\xcoord}{0};\bdot{-\xcoord}{-2};\wdot{\xcoord}{-2};\bdot{-\xcoord}{2};\wdot{\xcoord}{2};
        \wdot{-\xcoord}{1};\bdot{\xcoord}{1};
        \wdot{-\xcoord}{-1};\bdot{\xcoord}{-1};

        \scriptsize
        \node [right] at (3,0) {1};
        \node [below right] at (\xcoord,-2.89777747887) {2};
        \node [below left] at (-\xcoord,-2.89777747887) {3};
        \node [left] at (-3,0) {4};
        \node [above left] at (-\xcoord,2.89777747887) {5};
        \node [above right] at (\xcoord,2.89777747887) {6};
\end{tikzpicture}

}\vspace{0.2cm} \\ 
    {$n=2$} & & {$n=3$}
\end{tabular}
\end{center}
\normalsize
    \caption{The \tsym graphs parametrizing the top cells for $n=2$ and $3$.}
    \label{fig:top-cells-2-3}
\end{figure}
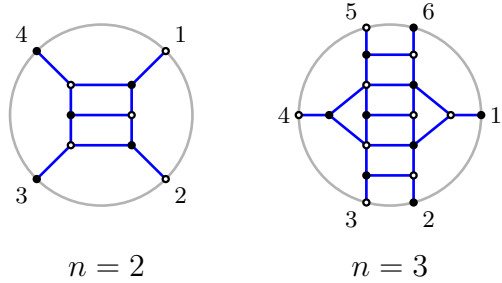

\begin{example}
	\label{alm-symm-vs-symm-example}
	To demonstrate the second part of \Cref{main-theorem}, consider
	$$X = \rowspan \begin{pmatrix}
		2 & 0 & 0 & -3 & 0 & 6\\
		0 & 2 & 2 & 6 & 0 & -9\\
		0 & 0 & 0 & 0 & 1 & 1
	 \end{pmatrix} \in \lgrnn(3,6).$$
    Following the standard approach of \cite{postnikov}, $X$ can be represented by a reduced plabic graph, which does not exhibit any obvious symmetry (see \Cref{two-graphs-figure} (left)). However, \Cref{main-theorem} allows us to instead represent $X$ with a symmetric plabic graph that has symmetric edge weights, as shown in \Cref{two-graphs-figure} (right). This alternative representation is \tsym but is no longer reduced.

    Another example illustrating this fact appears in \Cref{fig:top-cells-2-3}, where we show the \tsym graphs representing the top cells in $\lgnn$ for $n=2$ and $3$.

\end{example}

\subsection{\Tsym graphs are not always reduced}

One significant difference from previous results is that, in our setting, plabic graphs typically cannot be both reduced and \tsym (see \Cref{alm-symm-vs-symm-example}). Specifically, \Cref{prop:condition-symm-and-reduced} shows that a graph can be both reduced and \tsym only if it corresponds to a permutation $f$ satisfying
$$
\#\{i \in [2n] \mid f(i) = i+n\} \in \{0, 2\},
$$
where $[2n] = \{1,2,\ldots, 2n\}$. In particular, for $n > 1$, the top cell cannot be represented by a graph that is simultaneously reduced and \tsym. Further details appear in \Cref{subsection: tsym-reduced}.

In contrast, earlier works (for example, \cite{chepuri-george-speyer, pashas, karpman-su-combinatorics}) relied on plabic graphs that are both reduced and exhibited certain symmetries. To accommodate the failure of reducedness in our setting, we develop new techniques (for example, see \Cref{cell-structure-section}).

\subsection{$\lgnn$ through totally nonnegative symmetric matrices and Pl\"ucker coordinates}

We show in \Cref{minor-condition-on-lg} that the TNN Lagrangian Grassmannian can be also defined through a certain condition on minors. Specifically,
$$\lgnn = \{X \in \grnn\st \Delta_I(X) = \Delta_{I^C + n}(X) \text{ for all } I \subset [2n], |I| = n\},$$
where $I^C = [2n] \setminus I$ and $I+n = \{i + n \pmod {2n}\st i \in I\} \subset [2n]$.

It also turns out that an open dense subset of $\lgnn$ --- specifically, the Schubert cell $\{X\in \lgnn\st \Delta_{[n]}(X) \neq 0\}$ --- is in bijection with the set of $n\times n$ matrices that have all nonnegative minors and are symmetric with respect to reflection across the anti-diagonal. Specifically, $X = \rowspan \begin{pmatrix}
    I \st M
\end{pmatrix} \in \lgnn$ if and only if the matrix $M'$ obtained from $M$ by reversing the order of rows and then changing the signs of even rows, is symmetric with respect to the anti-diagonal and is totally nonnegative (see \Cref{prop: lgnn-vs-persymm-matrices}). Thus, $\lgnn$ may be viewed as a compactification of the space of such matrices.


\section{Background}\label{background-section}

\subsection{Totally nonnegative Grassmannian}
    The \term{Grassmannian} $\grkn$ is the set of all $k$-dimensional subspaces of $\RR^n$. Each $X \in \grkn$ can be represented as a rowspan of a $k \times n$ matrix, the rows of which form a basis of $X$. That is,
	$$\grkn = \{X \subset \RR^n\st \dim X = k\} \simeq \mat(k,n)/\text{row operations}.$$
	
    Each point in $\grkn$ is uniquely determined by its \term{Pl\"ucker coordinates}, determinants of the maximal $k \times k$ minors of the corresponding matrix (defined up to common rescaling). These coordinates are denoted by $\Delta_I$, where $I \in {[n]\choose k}$ is a column index set.
	
	Finally, we define the \term{totally nonnegative Grassmannian} $\grknnn$ as the subspace of the Grassmannian with all non-negative Pl\"ucker coordinates. That is,
    $$\grknnn = \{X \in \grkn\st \Delta_I(X) \geq 0 \text{ for all } I \in {[n] \choose k}\}.$$
	\begin{example} Take
		$$X = \rowspan \begin{pmatrix}
			1 & 1 & 3 & 2\\
			0 & 1 & 4 & 5
		\end{pmatrix}.$$
		Then 
		$$\Delta_{12} = 1, \Delta_{13} = 4, \Delta_{14} = 5, \Delta_{23} = 1, \Delta_{24} = 3, \Delta_{34} = 7,$$
		so $X \in \grrnn(2,4)$.
	\end{example}

\subsection{Permutations and positroid cells}

	For each point $X \in \grknnn$ we can define the corresponding bounded affine permutation $f_X: \ZZ \to \ZZ$ (see \cite{knutson-lam-speyer})
	$$f_X(i)= \min\{j\geq i \st x_i \in \spn\{x_{i+1}, \ldots, x_j\}\},$$
	where $x_i$ is $i$-th column of a matrix representing $X$ and we use the convention $x_{i+n} = x_i$.

	The permutation $f_X$ captures the structure of column dependencies in $X$, uniquely determining which of its Pl\"ucker coordinates are zero. We combine all points 	$X$ corresponding to one permutation $f$ into a \term{positroid cell}, defined as
	$$\cell{f} = \{X \in \grknnn \st f_X = f\}.$$

    It is known that the codimension of $\cell{f}$ is $\ell(f)$, the number of inversions of $f$.
    \begin{equation*}
        \ell(f) = |\{(i,j) \in [n] \times \ZZ\st i<j, f(i)>f(j)\}|.
    \end{equation*}
    We also define 
    $$d(f) = \dim \cell{f} = k(n-k) - \ell(f).$$

    These positroid cells partition $\grknnn$ and play a central role in the combinatorial structure of the totally nonnegative Grassmannian.

\subsection{Bruhat order on permutations}
\label{subsection: background-bruhat-order}

    Call a bijection $f: \ZZ \to \ZZ$ an affine permutation if for each $i \in \ZZ$, $f(i + n) = f(i) + n$.

    Consider the set $\affinekn$ of $(k,n)$-affine permutations, i.e. bounded affine permutations $f: \ZZ \to \ZZ$ such that $i \leq f(i) \leq i+n$ and $\sum\limits_{i=1}^{n} (f(i) - i) = kn$. We equip it with the \term{circular Bruhat order} $\leq$ (see \cite{postnikov}), dual to the usual order in $\affinekn$. We then consider the subset $\Bkn$ of the bounded affine permutations, satisfying $i \leq f(i) \leq i+n$, $f(i + n) = f(i) + n$ for each $i$.

    Finally, let $\naffine = \affinezeron$ be the affine Coxeter group. 

    \begin{remark}
        Notably, since the circular Bruhat order is dual to the usual order in $\affinekn$, this means that $\ell(f) \geq \ell(g)$ whenever $f \leq g$. The maximal permutation in $\Bkn$ is $f: i \mapsto i + k$. It is also important to point out that $\Bkn$ is the upper order ideal in $\affinekn$.
    \end{remark}
    The following is a simple property of affine permutations.
    \begin{lemma}
        \label{condition-fsi-greaterthan-f-lemma}
        Suppose $f \in \Bkn$ and $i \in \ZZ$. Then $fs_i > f$ if and only if we have $i < i+1 \leq f(i+1) < f(i) \leq i+n$.
    \end{lemma}

    \begin{definition}

        Suppose $f \in \Bkn$ and $i < j$.
        \begin{itemize}
            \item Say that $i$ and $j$ form a \term{crossing} if 
            \begin{equation}\begin{aligned}
            \label{eq: crossing-def}
             i < j \leq f(i) < f(j) \leq i + n.
            \end{aligned}
            \end{equation}

            \item Call a crossing between $i$ and $j$ \term{simple} if 
             \begin{equation}
                \begin{aligned}
                    \label{eq:simple-crossing-def}
                    \text{there are no } \ell  \text{ such that } i < \ell < j \leq f(i) < f(\ell) < f(j) \leq i + n.
                \end{aligned}
            \end{equation}

            \item Say that $i$ and $j$ form an \term{alignment} in $f$ if
            \begin{equation}\begin{aligned}
                \label{eq: alignment-def}
                 i < j \leq f(j) < f(i).
            \end{aligned}
            \end{equation}

            \item Finally, we say that an $i$ and $j$ form a \term{simple alignment} if they form an alignment and 
            \begin{equation}
                \begin{aligned}
                    \label{eq:simple-alignment-def}
                    \text{there are no } \ell  \text{ such that } i < \ell < j \leq f(j) < f(\ell) < f(i).
                \end{aligned}
            \end{equation}
        \end{itemize}

        These are defined similarly in the case $i > j$.
    \end{definition}

    \begin{lemma}[\cite{postnikov}, Theorem~17.8]
        \label{lemma:covering-relation}
        Suppose $f, g \in \Bkn$. Then $f \lessdot g$ if and only if $g = f \trsp{i,j}$, where $i$ and $j$ form a simple alignment in $f$ (or, equivalently, a simple crossing in $g$).
    \end{lemma}

    \begin{definition}
        Suppose $f \in \Bkn$ and $i \in \ZZ$. Define 
        $$f \star s_i = \max\{f, f  s_i\} \text{ and } s_i \star f = \max\{f, s_i  f\}.$$
    \end{definition}
    \begin{remark}
        Notice that since $\Bkn$ is an upper order ideal (see \cite[Lemma~3.6]{knutson-lam-speyer}), $f \star s_i$ and $s_i \star f$ are always in $\Bkn$ (even though $fs_i$ might be not).
    \end{remark}

\subsection{Plabic graphs}\label{plabic-graphs-subsection}

    A \term{plabic graph} is a planar bipartite graph embedded in a disk. In these graphs, the boundary vertices are labeled from $1$ to $n$ in clockwise order, each of degree at most $1$. A \term{plabic network} is a plabic graph with positive edge weights.
	
    \begin{theorem}[\cite{postnikov}, Theorem 4.8]
        Every point in $\grknnn$ can be represented by a plabic network with $n$ boundary vertices.
    \end{theorem}
	
    The map from weighted plabic graphs to $\grknnn$ is called the \term{boundary measurement map}. We denote the point in $\grknnn$ represented with a network $N = (G,w)$ as $X(N) = X(N,w)$. We refer the reader to \cite{lam-textbook} for the construction of the map.

    \begin{proposition}
        Let $G$ be a plabic graph. Then all points $X(G,w)$ for different weightings $w$ of edges of $G$ correspond to the same permutation $f$. We denote $f_G = f$ and $\cell{G} = \cell{f_G}$
    \end{proposition}

\begin{figure}
        \centering
        \newcommand{\movesline}{1}
\newcommand{\movesscale}{0.6}
\begin{tabular}{ccccccc}
   {    \begin{tikzpicture}[scale = 0.4,baseline={(0,0)}]
        \draw[line width=\movesline, blue] (2,2) -- (0.85, 0.85);
        \draw[line width=\movesline, blue] (2,-2) -- (0.85, -0.85);
        \draw[line width=\movesline, blue] (-2,2) -- (-0.85, 0.85);
        \draw[line width=\movesline, blue] (-2,-2) -- (-0.85, -0.85);
        \draw[line width=\movesline, blue] (0.85,-0.85) -- (-0.85, -0.85);
        \draw[line width=\movesline, blue] (0.85,0.85) -- (0.85, -0.85);
        \draw[line width=\movesline, blue] (0.85,0.85) -- (-0.85, 0.85);
        \draw[line width=\movesline, blue] (-0.85,0.85) -- (-0.85, -0.85);
        \bdot{0.85}{0.85}; \wdot{-0.85}{0.85}; \wdot{0.85}{-0.85}; \bdot{-0.85}{-0.85};
        \bdot{{(0.85 + 2)/2}}{{-(0.85 + 2)/2}};\bdot{{-(0.85 + 2)/2}}{{(0.85 + 2)/2}};
    \end{tikzpicture}}  & $\longleftrightarrow$ & {        \begin{tikzpicture}[scale = 0.4,baseline={(0,0)}]
        \draw[line width=\movesline, blue] (2,2) -- (0.85, 0.85);
        \draw[line width=\movesline, blue] (2,-2) -- (0.85, -0.85);
        \draw[line width=\movesline, blue] (-2,2) -- (-0.85, 0.85);
        \draw[line width=\movesline, blue] (-2,-2) -- (-0.85, -0.85);
        \draw[line width=\movesline, blue] (0.85,-0.85) -- (-0.85, -0.85);
        \draw[line width=\movesline, blue] (0.85,0.85) -- (0.85, -0.85);
        \draw[line width=\movesline, blue] (0.85,0.85) -- (-0.85, 0.85);
        \draw[line width=\movesline, blue] (-0.85,0.85) -- (-0.85, -0.85);
        \wdot{0.85}{0.85}; \bdot{-0.85}{0.85}; \bdot{0.85}{-0.85}; \wdot{-0.85}{-0.85};
        \bdot{{(0.85 + 2)/2}}{{(0.85 + 2)/2}};\bdot{{-(0.85 + 2)/2}}{{-(0.85 + 2)/2}};
    \end{tikzpicture}} 
     & {\hspace{7ex}} & {    \begin{tikzpicture}[scale= 0.4,baseline={(0,0)}]
        \draw[line width=\movesline, blue] (-3,1) -- (-1,0);
        \draw[line width=\movesline, blue] (-3,-1) -- (-1,0);
        \draw[line width=\movesline, blue] (-1,0) -- (1,0);
        \draw[line width=\movesline, blue] (3,1) -- (1,0);
        \draw[line width=\movesline, blue] (3,-1) -- (1,0);
        \draw[line width=\movesline, blue] (3,0) -- (1,0);
        \bdot{-1}{0}; \wdot{0}{0}; \bdot{1}{0};
    \end{tikzpicture}}
    &$\longleftrightarrow$ & {    \begin{tikzpicture}[scale= 0.4,baseline={(0,0)}]
        \draw[line width=\movesline, blue] (-2,1) -- (0,0);
        \draw[line width=\movesline, blue] (-2,-1) -- (0,0);
        \draw[line width=\movesline, blue] (2,1) -- (0,0);
        \draw[line width=\movesline, blue] (2,-1) -- (0,0);
        \draw[line width=\movesline, blue] (2,0) -- (0,0);\bdot{0}{0};
    \end{tikzpicture}}
\end{tabular}
        \caption{Some local moves for plabic graphs.}
        \label{local-moves-figure}
    \end{figure}
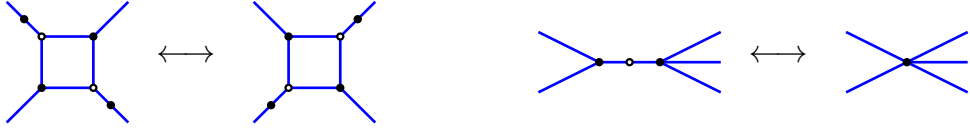

    It is important to note that while every plabic graph corresponds to only one point $X \in \grknnn$, the plabic graph representing $X$ is not unique. Instead, there is a family of plabic graphs representing the same point, and these plabic graphs are connected to one another by a series of certain local moves (see \Cref{local-moves-figure} for a subset of these moves). Graphs that have the minimal number of faces in the equivalence class are called \term{reduced}.  

    \begin{lemma}[\cite{oh-postnikov-speyer}, Theorem~6.8]
        \label{reduced-number-of-faces-lemma}
        A graph $G$ corresponding to a permutation $f= f_G$ is reduced if and only if 
        $$\# \text{faces} (G) = d(f) + 1.$$
    \end{lemma}

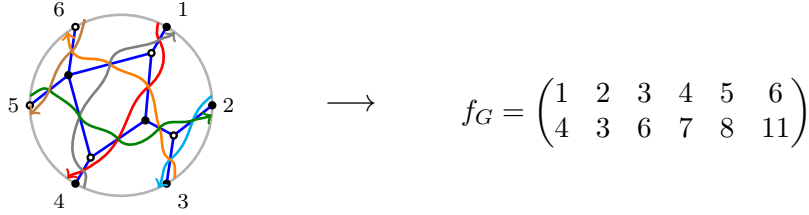
\begin{figure}
    \centering
    \begin{minipage}{0.1\textwidth}
        
    \end{minipage}
    \begin{minipage}{0.3\textwidth}
        \centering
        \newcommand{\midpoint}[4]{(0.5*(#1 + #3), 0.5*(#2 + #4))}
\newcommand{\midpointx}[4]{0.5*(#1 + #3)}
\newcommand{\midpointy}[4]{0.5*(#2 + #4)}
\newcommand{\bpointx}[1]{3*cos(60*(2-#1))}
\newcommand{\bpointy}[1]{3*sin(60*(2-#1))}
\newcommand{\bpoint}[1]{%
  ({\bpointx{#1}}, {\bpointy{#1}})
}

\newcommand{\bdott}[1]{
\draw[line width=0.3 mm,blue] [black, fill = black] #1 circle [radius = 0.1];
}
\newcommand{\wdott}[1]{
\draw[line width=0.3 mm,blue] [black, fill = white] #1 circle [radius = 0.1];
}

\newcommand{\ep}{0.1}

\begin{tikzpicture}[scale = 0.4, inner sep=0pt, outer sep=0pt]
    \tiny
    \node [inner sep=2pt, outer sep=2pt, above right] at \bpoint{1} {1};
    \node [inner sep=2pt, outer sep=2pt, right] at \bpoint{2} {2};
    \node [inner sep=2pt, outer sep=2pt, below right] at \bpoint{3} {3};
    \node [inner sep=2pt, outer sep=2pt, below left] at \bpoint{4} {4};
    \node [inner sep=2pt, outer sep=2pt, left] at \bpoint{5} {5};
    \node [inner sep=2pt, outer sep=2pt, above left] at \bpoint{6} { 6};
    
    \node (1) at \bpoint{1} {};
    \node (2) at \bpoint{2} {};
    \node (3) at \bpoint{3} {};
    \node (4) at \bpoint{4} {};
    \node (5) at \bpoint{5} {};
    \node (6) at \bpoint{6} {};

    \node (11) at ({2*cos(60)},{2*sin(60)}) {};
    \node (56) at ({-2*cos(-30)},{-2*sin(-30)}) {};
    \node (44) at ({-2*cos(60)},{-2*sin(60)}) {};
    \node (2323) at (.8,-.5) {};
    \node (23) at ({2*cos(-30)},{2*sin(-30)}) {};

    \node[inner sep=-0.5pt, outer sep=-0.5pt] (1-11) at ({2.5*cos(60)},{2.5*sin(60)}) {};
    \node[inner sep=-0.5pt, outer sep=-0.5pt] (2-23) at ($ (2) !0.5! (23) $) {};
    \node[inner sep=-0.5pt, outer sep=-0.5pt] (3-23) at ($ (3) !0.5! (23) $) {};
    \node[inner sep=-0.5pt, outer sep=-0.5pt] (4-44) at ($ (4) !0.5! (44) $) {};
    \node[inner sep=-0.5pt, outer sep=-0.5pt] (5-56) at ($ (5) !0.5! (56) $) {};
    \node[inner sep=-0.5pt, outer sep=-0.5pt] (6-56) at ($ (6) !0.5! (56) $) {};
    \node[inner sep=-0.5pt, outer sep=-0.5pt] (11-2323) at ($ (11) !0.5! (2323) $) {};
    \node[inner sep=-0.5pt, outer sep=-0.5pt] (23-2323) at ($ (23) !0.5! (2323) $) {};
    \node[inner sep=-0.5pt, outer sep=-0.5pt] (44-2323) at ($ (44) !0.5! (2323) $) {};
    \node[inner sep=-0.5pt, outer sep=-0.5pt] (44-56) at ($ (44) !0.5! (56) $) {};
    \node[inner sep=-0.5pt, outer sep=-0.5pt] (56-11) at ($ (11) !0.5! (56) $) {};

    \node[inner sep=-0.5pt, outer sep=-0.5pt] (1r) at \bpoint{1 + \ep} {};
    \node[inner sep=-0.5pt, outer sep=-0.5pt] (2r) at \bpoint{2 + \ep} {};
    \node[inner sep=-0.5pt, outer sep=-0.5pt] (3r) at \bpoint{3 + \ep} {};
    \node[inner sep=-0.5pt, outer sep=-0.5pt] (4r) at \bpoint{4 + \ep} {};
    \node[inner sep=-0.5pt, outer sep=-0.5pt] (5r) at \bpoint{5 + \ep} {};
    \node[inner sep=-0.5pt, outer sep=-0.5pt] (6r) at \bpoint{6 + \ep} {};
    \node[inner sep=-0.5pt, outer sep=-0.5pt] (7r) at \bpoint{7 + \ep} {};
    \node[inner sep=-0.5pt, outer sep=-0.5pt] (1l) at \bpoint{1 - \ep} {};
    \node[inner sep=-0.5pt, outer sep=-0.5pt] (2l) at \bpoint{2 - \ep} {};
    \node[inner sep=-0.5pt, outer sep=-0.5pt] (3l) at \bpoint{3 - \ep} {};
    \node[inner sep=-0.5pt, outer sep=-0.5pt] (4l) at \bpoint{4 - \ep} {};
    \node[inner sep=-0.5pt, outer sep=-0.5pt] (5l) at \bpoint{5 - \ep} {};
    \node[inner sep=-0.5pt, outer sep=-0.5pt] (6l) at \bpoint{6 - \ep} {};
    \node[inner sep=-0.5pt, outer sep=-0.5pt] (7l) at \bpoint{7 - \ep} {};

    \draw[line width=\graphcircle,opacity=0.3] (0,0) circle (3);
    
    \draw[line width=\graphline, blue] (1) -- (11)--(56)--(44);
    \draw[line width=\graphline, blue] (4) -- (44)--(2323)--(11);
    \draw[line width=\graphline, blue] (2323) -- (23);
    \draw[line width=\graphline, blue] (2) -- (23) -- (3);
    \draw[line width=\graphline, blue] (5) --(56) -- (6); 
    \bdott{(1)};
    \bdott{(2)};
    \bdott{(3)};
    \bdott{(4)};
    \wdott{(5)};
    \wdott{(6)};
    \wdott{(23)};
    \bdott{(2323)};
    \wdott{(44)};
    \bdott{(56)};
    \wdott{(11)};

    \draw [line width=\graphline, red, ->]
    (1r) to [out=-100, in=120] (1-11)
    to [out=-60, in=40] (11-2323)
    to [out=-140, in=60] (44-2323)
    to [out=-120, in=30] (4-44) 
    to [out=-150, in=50] (4l); 

    \draw [line width=\graphline, cyan, ->]
    (2r) to [out=220, in=70] (2-23)
    to [out=-110, in=40] (3-23)
    to [out=-140, in=60] (3l);

    \draw [line width=\graphline, orange, ->]
    (3r) to [out=90, in=-30] (3-23)
    to [out=150, in=-70] (23-2323)
    to [out=110, in=-60] (11-2323)
    to [out=120, in=-40] (56-11)
    to [out=140, in=-60] (6-56)
    to [out=120, in=-80] (6r);

    \draw [line width=\graphline, gray, ->]
    (4r) to [out=80, in=-40] (4-44)
    to [out=140, in=-110] (44-56)
    to [out=70, in=-100] (56-11)
    to [out=80, in=-150] (1-11)
    to [out=30, in=-130] (1l);

    \draw [line width=\graphline, green!50!black, ->]
    (5l) to [out=30, in=140] (5-56)
    to [out=-40, in= 120] (44-56)
    to [out=-60, in=130] (44-2323)
    to [out=-50, in=-130] (23-2323)
    to [out=50, in=-160] (2-23)
    to [out=20, in=-150] (2l);

    \draw [line width=\graphline, brown, ->]
    (6l) to [out=80, in=40] (6-56)
    to [out=-140, in=70] (5-56)
    to [out=-110, in=50] (5r);



    \normalsize
\end{tikzpicture}
    \end{minipage}
    \begin{minipage}{0.1\textwidth}
        $\longrightarrow$
    \end{minipage}
    \begin{minipage}{0.25\textwidth}
        \centering
        \small
        \begin{equation*}
                f_G = \begin{pmatrix}
                1 & 2 & 3 & 4 & 5 & 6 \\
                4 & 3 & 6 & 7 & 8 & 11
            \end{pmatrix}
            \end{equation*} 
        \normalsize
    \end{minipage}
    \begin{minipage}{0.1\textwidth}
        
    \end{minipage}
    \caption{Strand permutation for a reduced plabic graph.}
    \label{permutation-figure}
    \end{figure}

    The advantage of reduced graphs is that if a plabic graph $G$ is reduced, then the corresponding permutation $f_G$ can be also computed through \term{strands}, that start at boundary vertices, go along edges and turn maximally right (resp. left) at black (resp. white) vertices.
    In this case, for each $i \in [n]$, $f(i)$ is congruent modulo $n$ to the number of the ending vertex of the strand that starts at $i$. See \Cref{permutation-figure} for an example.

    The relationship between $\grknnn$ and its associated combinatorial objects -- plabic graphs, bounded affine permutations, and positroid cells -- is captured by the following theorem of Postnikov.

    \begin{theorem}[\cite{postnikov}]
    \label{postnikov-theorem}
    \begin{enumerate}
        \item \( \grknnn = \bigsqcup\limits_{f \in \Bkn} \cell{f} \) is a union of \term{positroid cells} $\cell{f}$ labeled by bounded affine permutations $f$.
        \item Each positroid cell \( \cell{f} \) is homeomorphic to \( \mathbb{R}^{k(n-k) - \ell(f)} \) and can be parametrized by a reduced plabic graph.
        \item For two permutations $f$ and $g$, 
        $$\cell{f} \subset \overline{\cell{g}} \quad \text{if and only if } \quad f \leq g \text{ in the affine Bruhat order}.$$
        That is,
        \[
        \overline{\cell{g}} = \bigsqcup_{f \leq g} \cell{f}.
        \]
    \end{enumerate}
    \end{theorem}

\subsection{Adding bridges}
\label{subsection: background-adding-bridges}

\begin{figure}
    \centering
    \begin{center}

\begin{tabular}{ccccc}

    {
\begin{tikzpicture}[baseline={(0,0)},scale = \bridgescale]

    \draw [line width=\bridgecircle, opacity=0.3,domain=-60:60] plot ({3*cos(\x)}, {3*sin(\x)});
    
    \node [above right] at ({3*cos(30)},{3*sin(30)}) {$i$};
    \node [below right] at ({3*cos(-30)},{3*sin(-30)}) {$i+1$};
    
    \draw[line width=\bridgeline, blue] ({3*cos(30)},{3*sin(30)}) -- ({0.6},{3*sin(30)});
    \draw[line width=\bridgeline, blue] ({3*cos(30)},{3*sin(-30)}) -- ({0.6},{3*sin(-30)});
    
    \draw[line width=\bridgeline, blue] ({0.6},{3*sin(30)}) -- ({0},{3*sin(30) + 0.3});
    \draw[line width=\bridgeline, blue] ({0.6},{3*sin(30)}) -- ({0},{3*sin(30) - 0.3});
    
    \draw[line width=\bridgeline, blue] ({0.6},{-3*sin(30)}) -- ({0},{-3*sin(30) + 0.4});
    \draw[line width=\bridgeline, blue] ({0.6},{-3*sin(30)}) -- ({0},{-3*sin(30)});
    \draw[line width=\bridgeline, blue] ({0.6},{-3*sin(30)}) -- ({0},{-3*sin(30) - 0.4});

    \wdot{{0.6}}{{3*sin(30)}};
    \wdot{{0.6}}{{-3*sin(30)}};
    \bdot{{3*cos(30)}}{{3*sin(30)}}; 
    \bdot{{3*cos(-30)}}{{3*sin(-30)}}; 

    
\end{tikzpicture}

    }
    &
    {$\rightarrow$ \quad}
    &
    {
\begin{tikzpicture}[baseline={(0,0)}, scale = \bridgescale]

\draw [line width=\bridgecircle,opacity=0.3,domain=-60:60] plot ({3*cos(\x)}, {3*sin(\x)});

\node [above right] at ({3*cos(30)},{3*sin(30)}) {$i$};
\node [below right] at ({3*cos(-30)},{3*sin(-30)}) {$i+1$};

\draw[line width=\bridgeline, blue] ({3*cos(30)},{3*sin(30)}) -- ({0.6},{3*sin(30)});
\draw[line width=\bridgeline, blue] ({3*cos(30)},{3*sin(-30)}) -- ({0.6},{3*sin(-30)});

\draw[line width=\bridgeline, red] ({1.5*cos(30) + 0.6/2},{3*sin(30)}) -- ({1.5*cos(30) + 0.6/2},{-3*sin(30)});

\draw[line width=\bridgeline, blue] ({0.6},{3*sin(30)}) -- ({0},{3*sin(30) + 0.3});
\draw[line width=\bridgeline, blue] ({0.6},{3*sin(30)}) -- ({0},{3*sin(30) - 0.3});

\draw[line width=\bridgeline, blue] ({0.6},{-3*sin(30)}) -- ({0},{-3*sin(30) + 0.4});
\draw[line width=\bridgeline, blue] ({0.6},{-3*sin(30)}) -- ({0},{-3*sin(30)});
\draw[line width=\bridgeline, blue] ({0.6},{-3*sin(30)}) -- ({0},{-3*sin(30) - 0.4});

\wdot{{1.5*cos(30) + 0.6/2}}{{3*sin(30)}};
\bdot{{1.5*cos(30) + 0.6/2}}{{-3*sin(30)}};

\wdot{{0.6}}{{3*sin(30)}};
\wdot{{0.6}}{{-3*sin(30)}};
\bdot{{3*cos(30)}}{{3*sin(30)}}; 
\bdot{{3*cos(-30)}}{{3*sin(-30)}}; 

\end{tikzpicture}

}
    &
    {$\rightarrow$ \quad}
    &
    {
\begin{tikzpicture}[baseline={(0,0)}, scale = \bridgescale]

\draw [line width=\bridgecircle, opacity=0.3,domain=-60:60] plot ({3*cos(\x)}, {3*sin(\x)});

\node [above right] at ({3*cos(30)},{3*sin(30)}) {$i$};
\node [below right] at ({3*cos(-30)},{3*sin(-30)}) {$i+1$};

\draw[line width=\bridgeline, blue] ({3*cos(30)},{3*sin(30)}) -- ({0.6},{3*sin(30)});
\draw[line width=\bridgeline, blue] ({3*cos(30)},{3*sin(-30)}) -- ({0.6},{3*sin(-30)});

\draw[line width=\bridgeline, red] ({1.5*cos(30) + 0.6/2},{3*sin(30)}) -- ({1.5*cos(30) + 0.6/2},{-3*sin(30)});

\draw[line width=\bridgeline, blue] ({0.6},{3*sin(30)}) -- ({0},{3*sin(30) + 0.3});
\draw[line width=\bridgeline, blue] ({0.6},{3*sin(30)}) -- ({0},{3*sin(30) - 0.3});

\draw[line width=\bridgeline, blue] ({0.6},{-3*sin(30)}) -- ({0},{-3*sin(30) + 0.4});
\draw[line width=\bridgeline, blue] ({0.6},{-3*sin(30)}) -- ({0},{-3*sin(30)});
\draw[line width=\bridgeline, blue] ({0.6},{-3*sin(30)}) -- ({0},{-3*sin(30) - 0.4});

\wdot{{1.5*cos(30) + 0.6/2}}{{3*sin(30)}};
\bdot{{1.5*cos(30) + 0.6/2}}{{-3*sin(30)}};

\bdot{{0.75*cos(30) + 3*0.6/4}}{{3*sin(30)}};

\wdot{{0.6}}{{3*sin(30)}};
\wdot{{0.6}}{{-3*sin(30)}};
\bdot{{3*cos(30)}}{{3*sin(30)}}; 
\wdot{{3*cos(-30)}}{{3*sin(-30)}}; 


\end{tikzpicture}

}
\end{tabular}

\end{center}
    \caption{The process of adding a bridge between $i$ and $i+1$, white at $i$. We add a new edge on the first step and then adjust the colors to make sure the resulting graph is bipartite on the second step.}
    \label{adding-bridges-figure}
\end{figure}
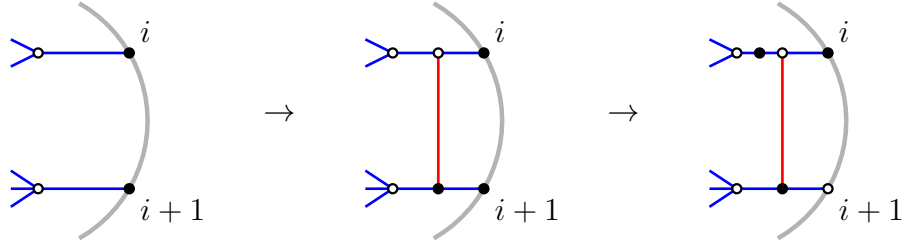

The main tool that we use for constructing plabic graphs is \textbf{bridge addition}. For a plabic graph $G$, bridge addition is the operation of attaching a bicolored edge between two boundary vertices $i$ and $i+1$ to create a new graph $G'$, as demonstrated in \Cref{adding-bridges-figure}. The bridge can be either white at $i$ or black at $i$. See \cite[section 7.4]{lam-textbook} for a full description.

To describe the process of adding the bridge in matrix language, take a point $X$ corresponding to a weighted plabic graph $G$. Then if $G'$ is obtained from $G$ by adding a bridge between $i$ and $i+1$, white at $i$, of weight $a$, then $Y$, the point corresponding to $G'$, is related to $X$ by
\begin{align}
    \label{formula-adding-bridge-x}
    Y = X \cdot x_i(a),
\end{align}
where $x_i(a)$ is the elementary matrix with $1$'s on the diagonal and $a$ in position $(i,i+1)$ for $i \in [n-1]$, and $x_n(a)$ is the same except with $(-1)^{\,n-1}a$ in position $(n,1)$.

Similarly, if $Y$ is obtained from $X$ by adding a bridge between $i$ and $i+1$, black at $i$, of weight $b$, then 
\begin{align}
    \label{formula-adding-bridge-y}
    Y = X \cdot y_i(b),
\end{align}
where $y_i(b)$ is the elementary matrix with $1$'s on the diagonal and $b$ in position $(i+1,i)$ for $i \in [n-1]$, and $y_n(b)$ is the same except with $(-1)^{\,n-1}b$ in position $(1,n)$. For more details on this construction, see \cite[Section~7.4]{lam-textbook}.

\begin{lemma}[\cite{lam-textbook}, Section 7.6]
    \label{add-bridges-permutations-background}

    Suppose $X\in \grknnn$ corresponds to the permutation $f$. 

    Then for every $a> 0$, the point $X \cdot x_i(a)$ corresponds to the permutation $f \star s_i$, and $X \cdot y_i(a)$ corresponds to $s_i \star f.$
 
\end{lemma}

To describe how adding a bridge changes the Pl\"ucker coordinates of a point, we introduce new notation.
\begin{definition}
    \label{alp-beti-definition}
    Define two operations $\alp$ and $\beti$ that acts on index sets in the following way:
    $$\alp(I) = \begin{cases}
        I - \{i+1\} \cup \{i\}, i \not\in I, i+1 \in I\\
        \varnothing, \text{ otherwise}
    \end{cases}.$$
    and
    $$\beti(I) = \begin{cases}
        I - \{i\} \cup \{i+1\}, i \in I, i+1 \not\in I\\
        \varnothing, \text{ otherwise}
    \end{cases}.$$

    To make this definition make sense, we say that $\Delta_{\varnothing}(X) = 0$ for every $X$.
\end{definition}
\begin{lemma}[\cite{lam-textbook}, Lemma 7.6]
    \label{add-bridges-minors-background}
    Suppose $Y \in \grkn$ is obtained from $X \in \grkn$ by adding a bridge between $i$ and $i+1$, white at $i$, of weight $a$ (i.e. $Y = X \cdot x_i(a))$. Then 
    $$\Delta_I(Y) = \Delta_I(X) + a \Delta_{\alp I}(X).$$
    Similarly, if $Y$ is obtained from $X$ by adding a bridge between $i$ and $i+1$, black at $i$, of weight $a$ (i.e. $Y = X \cdot y_i(a)$), then 
    $$\Delta_I(Y) = \Delta_I(X) + a \Delta_{\beti I}(X).$$
\end{lemma}

We can also delete bridges.
\begin{lemma}[\cite{lam-textbook}, Section 7.6]
    \label{delete-bridges}
    Suppose $X\in \grknnn$ is a point corresponding to a permutation $f$ with $f' = f  s_i < f$. Denote $c = \frac{\Delta_J(X)}{\Delta_{\alpha J} (X)} > 0$, where $J = \Imin{j}(f)$. 

    Let $Z \in \grkn$ be obtained from $X$ by adding a bridge between $i$ and $i+1$, white at $i$, of weight $-a$ (for some $a\geq 0$), i.e. $Z = X \cdot x_i(-a)$. Then:
    \begin{itemize}
        \item If $a < c$, then $Z \in \cell{f}$.
        \item If $a = c$, then $Z \in \cell{f'}$.
        \item If $a > c$, then $Z \not\in \grknnn$. 
    \end{itemize}
\end{lemma}

\subsection{Adding lollipops}
\label{subsection: background-adding-lollipops}

Another basic local operation on plabic graphs is the addition of a \term{lollipop}.  
Given a plabic graph $G$ with $n-1$ boundary vertices, we may insert a new boundary vertex between vertices $i$ and $i+1$, then renumber so that the new vertex receives label $i$.  
A leaf is attached at this vertex, colored either black or white.  
This operation corresponds to adding a fixed point to the permutation associated with $G$: if the lollipop is black, the new fixed point satisfies $i \mapsto i$, and if it is white, it satisfies $i \mapsto i+n$.

\subsection{Grassmannian necklace}
\begin{definition}
    For $X \in \grknnn$, define the matroid $\M_X = \{I\in {{[n] \choose k}} \st \Delta_I(X) > 0\}$.

    Similarly, for every $f \in \Bkn$, define $\M_f = \M_X$, where $X \in \cell{f}$. (It is easy to see that $\M_f$ does not depend on the choice of $X$.)
\end{definition}

\begin{definition}
    For each $i$, let $<_i$ be the order on $[n]$ given by $i <_i i+1 <_i \ldots <_i n <_i 1 <_i \ldots <_i i-1$. Similarly for $I, J \in {{[n]} \choose k}$, say that $I <_i J$ if $i_\ell <_i j_{\ell}$ for each $\ell$, where $I = \{i_1 <_i \ldots <_i i_k\}$ and $J = \{j_1 <_i \ldots <_i j_k\}$.

    It turns out that for each $X \in \grknnn$, $\M_X$ has a unique $<_i$-minimal element, denoted $\Imin{i}(X)$, and a unique maximal element, denoted $\Imax{i}(X)$.
\end{definition}

\begin{remark}
    Suppose $X \in \grknnn$ and $f = f_X$. Then $\Imin{i}(X)$ and $\Imax{i}(X)$ are uniquely determined by $f$ (since $\M_X$ is determined by $f$), so we may use the notation $\Imin{i}(f) = \Imin{i}(X)$ and $\Imax{i}(f) = \Imax{i}(X)$.
\end{remark}

\begin{lemma}[\cite{pashas}, Remark 4.3]
    \label{lemma: f-from-grassmann-necklace}
    If $f = f_X$, then for each $i$, $f(i)$ is determined by the following:
    \begin{itemize}
        \item If $i$ is not a fixed point, then $\Imin{i+1}(f) \setminus \Imin{i}(f) = \{f(i) \mod n\}$.

        \item If $f(i) = i$, then $\Imin{i+1}(f) = \Imin{i}(f)$ and $i \not\in \Imin{i}(X)$.

        \item If $f(i) = i+n$, then $\Imin{i+1}(f) = \Imin{i}(f)$ and $i \in \Imin{i}(X)$.
    \end{itemize}

    Similarly, $f^{-1}(i)$ is determined by the following:
    \begin{itemize}
        \item If $i$ is not a fixed point, then $\Imax{i}(f) \setminus \Imax{i+1}(f) = \{f^{-1}(i) \mod n\}$.

        \item If $f^{-1}(i) = i$, then $\Imax{i+1}(f) = \Imax{i}(f)$ and $i \not\in \Imax{i}(f)$.

        \item If $f(i) = i+n$, then $\Imax{i+1}(f) = \Imax{i}(f)$ and $i \in \Imax{i}(f)$.
    \end{itemize}
\end{lemma}

We can also define $\Imin{i}$ directly from $f$.
\begin{lemma}[\cite{lam-textbook}, Section 6.3]
    \label{Imin-def-through-f-lemma}
    Suppose $X \in \grknnn$ and $f = f_X$. Then for each $i \in [n],$ 
    $$\Imin{i}(f) = \{f(j)\st j< i \text{ and } f(j) \geq i\} \mod n.$$
\end{lemma}

$\Imax{i}(f)$ can be defined through a similar formula.

\begin{corollary}
    \label{Imin-gets-smaller-cor}
    Suppose $f \in \Bkn$ and $fs_i > f$. Then 
    $$\Imin{i+1}(fs_i) = \Imin{i+1}(f) \setminus \{f(i)\} \cup \{f(i+1)\} <_i \Imin{i+1}(f)$$

    Similarly, for $j \neq i$, if $s_j f > f$ for some $j$, then 
    $$\Imin{i+1}(s_j f) \leq_i \Imin{i+1}(f).$$
\end{corollary}
\begin{proof}
    First, by \Cref{condition-fsi-greaterthan-f-lemma}, we get
    $$i < i+1 \leq f(i+1) < f(i) \leq i + n.$$
    Now, comparing the expressions of $\Imin{i+1}(f)$ and $\Imin{i+1}(s_if)$, we get  
    $$\Imin{i+1}(fs_i) = \Imin{i+1}(f) \setminus \{f(i)\} \cup \{f(i+1)\}.$$
    And since $f(i+1) < f(i)$, we get 
    $$\Imin{i+1}(fs_i) <_i \Imin{i+1}(f).$$
    The second part follows similarly. 
\end{proof}

Throughout this paper, we use the notation ${{[n]} \choose k} = \{I \subset [n]\st |I| = k\}$.

\begin{lemma}[\cite{lam-textbook}, Section 6.4]
    \label{lemma-small-I-zero-minor}
    Suppose $J \in {n \choose k}$ and $J <_i \Imin{i}(f)$. Then $\Delta_J(X) = 0$ for every $X \in \cell{f}$.
\end{lemma}

\begin{lemma}
    \label{bruhat-order-matroids}
    Suppose $f, g \in \Bkn$. Then $f \leq g$ if and only if $\M_f \subset \M_g$.
\end{lemma}
\begin{proof}
    First, suppose $f \leq g$ and show that $\M_f \subset \M_g$. Notice that by \Cref{postnikov-theorem}, since $f \leq g$, we have $$\cell{f} \subset \overline{\cell{g}}.$$
    On the other hand,
    $$\cell{f} = \{X\st \M_X = \M_f\}$$
    and $\cell{g} = \{X\st \M_X=  \M_g\},$ so 
    $$\overline{\cell{g}} \subset \{X\st \M_X \subset \M_g\}$$
    Thus, we get 
    $$\{X\st \M_X = \M_f\} \subset \{X\st \M_X \subset \M_g\},$$
    which implies that $\M_f \subset \M_g$.
    
    Now, suppose $\M_f \subset \M_g$. We wish to show that $f \leq g$. Since $\M_f \subset \M_g$, we get 
    $$\Imin{i}(f) \leq \Imin{i}(g) \text{ for all } i \in [n].$$
    By \cite{lam-textbook}, Theorem 6.2, this implies that $f \leq g$.

\end{proof}


\section{The totally nonnegative Lagrangian Grassmannian} \label{lagrangian-section}


Recall that we defined the bilinear form $\form$ (\Cref{Q-def-equation}), that corresponds to the matrix
    $$\form = \begin{pmatrix}
         &\begin{matrix}
        1 & &  \\
        & -1 & \\
        & &\ddots  
      \end{matrix}\\
      
      \begin{matrix}
        -1 & &  \\
        & 1 & \\
        & &\ddots  
      \end{matrix} & 
    \end{pmatrix}$$
and the associated \term{Lagrangian Grassmannian}
    $$\lgr = \{X \in \gr| \form(x,y) = 0 \text{ for all } x,y \in X\}.$$
From that, we defined the \term{totally nonnegative Lagrangian Grassmannian} 
    $$\lgnn = \lgr \cap \grnn.$$

\subsection{A Pl\"ucker description of $\lgr$}
We show that the Lagrangian Grassmannian can be also defined through a certain condition on minors.

First, we define several operations.

\begin{definition}
    Define the operator 
    $$\rot: [2n] \to [2n]: i \mapsto i + n \mod 2n$$
    that rotates the $2n$ vertices placed on a circle by $180^\circ$. This map induces 
    $$\rot:{[2n] \choose n} \to {[2n]\choose n}: I \mapsto \rot(I)$$ 
    that rotates index sets. We also define
    $$\mainop: I \mapsto (\rot(I))^C = \rot (I^C),$$
    where $I^C = [2n] \setminus I$.
\end{definition}

\begin{definition}
    \label{def: graph-mainop-rot-swap}
    For a plabic graph $G$, define $\rot(G)$ to be the $180^\circ$-rotation of $G$, where we relabel each vertex $i$ to $i+n$.
     Define also $\swap(G)$ to be the color swap of $G$, i.e. the plabic graph with the same edges and the opposite vertex colors. Finally, let $$\mainop(G) = \swap(\rot(G)) = \rot(\swap(G)).$$
    
    Similarly define $\swap(N)$, $\rot(N)$ and $\mainop(N)$ for a plabic network $N$.

\end{definition}

Finally, we define similar operations for a point in $\grnn$.

\begin{definition}
    Suppose $X \in \grnn$. Define
    $$\rho(X) = \{\rho(x)\st x \in X\},$$
    where $\rho(x_1, \ldots, x_{2n}) = ((-1)^{n-1}x_{n+1}, \ldots, (-1)^{n-1}x_{2n}, x_1, \ldots, x_n).$
    Also, let
    $$\swap(X) := \alt(X^\perp),$$
    where $\alt: (x_1, \ldots, x_{2n}) \mapsto (-x_1, x_2, -x_3, \ldots, x_{2n})$. 

    We let $\mainop(X) = \rot(\swap(X)).$
    
\end{definition}

It is then easy to see that in each of the settings $\rot$ and $\swap$ (or the complement map) are involutions that commute with each other, so $\mainop$ is an involution as well.

The repeated notations can be justified with the following lemma. 

\begin{lemma}
    \label{lemma:swap-on-I-equals-swap-on-X}
    Suppose $X = X(N) \in \grnn$. Then 
    $$\swap(X(N)) = X(\swap(N)) \text{ and } \rot(X(N)) = X(\rot(N)).$$
    Moreover, there exist constants $a,b \neq 0$ such that for any $I \in {{[2n]}\choose n}$, we have 
    $$\Delta_{I^C}(X) = a \Delta_I(\swap(X)) \text{ and } \Delta_{\rot(I)}(X) = b \Delta_I(\rot(X)). $$
\end{lemma}

\begin{remark}
    Here, we need constants $a$ and $b$ because Pl\"ucker coordinates of each point are defined up to common rescaling.
\end{remark}

\begin{proof}
    Throughout this proof, let $I \in \twonchoosen$, $X \in \grnn$ and $N$ be a plabic network.

    First, notice that it follows from the definition of $\rot(X)$ and $\rot(N)$ that 
    $$\rot(X(N)) = X(\rot(N))$$
    Similarly, we get
    $$\Delta_{\rot(I)}(N) = \Delta_I(\rot(N)).$$
    Combining, we obtain
    $$\Delta_{\rot(I)} (X) = b\Delta_I(\rot(X))$$
    for some $b \neq 0$. (We just need to take $N$ to be a plabic network corresponding to $X$.) Since we get $\Delta_I(X(\rot(N))) = b \Delta_I(\rot(X(N))$, this means that
    $$\rot(X(N)) = X(\rot(N)).$$
    
    Now, notice that if $I^C$ is the boundary of an almost perfect matching in $N$, then $I$ is the boundary of the corresponding matching in $\swap(N)$. Therefore, we get
    $$\Delta_{I^C}(N) = \Delta_I(\swap(N)).$$
    Finally, it is a classical result (for example, see \cite{hochster1975topics}) that 
    $$\Delta_{I^C}(X) = a\Delta_I(\alt(X^\perp)) = a\Delta_I(\swap(X))$$
    for some $a \neq 0$.
    Combining with the equality above, we get 
    $$\swap(X) = X(\swap(N)).$$
\end{proof}

\begin{definition}
    Call a point $X$ in the Grassmannian \term{\tsym} if its Plu\"cker coordinates are invariant under $\mainop$, i.e.
    $$\Delta_I(X) = \Delta_{\mainop(I)} (X) \text{ for all } I.$$
\end{definition}

\begin{remark}
    \label{remark: tsym-equals-R-invar}
    From \Cref{lemma:swap-on-I-equals-swap-on-X}, we get that $X$ is \tsym if and only if $X = \mainop(X)$.
\end{remark}

\begin{proposition}
    \label{minor-condition-on-lg}
    The Lagrangian Grassmannian is the space of \tsym points. That is, 
    \begin{align}
        \label{minor-condition-equation}
        \lgr = \{X \in \gr\st \Delta_I(X) = \Delta_{\mainop(I)} (X) \text{ for all } I \in {{[2n] \choose n}} \}.
    \end{align}
    Similarly, 
    \begin{align}
        \label{tnn-minor-condition-equation}
        \lgnn = \{X \in \gr\st \Delta_I(X) = \Delta_{\mainop(I)} (X) \geq 0 \text{ for all } I \in {{[2n] \choose n}} \}.
    \end{align}
    
\end{proposition}

\begin{proof}
    Following a similar proof in \cite{pashas}, it can be computed that for every $x, y \in \RR^{2n}$, we have
    $$\form(x,y) = (-1)^n (x, \alt(\rot(y))).$$
    
    Therefore for $X \in \grnn$, we have that $X \in \lgnn$ if and only if for every $x, y \in X$
    $$(x, \alt(\rot(y)))= 0 .$$
    In turn, since $\dim(X) = \dim(X^\perp) = n$, this is equivalent to 
    $$X^\perp = \alt(\rot(X)).$$
    This is the same as
    $$X = \rot(\alt(X^\perp)) = \rot(\swap(X)) = \mainop(X).$$
    Therefore, we get that 
    \begin{align}
        X \in \lgr \iff X = \mainop(X).
        \label{eq: lgr-X-equals-RX}
    \end{align}
    We can now prove the first part of the statement in two parts.

    On the one hand, suppose $\Delta_I(X) = \Delta_{\mainop(I)}(X)$ for all $I$. Then by \Cref{lemma:swap-on-I-equals-swap-on-X}, there exists $c = ab \neq 0$ such that 
    $$\Delta_{\mainop(I)} (X) = c \Delta_I(\mainop(X)).$$
    Therefore, since $\Delta_I(X) = \Delta_{\mainop(I)}(X)$, we get 
    $$\Delta_I(X) = c \Delta_I{\mainop(X)} \text{ for all } I.$$
    And since Pl\"ucker coordinates determine a point in $\grnn$, we get 
    $$X = \mainop(X).$$
    Therefore by \Cref{eq: lgr-X-equals-RX}, we get $X \in \lgr$. 

    On the other hand, suppose $X \in \lgr$, i.e. $X = \mainop(X)$. Applying \Cref{lemma:swap-on-I-equals-swap-on-X}, we get that there exists $c \neq 0$ such that 
    $$\Delta_I(X) =  c \Delta_{\mainop(I)}(X)= c^2 \Delta_{\mainop^2 (I)}(X) = c^2 \Delta_I(X).$$
    Therefore $c = \pm 1$. We will show in \Cref{minus-not-possible} below that the case $c  = -1$ is impossible, which implies that $c = 1$.

\end{proof}
\begin{example}
    To demonstrate \Cref{minor-condition-on-lg} for $n = 2$, notice that for $I = \{2,4\}$, we get
    $$\mainop(I) = \rot(I^C) = \rot(\{1,3\}) = \{2,5\} \mod 4 = \{1,3\}.$$
    (Consequently, we also have $\mainop(\{1,3\}) = \{2,4\}$. ) For the other index sets $I$, we get $\mainop(I) = I$. Therefore 
    $$\lgrr(2,4) = \{X \in \grrrr(2,4) \st \Delta_{13}(X) = \Delta_{24}(X)\}.$$
\end{example}

\begin{lemma}
    \label{minus-not-possible}
    There does not exist $X \in \grnn$ such that 
    $$\Delta_I(X) = -\Delta_{\mainop(I)} (X) \text{ for all } I.$$
\end{lemma}
\begin{proof}
We will first study $\mainop(I)$. For $i \in I$, say that $i \in \pairs(I)$ if $n+i \in I$. Otherwise say that $i \in \sing(I)$. Notice that if $i \in \sing(I)$, then $i \in \mainop(I)$ and $i\in \sing(\mainop(I))$. On the other hand, if $i \in \pairs(I)$, then $i, n + i \not\in \mainop(I)$. That is, $S$ preserves the single terms of $I$, but changes the terms that come in pairs.

    Now, suppose that $X$ satisfies $\Delta_I(X) = - \Delta_{\mainop(I)}(X)$ for all $I$. Show that this means that $\Delta_I(X) = 0$ for all $I$ by induction on $|\pairs(I)|$. 

    First, for the base of induction, suppose $\pairs(I) = \varnothing$. That means that $I = \sing(I)$, so $\mainop(I) = I$. But then $\Delta_I(X) = - \Delta_I(X)$, so $\Delta_I(X) = 0$.

    Next, suppose we have $\Delta_I(X) = 0$ for all $I$ with $|\pairs(I)| < 2k$ and fix $I$ with $|\pairs(I)| = 2k$. Notice that this also means that $|\pairs(\mainop(I))| = 2k$. Choose $j \in \pairs(\mainop(I))$. Then $j, n+j \not\in I$. 
    
    Now use the Pl\"ucker relation (see \cite{lam-textbook}, Proposition 3.1) to write 
    $$\Delta_I(X) \Delta_{\mainop(I)}(X) = - \Delta_I^2(X)= \sum_{i \in I} \pm \Delta_{I \setminus i \cup j}(X) \Delta_{\mainop(I) \setminus j \cup i} (X)$$
    (i.e. we replace $j \in \mainop(I)$ with all terms $i \in I$).

    Consider all terms in the right hand side of the equality separately. If $i \in \sing(I)$, then $i \in \mainop(I)$, so it will appear twice in $\mainop(I) \setminus j \cup i$, so $\Delta_{\mainop(I) \setminus j \cup i}(X) = 0$. On the other hand, if $i \in  \pairs(I)$, then $\pairs(I \setminus i \cup j) = \pairs(I) \setminus\{i,\plusn{i}\}$. In that case $|\pairs(I \setminus i \cup j)| < 2k$, so by induction hypothesis $\Delta_{I \setminus i \cup j} (X)= 0$. Therefore the right hand side of the equality is equal to $0$, which means that 
    $$\Delta_I(X) = 0.$$
    But that would mean that $X$ is the zero space, which is not in $\gr$. This shows that there is no $X$ that satisfies the given condition!
\end{proof}

\subsection{The first Schubert cell of $\lgnn$}

\begin{definition}
    For an $n \times n$ matrix $M$, define 
    $$\embed(M) = (I \st M') \in \mat_{n \times 2n},$$
    where
    \begin{equation}
    \label{eq:M'-def}
        (M')_{ij} = (-1)^{j-1} M_{i (n+1-j)},
    \end{equation}
    i.e. $M'$ is the matrix obtained from $M$ by first flipping the order of rows and then swapping the sign of even rows.
\end{definition}

\begin{proposition}
    \label{prop: lgnn-vs-persymm-matrices}
    \begin{enumerate}
        \item $\embed(M) \in \lgr$ if and only if $M \in \antisymmat_n$, i.e. $M$ is symmetric with respect to the antidiagonal.
        \item $\embed(M) \in \grnn$ if and only if $M$ is totally nonnegative.
    \end{enumerate}
    In other words, the restriction $\embed: \antisymmat_n \cap \tnnmat_n \to \lgnn$ is an embedding with a dense image.
\end{proposition}
\begin{proof}
    Let $X = \rowspan \embed(M) = \rowspan (I \st M')$. 
    \begin{enumerate}
        \item Notice that since the form $\form$ is skew-symmetric, $X$ is in $\lgr$ if and only if $\form$ evaluates to $0$ on all pairs of rows. Since $\form(x,y) = \sum_{i=1}^n (-1)^{i-1} (x_i y_{n+i} - x_{n+i}y_i)$, we have 
        $$\form(\rm{row}_i, \rm{row}_j) = (-1)^{i-1} M'_{ji}- (-1)^{j-1} M'_{ij}.$$
        Applying \Cref{eq:M'-def}, we obtain that $X \in \lgr$ if and only if 
        $$M_{i (n+1-j)} = M_{j (n+1-i)},$$
        i.e. $M$ is symmetric with respect to the antidiagonal.
    
        \item The second part of the statement follows from Lemma 3.9 in \cite{postnikov}.
    \end{enumerate}

    Finally, for the last statement, notice that the map $\embed: \antisymmat_n \cap \tnnmat_n \to \lgnn$ is continuous, injective, and its image is the first Schur cell $\lgnn \cap \{X: \Delta_{[n]} (X) \neq 0\}$, which is dense in $\lgnn$.
\end{proof}

\begin{example}
	For $n=2$, we get that 
	$$X = \embed(\begin{pmatrix}
		a & b\\
            c & a
	\end{pmatrix}) = \rowspan \begin{pmatrix}
		 1 & 0 & c & a\\
        		0 & 1 & -a & -b
	\end{pmatrix} \in \lgrr(2,4).$$
	Moreover, $X \in \lgrnn(2,4)$ if and only if $a, b,c, a^2 - bc \geq 0$.
\end{example}


\section{\Tsym permutations} \label{symm-perm-section}

Let us study the bounded affine permutations that correspond to \tsym points in $\grnn$, i.e. points in $\lgnn$. 

\begin{definition}
    \label{def: tsym-permutations}
    Call a permutation $f \in \Bnn$ \term{\tsym} if the permutation $\rotaff^{-1} f: i \mapsto f(i) - n$ (where $\rotaff: i \mapsto i + n$) is an involution. In other words, we say that $f$ is \tsym if
    $$f = \rotaff  f^{-1}  \rotaff.$$

    We denote $\sbn$ to be the set of all \tsym permutations.
\end{definition}

\subsection{Permutations corresponding to $\lgnn$}

\begin{theorem}
    \label{permutations-are-symmetric}
    If $X \in \lgnn$, then the corresponding permutation $f = f_X$ is \term{\tsym}.
\end{theorem}
\begin{proof}
    Let $\M := \{I\st \Delta_I(X) \neq 0\}$ be the matroid that corresponds to $X$. Since $X$ is \tsym, for each $I \in \M$, we have $\Delta_{\mainop{I}}(X) = \Delta_I(X) \neq 0$ (see \Cref{minor-condition-on-lg}). Therefore for $I \in M$, we have $\rotcomp{I} \in \M$.
    
    Recall also that for $I = \Imin{i} (f) = \Imin{i}(X)$, we have $I \leq_i J$ for all $J \in M$, which means $I^C \geq_i J^C$ and $\mainop{I} = \rot(I^C) \geq_{\rot(i)} \rot(J^C) = \mainop{J}$ for all $J \in \M$. For \tsym $X$, this means that 
    \begin{equation}
        \label{eq:imin-imax-sym}
        \Imax{\rot(i)} (f) = \rot((\Imin{i}(f))^C).
    \end{equation} 

    Next, recall that by \Cref{lemma: f-from-grassmann-necklace}, we have $\{f(i)\} \equiv \Imin{i+1} (f) \setminus \Imin{i} (f)\pmod{2n}$, unless that set is empty (in which case $i$ is a fixed point, i.e. $f(i) = i$ or $i+2n$). Similarly, $\{f^{-1}(i)\} \equiv \Imax{i} (f)\setminus \Imax{i+1}(f) \pmod{2n}$ (or $i$ is a fixed point of $f^{-1}$). By \Cref{eq:imin-imax-sym} we obtain
    \begin{align*}
        \Imax{\rot(i)} \setminus \Imax{\rot(i)+1} = \rot((\Imin{i})^C) \setminus \rot((\Imin{i+1})^C) = \rot(\Imin{i+1} \setminus \Imin{i}),
    \end{align*}
    which shows that $f^{-1}(\rot(i)) \equiv \rot(f(i))$. This also demonstrates that $i$ is a fixed point of $f$ if and only if $\rot(i)$ is a fixed point of $f^{-1}$ (in that case both sides of the equality above are empty). Consider two options:
    \begin{enumerate}
        \item $i$ is not a fixed point of $f$, i.e. $f(i) \neq i, i+2n$. Then since $f \in \Bnn$, 
        $$\rot(i) - 2n \leq f^{-1}(\rot(i)) \leq \rot(i),$$
        $$\rot(i) < \rot(f(i)) < \rot(i) + 2n.$$
        Thus, the only possible option to have $f^{-1}(\rot(i)) \equiv \rot(f(i))$ is if
        $$f^{-1}(\rot(i)) = \rot(f(i)) - 2n = \rot^{-1} (f(i)).$$

        \item If $i$ is a fixed point, suppose $f(i) = i$. Then $i \notin \Imin{i}(f)$. By \Cref{eq:imin-imax-sym}, this means that
        $$i+n = \rot(i) \in \Imax{\rot(i)} (f) = \rot((\Imin{i}(f))^C).$$
        Therefore, we get $f(i+n) = i+3n$. We can similarly show that if $f(i) = i+2n$, then $f(\rot(i)) = \rot(i)$.
    \end{enumerate}
    
    Combining both cases together, we get that $f(i) = (\rot  f^{-1}  \rot)(i)$ for any $i$.
\end{proof}

\begin{definition}
    For each $f \in \Bnn$, we can define 
    $$\celllgn{f} = \cell{f} \cap \lgnn.$$
\end{definition}

\begin{corollary}
    \label{cell-decomposition-corollary}
    $\lgnn = \bigsqcup\limits_{f \in \sbn} \celllgn{f}.$
\end{corollary}

\begin{remark}
    Notice that here we only prove that every \tsym point corresponds to a \tsym permutation. Later, in \Cref{symm-graphs}, we will show that every \tsym permutation corresponds to at least one \tsym point, i.e. for each $f \in \sbn$, we have $\celllgn{f} \neq \varnothing$.
\end{remark}

\subsection{Involutions}

\label{subsection: symm-perm-involutions}
Consider the poset $\affinen$ of $(n,2n)$-affine permutations and the affine Coxeter group $\affine$ (see \Cref{subsection: background-bruhat-order}).
    \begin{definition}
        Denote
        $$\sj = \{f \in \affinen \st f = \rot f^{-1} \rot\}.$$
        Then, in particular, we have $\sbn = \sj \cap \Bnn$.
    \end{definition}
    
    \begin{definition}
        Let $\invol = \{f \in \affine \st f^{-1} = f\}$ be the set of all affine involutions.  
    \end{definition}
    By \Cref{def: tsym-permutations}, up to a cyclic shift the \tsym permutations are precisely involutions.
    \begin{lemma}
        \label{lemma: symmetric-equals-involutions}
        Suppose $f \in \affinen$. Then $f \in \sbn$ if and only if $\wf = \rho^{-1} f \in \invol$. That is, $\sj = \rho (\invol)$. 
    \end{lemma}
    
    See \cite{springer1988-involutions, springer1990-involutions} (the finite case) and in \cite{hultman, hansson-hultman-involutions, hamaker-marberg-pawlowski-involutions, marberg-involutions} (the general case) for more details on the poset $\invol$.

    Notice that  the set $\sbn$ is an order ideal in $\rotaff(\invol) = \sj$ (since $\Bnn$ is an upper ideal in the set $\affinen$). Therefore the properties of $\invol$ transfer to $\sbn$.

\subsection{\Tsym Bruhat order}

    \begin{definition}
        Denote $\sless$ to be the restriction of Bruhat order onto bounded \tsym permutations. That is, $f \sleq g$ if $f,g \in \sbn$ and $f \leq g$ in Bruhat order. 
        
    \end{definition}

    Since $\sbn$ is an upper ideal in $\rotaff(J)$, the covering relation on $\sbn$ can be translated from the covering relation on $\invol$, described in \cite{marberg-involutions} (Corollary 8.12).

    In the standard Bruhat order the rank function is defined through the number of alignments (see \Cref{eq: alignment-def}). It can be seen that in the \tsym case alignments come in pairs: if $i$ and $j$ form an alignment in $f \in \sbn$, then $f(j) + n$ and $f(i) + n$ form an alignment as well. In this case we say $(i,j)$ and $(f(j) + n, f(i) + n)$ form a \term{\tpairal}.

    \begin{definition}
        \label{def: symmell-symmdim}
        For each $f \in \sbn$, define 
        $$\symmell(f) = \# \{\text{\tpairal}\}.$$

        We also define 
        $$\symmdim(f) = {{n+1} \choose 2} - \symmell(f).$$
    \end{definition}

    \begin{remark}
        Notice that $\symmdim(f) \geq 0$ for all $f$. Moreover, $\symmdim(f) = 0$ if $f$ only has fixed points and $\symmdim(f) = {{n+1} \choose 2}$ for the longest permutation $f = \rot: i \mapsto i+n$.
    \end{remark}

    Alternatively, we can define $\symmdim(f)$ through $d(f)$.

    \begin{lemma}
        \label{d-vs-symmdim}
        Suppose $f \in \sbn$. Then 
        $$d(f) = 2 \symmdim(f) - \frac{\# \{i \in [2n] \st f(i) = i+n\}}{2}.$$
    \end{lemma}
    \begin{proof}
        Recall that we can define $\ell(f)$ as the number of alignments in $f$ and $\symmell(f)$ as the number of \tpairal. Therefore, every \tpairal is counted twice in $\ell(f)$ except for the alignments that are symmetric to themselves (call them \textbf{\tcenterals}). Then 
        \begin{align}
            \label{eq:ell-symell-symm-alignments}
            \ell(f) = 2 \symmell(f) - \#\{\text{\tcenterals}\}.
        \end{align}
        Now, notice that every $i$ and $j$ form a \tcenteral if and only if $f(i) = j+n$ and $j \neq i$. That is, every $i \in [2n]$ is in (exactly one) \tcenteral unless $f(i) = i+n$. Therefore
        \begin{align}
            \label{eq:symm-alignments-i-in}
            \#\{\text{\tcenterals}\} = \frac{2n - \#\{i \st f(i) = i+n\}}{2}.
        \end{align}
        Combining \Cref{eq:ell-symell-symm-alignments} and \Cref{eq:symm-alignments-i-in} we get
        \begin{align*}
            \ell(f) = 2 \symmell(f) - n + \frac{\#\{i \st f(i) = i+n\}}{2}.
        \end{align*}
        Finally, recall that $d(f) = n^2 - \ell(f)$ and $\symmdim(f) = {{n+1}\choose 2} - \symmell(f)$. Therefore
        \begin{align*}
            d(f) = n^2 - \ell(f) =  n^2 - 2 \symmell(f) + n - \frac{\#\{i \st f(i) = i+n\}}{2} =\\
            2({{n+1}\choose 2} - \symmell(f)) -  \frac{\#\{i \st f(i) = i+n\}}{2} = 2 \symmdim(f) -  \frac{\#\{i \st f(i) = i+n\}}{2}.
        \end{align*}
    \end{proof}

    \begin{proposition}
        $\symmdim$ is the rank function for the \tsym Bruhat order $\sleq$. 
    \end{proposition}
    \begin{proof}
        First, for $f \in \Bnn$, by \Cref{lemma: symmetric-equals-involutions}, $f \in \sbn$ if and only if $\wf = \rotaff^{-1}f \in \invol$, and the order $\sless$ is the opposite to the order in $\invol$. Moreover, since $\rotaff^{-1}(\sbn)$ forms an order ideal in $\invol$, the covering relations in $\sbn$ come from the covering relations in $\invol$. Therefore, the rank in $\sbn$ is complementary to the rank in $\invol$.

        By Exercise 35 in chapter 2 of \cite{bjorner-brenti}, the rank function in $\invol$ is 
        $$r(\wf) = \frac{\ell(\wf) + \absell(\wf)}{2},$$
        where $\absell(\wf)$ is the minimal number of transpositions $\wf$ can be represented in terms of. Since $\wf$ is an involution, we get
        $$\absell(\wf) = n - \frac{\#\{i \st \wf(i) = i\}}{2} =  n - \frac{\#\{i \st f(i) = i+n\}}{2}.$$
        Substituting this above and applying \Cref{d-vs-symmdim}, we get
        \begin{align*}
            r(\wf) = \frac{\ell(\wf) + n - \frac{\#\{i \st f(i) = i+n\}}{2}}{2} = \frac{\ell(f) + n - (2 \symmdim(f) - d(f)))}{2} =\\
            \frac{n^2 - d(f) + n - (2 \symmdim(f) - d(f)))}{2} = \frac{n^2 + n}{2} - \symmdim(f) = \symmell(f).
        \end{align*}
        This shows that $\symmell$ is the rank function of the order dual to $\sless$. Therefore, $\symmdim$ is the rank function of $\sless$.
    \end{proof}

\subsection{Semidirect product}

    We now introduce the analog of multiplying by a simple transposition
    
    \begin{definition}
        \label{def: semidirect-product}
        For each $f \in \sj$ and $i \in \mathbb{Z}$, we define 
        $$f \twist{i} =  \begin{cases}
        s_{\plusn{i}} f s_i, \quad \text{if } s_{\plusn{i}} f \neq f s_i,\\
        f s_i, \quad \text{otherwise}.
    \end{cases}$$
    \end{definition}
    Notice for every $f \in \sj$, we have $f \twist{i} \in \sj$.
        \begin{remark}
            This definition is parallel to the semidirect product defined in \cite{springer1988-involutions} and \cite{springer1990-involutions} and further studied in \cite{hamaker-marberg-pawlowski-involutions}, \cite{hultman} and \cite{hansson-hultman-involutions}. The notation we use is similar to the one in \cite{hamaker-marberg-pawlowski-involutions}.
        \end{remark}
    \begin{lemma}
        \label{hultman-cross-equals-sym-cross}
        Suppose $f \in \sbn$. Then $f \twist{i}$ and $f$ are always comparable. Moreover,
        \begin{enumerate}
            \item $f \twist{i} < f$ if and only if $f s_i < f$. In this case, we get $f \twist{i} \leq fs_i < f$.
            \item $f \twist{i} > f$ if and only if $f s_i > f$. In this case, we have $f < fs_i \leq f \twist{i}$.
        \end{enumerate}
    \end{lemma}

    \begin{proof}
        By \Cref{lemma: symmetric-equals-involutions}, if $f \in \sbn$, then $\wf = \rotaff^{-1} f \in \affine$. We can also notice that the semidirect product $f \twist{i}$ is similar to the semidirect product $\wf\hultmantwist{i}$ introduced in \cite{hultman}. The two semidirect products are related by
        $$f \twist{i} = \rotaff(\wf \hultmantwist{i}).$$
        Now, it was explained in the proof of Lemma 3.8 in \cite{hultman}, 
        $$\wf \hultmantwist{i} < \wf \text{ if and only if } \wf s_i < \wf.$$
        Applying $\rotaff$ to both sides we get exactly the first part of the statement. The second part follows similarly.
    \end{proof}

    \begin{lemma}
    \label{ftwist-is-bounded}
        Suppose $f \in \sbn$. Then 
        $$f \twist{i} \in \sbn \text{ if and only if } f(i) \neq i \text{ and }f(i+1) \neq i+1 + 2n.$$
        
    \end{lemma}

    \begin{proof}

        It is easy to see that since $f \in \Bnn$,
        \begin{equation}
            \label{eq:fsi-not-bounded-fixed-points}
            fs_i \in \Bnn \text{ if and only if } f(i) \neq i \text{ and } f(i+1) \neq i + 2n + 1.
        \end{equation}
        
        That is, we need to show that $f\twist{i} \in \Bnn$ if and only if $fs_i \in \Bnn$. We will proceed by assuming the contrary. Then there exists $f \in \sbn$ such that only one of the two statements is true.

        First, suppose $fs_i \in \Bnn,$ but $f\twist{i} \notin \Bnn$. In particular, this means that $f\twist{i} \notin \{f, fs_i\}$, so $f \twist{i} = s_{i+n} (fs_i)$. Similarly to \Cref{eq:fsi-not-bounded-fixed-points}, we get that since $fs_i \in \Bnn$, but $s_{i+n} (fs_i) \notin \Bnn$, 
        $$(fs_i)(i+n+1) = i+n+1 \text{ or } (fs_i)(i+n) = i + 3n.$$
        This is equivalent to 
        $$f(i+n+1) = i+n+1 \text{ or } f(i+n) = i + 3n.$$
        But since $f \in \sbn$, we get that $f(i+1) = i+2n+1$ or $f(i) = i$, which by \Cref{eq:fsi-not-bounded-fixed-points} means that $fs_i \not\in \Bnn$, which contradicts our assumption.

        Second, suppose $fs_i \notin \Bnn$ and $f\twist{i} \in \Bnn$. Again, since $f \in \Bnn$ and $\Bnn$ is an upper order ideal, we get $fs_i < f$. In this case, by \Cref{hultman-cross-equals-sym-cross}, we have to have $f\twist{i} \leq fs_i < f$. But then since $f\twist{i} \in \Bnn$, we get $fs_i \in \Bnn$ as well. This again leads to contradiction. 
        
        Since neither of the cases are possible, we get that $fs_i \in \Bnn$ if and only if $f \twist{i} \in \Bnn$, which concludes the proof.
    \end{proof}
    
    \subsection{The operation $\beta_i$}

    \begin{definition}  
    \label{def: symcross}
        For $f\in \sbn$ and $i \in \ZZ$, let
        $$\symcross{i}{f} = \begin{cases}
            f \twist{i}, \quad \text{if } 
            f\twist{i} > f;\\
            f \quad \text{otherwise}.
         \end{cases}$$
    \end{definition}

    \begin{remark}
         One can see that $\symcross{i}{f} \geq f$ and $\symcross{i}{f} \in \sbn$ whenever $f \in \sbn$.
    \end{remark}

    We can alternatively define $\symcross{i}{f}$ using the following.
    \begin{lemma}
        \label{symcross-through-star-lemma}
        Suppose $f \in \sbn, i \in \ZZ$. Then 
        $$\symcross{i}{f} = s_{i+n} \star (f \star s_i),$$
        where $f \star s_i = \max\{f, fs_i\}$ and $s_i \star f = \max \{f, s_i f\}$ 
    \end{lemma}
    \begin{proof}
        By \Cref{def: semidirect-product}, $\symcross{i}{f}$ falls into one of the following three cases.
        \begin{enumerate}
            \item $\symcross{i}{f} = s_{i+n}fs_i > f$. In that case, we have $s_{i+n}fs_i > fs_i > f$. Then we get $f \star s_i = \max\{f, fs_i\} =  fs_i$ and $s_{i+n} \star (f \star s_i) = \max\{fx_i, s_{i+n} f s_i\} = s_{i+n} f s_i = \symcross{i}{f}$.

            \item $\symcross{i}{f} = s_{i+n}f = fs_i > f$. In that case, we get $f \star s_i = fs_i$. But $s_{i+n} (fs_i) = f < fs_i$, so $s_{i+n} \star f \star s_i = \max\{fs_i, s_{i+n} (fs_i)\} = fs_i = \symcross{i}{f}$.

            \item $\symcross{i}{f} = f.$ This means that $f\twist{i} < f$ or $f\twist{i} \not\in \sbn$ (which by \Cref{ftwist-is-bounded} implies $f(i) = i$ or $f(i+1) = i+1+2n$). In both cases we get $fs_i < f$, so $f \star s_i = f$. Since $fs_i < f$, we can also apply \Cref{condition-fsi-greaterthan-f-lemma} to get 
            $$f(i) < f(i+1).$$
            Since $f \in \sbn$, we get $f(i) = f^{-1}(i+n) + n$ and $f(i+1) = f^{-1}(i+n+1) + n$, so we get
            $$f^{-1}(i+n) \leq f^{-1}(i+n + 1).$$
            From here, we can see that $s_{i+n} f < f$. Therefore $s_{i+n} \star (f \star s_i) = s_{i+n} \star f = f = \symcross{i}{f}.$
        \end{enumerate}
        We can now see that the equality holds in all three cases.
    \end{proof}

    \begin{lemma}
    \label{lemma: conditions-symmcross-i-f-greater-f}
        \begin{enumerate}
            \item For $f \in \sbn$, we have $f < \symcross{i}{f}$ if and only if $f < fs_i$.
            \item For $g \in \sbn$, there exists $f \in \sbn$ such that $g = \symcross{i}{f} > f$ if and only if $gs_i < g$ and $g(i) \neq i, g(i+1) \neq i+1+2n$ (or, equivalently, $gs_i < g$ and $gs_i \in \Bnn$).
        \end{enumerate}
    \end{lemma}
    \begin{proof}
        \begin{enumerate}
            \item By \Cref{def: symcross}, we get that $\symcross{i}{f} > f$ if and only if $f \twist{i} > f$ (since $\Bnn$ is an upper order ideal, we get that $f \twist{i} \in \Bnn$ automatically). And it was shown in the proof of \Cref{hultman-cross-equals-sym-cross} that $f \twist{i} > f$ if and only if $fs_i > f$.

            \item Notice that $g = \symcross{i}{f} < g$ if and only if $g = f \twist{i}$ (or, equivalently, $f = g \twist{i}$) and $f < g$ (or equivalently, $f \twist{i} < f$. Because of \Cref{hultman-cross-equals-sym-cross}, this is equivalent to $gs_i < g$.
        \end{enumerate}
    \end{proof}

    From the definition of $\beta_i$ and the definition of $\symmell$ and $\symmdim$, we get:
    
    \begin{lemma}
        If $g = \symcross{i}{f} > f$, then $\symmdim(g) = \symmdim(f) + 1$.
    \end{lemma}
    \begin{proof}

        Notice that for $j \notin \{i,i+1\}$, we have that $i$ and $j$ form a \tpairals in $f$ if and only $i+1$ and $j$ form a \tpairals in $\symcross{i}{f}$. Similarly, $i+1$ and $j$ form a \tpairals in $f$ if and only $i$ and $j$ form a \tpairals in $\symcross{i}{f}$. Finally, we can see that $i$ and $i+1$ form a \tpairals in $\symcross{i}{f}$, but not in $f$. We can then conclude by \Cref{def: symmell-symmdim}.
    \end{proof}

    We now add a lemma that allows us to go down in order.

    \begin{lemma}
        \label{going-down-bruhat}
        Suppose $f \in \sbn$. Then either $f$ has a pair of \tsym fixed points (that is, there exists $i$ such that $f(i) = i$ and $f(i+n) = i+3n$) or there exists $f' \in \sbn$ such that $f = \symcross{i}{f'} > f'$ for some $i$.
    \end{lemma}
    \begin{proof}

        First, notice that since $f \in \sbn$, $f(i) = i$ if and only if $f(i+n) = i+3n$. That is, if $f$ has a fixed point $i$, then $f$ automatically has a pair of symmetric points, so we may assume that $f$ is free of fixed points.

        Now, since $f$ cannot be non-increasing, there exists $i$ such that $f(i) < f(i+1)$. Since neither of $i$ and $i+1$ are fixed points, we get
        $$i < i + 1 \leq f(i) < f(i + 1) \leq i + n.$$
        By \Cref{condition-fsi-greaterthan-f-lemma}, we know that this is equivalent to $fs_i < f$. Moreover, since neither of $i$ and $i+1$ are fixed points, by \Cref{lemma: conditions-symmcross-i-f-greater-f}, we get that there exists $f' \in \sbn$ such that $f = \symcross{i}{f'} > f'$.
    \end{proof}


\section{\Tsym graphs} \label{symm-graphs-section}


Recall from \Cref{def: graph-mainop-rot-swap} that for a graph $G$, $\mainop(G) = \rot(\swap(G))$ is the $180^\circ$ rotation of the color swap of $G$. We call $G$ \term{\tsym} if $G = \mainop(G)$.

In this section we will prove one of the main theorems.

\begin{theorem}
    \label{symm-graphs}
    The cell $\celllgn{f}$ is nonempty if and only if $f \in \sbn$. 
    
    For $f \in \sbn$, 
        the cell $\celllgn{f}$ is homeomorphic to $\RR^{{{n+1} \choose 2}-\symmell(f)} = \RR^{\symmdim(f)},$ 
        and can be parametrized by a \tsym, non-reduced, plabic graph.
\end{theorem}

Here, we use the following definition
\begin{definition}
    Suppose $G$ is a \tsym plabic graph. We say that $G$ \textbf{parametrizes} the cell $\celllgn{f}$ if for every $X \in \celllgn{f}$, we can assign symmetric positive weights to edges of $G$ (i.e. edges symmetric to each other have the same weight) to create a plabic network $N$ such that $X = X(N)$.
\end{definition}

\subsection{Pseudo symmetric networks}

It is important to emphasize that not all plabic graphs representing a point in $\lgnn$ are \tsym, as demonstrated in \Cref{alm-symm-vs-symm-example}. However, every such graph is \tsym up to a sequence local moves (see \Cref{local-moves-figure} or \cite{lam-textbook} for a more detailed explanation of these moves).

\begin{proposition}
    \label{networks-are-almost-symmetric}
    For a network $N$, we have $X(N) \in \lgnn$ if and only if $N$ is \term{pseudo symmetric}, i.e. is equivalent $\mainop(G)$.
\end{proposition}
\begin{proof}
    By \Cref{remark: tsym-equals-R-invar}, $X(N) \in \lgnn$ if and only if $X(N) = \mainop(X(N))$. By \Cref{lemma:swap-on-I-equals-swap-on-X}, we also get $\mainop(X(N)) = X(\mainop(N))$. Therefore, $X(N) \in \lgnn$ if and only if 
    $$X(N) = X(\mainop(N)).$$
    By \cite[Proposition~4.8 and Theorem~7.15]{lam-textbook}, this is equivalent to $N$ and $\mainop(N)$ being connected by local moves, i.e. to $N$ being pseudo symmetric.
\end{proof}

\begin{remark}
    We call these networks pseudo symmetric, because they are only symmetric up to move-equivalence. In reality these networks may lack any visible symmetries (see \Cref{two-graphs-figure}(left)). For this reason it is generally hard to check whether a network is almost symmetric without finding the corresponding point in the Lagrangian Grassmannian.

    In \Cref{symm-graphs}, we prove that every \tsym point can be represented with a (non-reduced) weighted graph that is actually physically rotationally symmetric (see \Cref{two-graphs-figure}(right)). 
\end{remark}

\subsection{Adding \tsym bridges}

Here, we introduce an analog of adding bridges (see \Cref{subsection: background-adding-bridges}) in the \tsym case. The idea is to add two bridges symmetric to each other simultaneously so that the corresponding point remains \tsym. 

\begin{definition}
    \label{def-symbridge}
    Let $G$ be a plabic graph. Denote $\symbridge{i}(G)$ to be the graph obtained from $G$ by adding one bridge between $i$ and $i+1$, white at $i$, and another bridge between $\plusn{i}$ and $\plusn{i+1}$, white at $\plusn{i+1}$. We say that $\symbridge{i}(G)$ was obtained from $G$ by adding a pair of symmetric bridges at $i$. 

    Similarly, if $N$ is a plabic network, let $\symbridge{i}(N,a)$ be the network obtained from $N$ by adding a pair of symmetric bridges at $i$, both of weight $a$.

    We introduce one exception: if $n = 2$ and $N$ is a graph with two lollipops, we simplify $\symbridge{i}(N)$ using edge contraction and parallel edge reduction to add one bridge of weight $2a$ instead, since it creates the same graph but with double-edge reduced in the middle. Similarly we define $\symbridge{i}(G)$ for a graph $G$ with two boundary lollipops.

    Finally, for $X \in \gr$ corresponding to the network $N$, denote $\symbridge{i}(X,a)$ to be the point corresponding to the network $\symbridge{i}(N,a)$. Equivalently, $\symbridge{i}(X,a) = X \cdot x_i(a) \cdot y_{\plusn{i}}(a)$ (see \Cref{formula-adding-bridge-x} and \Cref{formula-adding-bridge-y}).
\end{definition}  

\begin{remark}
    To simplify the picture, in many of the examples in this paper we apply moves M2, R2 and construction by the boundary such that $\symbridge{i}(G)$ (or $\symbridge{i}(N)$ does not have vertices of degree $2$ or double edges. For example, see \Cref{fig:bridge-construction-odd-even-faces}.
\end{remark}
    We use the same notation $\symbridge{i}$ for all three cases, since it will be clear from context which case we refer to. Moreover, this notation is the same as we used for \tsym permutations. That can be justified by \Cref{sym-adding-bridges}.

    Notice also that networks with some negative weights also correspond to points in $\gr$. Because of this, the definition of $\symbridge{i}(X,a)$ or $\symbridge{i}(N,a)$ also makes sense for $a \leq 0$.

\begin{proposition}
    \label{sym-adding-bridges}
    Suppose $X \in \lgnn$ is a \tsym point corresponding to the permutation $f$. Then for every $a>0$ the point $\symbridge{i}(X,a)$ is \tsym and corresponds to the permutation $\symcross{i}{f}$. 
\end{proposition}
\begin{proof}

    First, recall that the point $X$ is in $\lgnn$ if and only if $X \form X^T = 0$. Thus, we need to verify that 
    $$\symbridge{i}(X,a) \form (\symbridge{i}(X,a))^T = 0.$$
    Using the definition of $\symbridge{i}$, we compute
    \begin{align}
        (\symbridge{i}(X,a))^T \form \symbridge{i}(X,a) = X \cdot (x_i(a) y_{\plusn{i}}(a)) \cdot \form \cdot (x_i(a) y_{\plusn{i}}(a))^T \cdot X.
        \label{eq: sym-adding-bridges-proof}
    \end{align}
    
    It can be verified through a matrix multiplication that
    $$(x_i(a) y_{\plusn{i}}(a)) \cdot \form \cdot (x_i(a) y_{\plusn{i}}(a))^T = \form.$$
    Substituting this into \Cref{eq: sym-adding-bridges-proof}, we get
    $$X \cdot (x_i(a) y_{\plusn{i}}(a)) \cdot \form \cdot (x_i(a) y_{\plusn{i}}(a))^T \cdot X = X \cdot \form \cdot X^T.$$
    Finally, since $X \in \lgnn$, we have $X \cdot \form \cdot X^T = 0$, which shows that $\symbridge{i}(X,a) \in \lgnn$.

    For the second part of the statement, since $X$ corresponds to $f$, then by \Cref{add-bridges-permutations-background}, $X \cdot x_i(a)$ corresponds to $f \star s_i$. But then $\symbridge{i}(X,a) = X \cdot x_i(a) \cdot y_{i+n}(a)$ corresponds to $s_{i+n} \star (f \star s_i)$, which by \Cref{symcross-through-star-lemma} is equal to $\symcross{i}{f}$.
\end{proof}

We now describe how adding symmetric bridges changes Pl\"ucker coordinates using the operations $\alp$ and $\beti$ (see \Cref{alp-beti-definition}).

\begin{proposition}
    \label{minors-after-adding-bridges}
    The minors of $\symbridge{i}(X,a)$ are determined by
    $$\Delta_I(\symbridge{i}(X,a)) = \Delta_I(X) + a \Delta_{\alp I} (X) + a \Delta_{\bet I} (X) + a^2 \Delta_{\alp \bet I} (X).$$
\end{proposition}
\begin{proof}
    This proposition follows directly from \Cref{add-bridges-minors-background}. First, we get that 
    $$\Delta_I(X \cdot x_i(a)) = \Delta_I(X) + a \Delta_{\alp I} (X).$$
    Next, since $\symbridge{i}(X,a) = (X \cdot x_i(a))\cdot y_{i+n}(a)$, we get
    $$\Delta_I(\symbridge{i}(X,a)) = \Delta_I(X \cdot x_i(a)) + a \Delta_{\bet I} (X \cdot x_i(a)).$$
    Combining the two formulas produces the desired expression.
\end{proof}

\subsection{Removing \tsym bridges}

Now, similarly to \Cref{delete-bridges}, we explain how to remove bridges.
\begin{lemma}
    \label{sym-deleting-bridges}
    Suppose $Y\in \celllgn{g}$ is a \tsym point and $f\in \sbn$ is such that $g = \symcross{i}{f} > g$. Then there exists $c>0$ such that
    \begin{enumerate}[(1)]
        \item For $a< c$ ($a>0$), $\symbridge{i}(Y,-a) \in \celllgn{g}$
        \item For $a = c$, $\symbridge{i}(Y,-a) \in \celllgn{f}$.
        \item For $a>c$, $\symbridge{i}(Y,-a) \not\in \grnn$.
    \end{enumerate}
    Moreover, $c$ is the smallest root of
    $$P(a) = \Delta_J(Y) - a(\Delta_{\alp J}(Y) + \Delta_{\bet J}(Y)) + a^2 \Delta_{\alp \bet J}(Y),$$
    where $J = \Imin{i+1} (g)$.
\end{lemma}
\begin{proof}
    Define $T  = \{a \geq 0\st \symbridge{i}(Y,-a) \in \grnn\}$. By the argument from the proof of \Cref{sym-adding-bridges}, all $\symbridge{i}(Y,-a) \in \lgr$. First, we will show that $T = [0,c]$ for some $c \in (0, + \infty)$ by proving several claims.
    \begin{enumerate}[1.]
        \item $T$ is closed.

        This follows simply from $\symbridge{i}(Y,\cdot)$ being continuous and $\grnn$ being closed.

        \item For every $a \in T$, we have $[0,a) \subset T$. 

        For that, suppose $a \in T$ and $a' \in [0,a)$. Then we get
        $\symbridge{i}(Y,-a) \in \grnn$ (since $a \in T$) and $a - a' > 0$, so we obtain that
        $$\symbridge{i}(Y,-a') = \symbridge{i}(\symbridge{i}(Y,-a), a - a') \in \grnn.$$ 
        Therefore $a' \in T$.

        \item $T$ is bounded.

        Since $g = \symcross{i}{f} > f$, by \Cref{lemma: conditions-symmcross-i-f-greater-f}, we have $g s_i < g$ and $gs_i \in \Bnn$. But then we can take $c_0>0$ given by \Cref{delete-bridges}. We will show that all $a \in T$ satisfy $a \leq c_0$. 

        For that, suppose $a \in T$. Then 
        $$\symbridge{i}(Y,-a) = Y\cdot x_i(-a) \cdot y_{\plusn{i}}(-a)  \in \grnn.$$ 
        Consider also $X' = \symbridge{i}(Y,-a) \cdot y_{\plusn{i}}(a) \in \grnn$. On the other hand, we can rewrite 
        $$X' = Y \cdot x_i(-a) \cdot y_{\plusn{i}}(-1) \cdot y_{\plusn{i}}(a) = Y\cdot x_i(-a).$$ 
        But then by \Cref{delete-bridges}, since $X' \in \grnn$, we have $a \leq c_0$. Therefore $T \subset [0, c_0]$.

        \item There is $a> 0$ in $T$.

        For this, consider $c_0$ from \Cref{delete-bridges} again. Let $X' = Y \cdot x_i(-\frac{c_0}{2})$. By \Cref{delete-bridges}, we have $X' \in \cell{g}$. Now, since by \Cref{lemma: conditions-symmcross-i-f-greater-f} $s_{\plusn{i}}  g < g$, so by the analog of \Cref{delete-bridges}, we can find $a > 0$ (we may assume $a < \frac{c_0}{2}$) such that 
        $$X'' = X' \cdot y_{\plusn{i}}(-a) = Y \cdot x_i(-\frac{c_0}{2}) \cdot y_{\plusn{i}}(-a) \in \grnn.$$
        Finally, we get that
        \begin{align*}
            \symbridge{i}(Y, -a) = Y \cdot x_i(-a) \cdot y_{\plusn{i}}(-a) =\\
            Y \cdot x_i(-\frac{c_0}{2}) \cdot y_{\plusn{i}}(-a) \cdot x_i(\frac{c_0}{2} - a) = X'' \cdot x_i(\frac{c_0}{2} - a) \in \grnn
        \end{align*}
        (since $X'' \in \grnn$ and $\frac{c_0}{2} - a \geq 0$. That is, $\symbridge{i}(Y,-a) \in \grnn$, so $a \in T$.
    \end{enumerate}

    These statements combined imply that $T = [0,c]$ for $c \in (0, + \infty)$. 
    
    Now, let $a \in T$ and suppose $X =\symbridge{i}(Y,-a) \in \celllgn{h}$. then since $Y = \symbridge{i}(X,a) \in \celllgn{g}$, by \Cref{sym-adding-bridges}, $g = \symcross{i}{h}$. On the other hand, we have $g = \symcross{i}{f}$, which means that $h  = f$ or $g$. Therefore 
    $$T = \{a \geq 0\st \symbridge{i}(Y,-a) \in \celllgn{f}\} \sqcup \{a \geq 0\st \symbridge{i}(Y,-a) \in \celllgn{g}\}.$$
    We can finally obtain the statements of the lemma.
    \begin{enumerate}
        \item[(3)] By the definition of $T$, for every $a > c$, we get $\symbridge{i}(Y,-a) \not\in \grnn$.
        
        \item[(2)] Show that $\symbridge{i}(Y,-c) \in \celllgn{f}$. 

        Assume the contrary. But then the point $\symbridge{i}(Y,-c) \in \celllgn{g}$ satisfies the conditions of the lemma, so we can use the above again to find a positive $c'$ such that $\symbridge{i}(\symbridge{i}(Y,-c), - c') = \symbridge{i}(Y, - (c +c')) \in \grnn$, which would contradict with $c$ being the maximal number in $T$.

        \item[(1)] Finally, show that for every $a< c$, we have $\symbridge{i}(Y,-a) \in \celllgn{g}$.

        For that, notice that $\symbridge{i}(Y,-c) \in \celllgn{f}$, so then the point $\symbridge{i}(Y,-a) = \symbridge{i}(\symbridge{i}(Y,-c), c -a) $ is going to correspond to $\symcross{i}{f} = g$ (by \Cref{sym-adding-bridges}).
    \end{enumerate}

    Now, show that $c$ can be found as the smallest root of $P$. 
    
    If $X = \symbridge{i}(Y,-c) \in \celllgn{f}$, then by  a version of \Cref{minors-after-adding-bridges}, we get
        \begin{align}
        \label{deleting-symmetric-bridges}
            \Delta_I(X) = \Delta_I(Y) - c\Delta_{\alp I}(Y) - c\Delta_{\bet I}(Y) + c^2 \Delta_{\alp \bet I}(Y).
        \end{align}
    
        Denote $J = \Imin{i+1} (g)$. $Y \in \celllgn{g}$, which implies $\Delta_J(Y) \neq 0$. On the other hand, recall that $g = \symcross{i}{f} > f$, so by \Cref{symcross-through-star-lemma}, we get $g = s_{i+n} \star (f \star s_i) = s_{i+n} \star (fs_i) \geq fs_i > f$. Therefore by \Cref{Imin-gets-smaller-cor}, we get $J = \Imin{i+1}(g) <_{i+1} \Imin{i+1}(f)$. Therefore since $X \in \cell{f}$, we get $\Delta_J (X) = 0$.
    
        Therefore plugging in $I = J$ above, we have
    
        $$P(c) = \Delta_J(Y) - c(\Delta_{\alp J}(Y) + \Delta_{\bet J}(Y)) + c^2 \Delta_{\alp \bet J}(Y) = 0.$$

        Therefore $c$ is a root of $P$. It is now left to show that $c$ is the smallest root of $P$. Suppose $P$ has two roots and $c$ is the biggest and $c'$ is the smallest. But then $c' < c$, so by part(1) of the Lemma, $\symbridge{i}(Y, -c') \in \celllgn{g}$. However, since $c'$ is also a root of $P$, we get that $\Delta_J(\symbridge{i}(Y,-c')) = P(c') = 0$, which contradicts with $\symbridge{i}(Y, -c') \in \celllgn{g}$. This proves that $c$ has to be the smallest root of $P$.
  
\end{proof}

\begin{proposition}
    \label{cell-topology}
    Suppose $f, g \in \sbn$ such that $g = \symcross{i}{f} > f$. Then the restriction
    $$\symbridge{i}: (\celllgn{f}) \times \RR_{>0} \to \celllgn{g}$$
    is a homeomorphism. 

    Additionally, if $G$ is a \tsym graph that parameterizes $\celllgn{f}$, then $\symbridge{i}(G)$ parameterizes $\celllgn{g}$.
\end{proposition}
\begin{proof}

    First of all, since $\symbridge{i}(X,a) = X \cdot x_i(a) \cdot y_{\plusn{i}}(a)$, $\symbridge{i}$ is continuous. 

    To find $\symbridge{i}^{-1}$, take $Y \in \celllgn{g}$. One can see that $\symbridge{i}(X,a) = Y$ is equivalent to $X = \symbridge{i}(Y,-a)$ (we can see the equivalence by applying $\symbridge{i}(\cdot, -a)$ or $\symbridge{i}(\cdot, a)$ from both sides). So we need to find $a>0$ such that $\symbridge{i}(Y,-a) \in \celllgn{f}$. By \Cref{sym-deleting-bridges}, such $a$ is exists, is unique and depends on $Y$ continuously (we take $a = c$). Therefore there is also a unique $X = \symbridge{i}(Y,-a) \in \celllgn{f}$ such that $Y = \symbridge{i}(X,a)$. This shows that $\symbridge{i}^{-1}$ is well-defined and continuous, so $\symbridge{i}$ is a homeomorphism.

    Now, for the second part, suppose $G$ parametrizes $\celllgn{f}$ and let $Y \in \celllgn{g}$. Then there exists $X \in \celllgn{f}, a > 0$ such that $Y = \symbridge{i}(X,a)$. Since $G$ parametrizes $\celllgn{f}$, we can find a symmetric weighting $w$ on $G$ such that $X = X(G,w)$. But then $Y = X(\symbridge{i}(G), w \cup \{a\}$) (where we just assign weights $a$ to the two newly added bridges in $Y$), so any $Y \in \celllgn{f}$ corresponds to a weighting of $\symbridge{i}(G)$. This shows that $\symbridge{i}(G)$ parametrizes $\celllgn{g}$.
    
\end{proof}

\subsection{Adding fixed points}
    We also introduce the operation of adding symmetric fixed points similar to the operation of adding lollipops, see \Cref{subsection: background-adding-lollipops}. 
    \begin{definition}
    \begin{itemize}
        \item For $f \in \sbb(n-1,2n-2)$, denote $\lollipop{i}(f) \in \sbn$ to be the permutation obtained from $f$ by
        adding
        fixed points between $i-1$ and $i$ and between $i+n-2$ and $i+n-1$ such that $\lollipop{i}(f)(i) = i$ and $\lollipop{i}(f)(i+n) = i+3n$. That is,
        $$\lollipop{i}(f)(j) = \begin{cases}
            f(j), & 1 \leq j < i,\\
            i, & j = i,\\
            f(j-1), & i < j < \plusn{i},\\
            i+3n, & j = \plusn{i},\\
            f(j-2), & \plusn{i} < j \leq 2n.
        \end{cases}$$
        We say that $\lollipop{i}(f)$ is obtained from $f$ by adding a pair of symmetric points at $i$.

        \item For a a plabic graph $G$ with $2n-2$ boundary vertices, we let $\lollipop{i}(G)$ be the graph obtained from $G$ by attaching two lollipops, one black between $i-1$ and $i$ (we renumber so that the new vertex has number $i$) and one white between $i+n-1$ and $i+n-2$. 

        \item For a plabic network $N$ with $2n-2$ boundary vertices, we let $\lollipop{i}(N)$ be the network obtained from $N$ by attaching two lollipops, one black between $i-1$ and $i$ and one white between $i+n-1$ and $i+n-2$, both of weight $1$. 

        \item Finally, for $X \in \lgrnn(n-1,2n)$, denote $\lollipop{i}(X)$ to be the point that corresponds to $\lollipop{i}(N)$, where $N$ is a network representing $X$.
    \end{itemize}   
    \end{definition}

    It is easy to see the following.

    \begin{lemma}
        \label{lemma: symm-adding-lollipops}
        \label{lollipop-homeomorphism}
        If $X \in \lgnn$ corresponds to $f \in \sbn$, then $\lollipop{i}(X)$ corresponds to $\lollipop{i}(f)$. Moreover, the restriction 
        $$\lollipop{i}:\celllgn{f} \to \celllgn{\lollipop{i}(f)}$$
        is a homeomorphism.

        Additionally, if $G$ is a \tsym graph parameterizing $\celllgn{f}$, then $\lollipop{i}(G)$ parametrizes $\celllgn{\lollipop{i}(f)}$.
    \end{lemma}

\subsection{Proof of the main theorem}

The idea of proving \Cref{symm-graphs} is to construct the \tsym graph $G$ inductively, moving up the Bruhat order. With every step up we will be adding a pair of symmetric edges to $G$. 

\begin{proof}[Proof of \Cref{symm-graphs}]
\label{proof-symm-graphs}

    It was already proven in \Cref{permutations-are-symmetric} that the cell $\celllgn{f}$ is empty unless $f \in \sbn$. For $f \in \sbn$, we show the statements by the induction on $\symmdim(f) + n$.

    The case $\symmdim(f) + n = 1$ is trivial: $n$ must be equal to $1$ and $f$ has to only have fixed points, so we can parametrize $\celllgn{f}$ with a (\tsym) graph with two lollipops of different colors. We can also see that $\celllgn{f} \simeq \RR^0$.

    Now, if $\symmdim(f) + n \geq 2$, we may use \Cref{going-down-bruhat} to either find a a pair of symmetric fixed points $i$ and $i+n$ in $f$ or to find $i$ and $f' \in \sbn$ such that $f' < f$ and $f = \symcross{i}{f'}$.

    \begin{enumerate}
        \item If there exists $i \in [2n]$ such that $f(i) = i$ and $f(i+n) = i+3n$, there is $f' \in \sbb(n-1,2n-2)$ such that $f = \lollipop{i}(f')$. By induction hypothesis there exists a \tsym graph $G'$ that parametrizes $\celllgn{f'}$ and $\celllgn{f'} \simeq \RR^{\symmdim(f')}$. But then by \Cref{lemma: symm-adding-lollipops}, we get that $G:= \lollipop{i}(G')$ parametrizes $\celllgn{f}$ and that $\celllgn{f} \simeq \celllgn{f'} \simeq \RR^{\symmdim(f')} = \RR^{\symmdim(f)}$.
        
        \item If $f = \symcross{i}{f'} > f'$ for some $f' \in \sbn$, then $\celllgn{f'} \simeq \RR^{\symmdim(f')}$ and there exists a \tsym graph $G'$ that parametrizes $\celllgn{f'}$. But then by \Cref{cell-topology}, we get that $\celllgn{f} \simeq \RR^{\symmdim(f') + 1} = \RR^{\symmdim(f)}$ and that $\symbridge{i}(G)$ parametrizes $\celllgn{f}$.
    \end{enumerate}

\end{proof}

\begin{figure}
    \centering
    
\begin{center}
\scalebox{0.75}{
\begin{tabular}{ccccccccccc}

    {\begin{tikzpicture}[scale = \constructionscale,baseline={(0,0)}]
    \tiny
    \draw[line width=\constructioncircle,opacity = 0.3] (0,0) circle (3);
    
    \wdot{{3*cos(60)}}{{3*sin(60)}}; \bdot{{3*cos(0)}}{{3*sin(0)}}; \bdot{{3*cos(-60)}}{{3*sin(-60)}}; \bdot{{3*cos(-120)}}{{3*sin(-120)}}; \wdot{{3*cos(180)}}{{3*sin(180)}}; \wdot{{3*cos(120)}}{{3*sin(120)}};
    \node [above right] at ({3*cos(60)},{3*sin(60)}) {1};
    \node [right] at ({3*cos(0)},{3*sin(0)}) {2};
    \node [below right] at ({3*cos(-60)},{3*sin(-60)}) {3};
    \node [below left] at ({3*cos(-120)},{3*sin(-120)}) {4};
    \node [left] at ({3*cos(180)},{3*sin(180)}) {5};
    \node [above left] at ({3*cos(120)},{3*sin(120)}) { 6};
    \normalsize
\end{tikzpicture}
    }
    &
    {\hspace{-0.6cm} $\rightarrow$ 
    \hspace{-0.4cm}}
    &
    {\begin{tikzpicture}[scale = \constructionscale,baseline={(0,0)}]
    \tiny
    \draw[line width=\constructioncircle,opacity = 0.3] (0,0) circle (3);
    \draw[line width=\constructionline, blue] ({3*cos(60)},{3*sin(60)}) -- ({1*cos(60)},{1*sin(60)});
    \draw[line width=\constructionline, blue] ({-3*cos(60)},{-3*sin(60)}) --({-1*cos(60)},{-1*sin(60)});
    \draw[line width=1.5*\constructionline, red] ({1*cos(60)},{1*sin(60)})--({-1*cos(60)},{-1*sin(60)});
    
    \wdot{{3*cos(60)}}{{3*sin(60)}}; \bdot{{3*cos(0)}}{{3*sin(0)}}; \bdot{{3*cos(-60)}}{{3*sin(-60)}}; \bdot{{3*cos(-120)}}{{3*sin(-120)}}; \wdot{{3*cos(180)}}{{3*sin(180)}}; \wdot{{3*cos(120)}}{{3*sin(120)}};
    \bdot{{cos(60)}}{{sin(60)}};\wdot{{-cos(60)}}{{-sin(60)}}; 
    \node [above right] at ({3*cos(60)},{3*sin(60)}) {1};
    \node [right] at ({3*cos(0)},{3*sin(0)}) {2};
    \node [below right] at ({3*cos(-60)},{3*sin(-60)}) {3};
    \node [below left] at ({3*cos(-120)},{3*sin(-120)}) {4};
    \node [left] at ({3*cos(180)},{3*sin(180)}) {5};
    \node [above left] at ({3*cos(120)},{3*sin(120)}) { 6};
    \normalsize
\end{tikzpicture}}
    &
    {\hspace{-0.6cm} $\rightarrow$ 
    \hspace{-0.4cm}}
    &
    {\begin{tikzpicture}[scale = \constructionscale,baseline={(0,0)}]
    \tiny
    \draw[line width=\constructioncircle,opacity = 0.3] (0,0) circle (3);
    \draw[line width=\constructionline, blue] ({3*cos(60)},{3*sin(60)}) -- (0,.7);
    \draw[line width=\constructionline, blue] ({-3*cos(60)},{-3*sin(60)}) --(0,-.7);
    \draw[line width=\constructionline, blue] (0,.7)--(0,-.7);
    \draw[line width=1.5*\constructionline, red] (0,-.7)--({3*cos(-60)},{3*sin(-60)});
    \draw[line width=1.5*\constructionline, red] (0,.7)-- ({-3*cos(-60)},{-3*sin(-60)});
    
    \wdot{{3*cos(60)}}{{3*sin(60)}}; \bdot{{3*cos(0)}}{{3*sin(0)}}; \bdot{{3*cos(-60)}}{{3*sin(-60)}}; \bdot{{3*cos(-120)}}{{3*sin(-120)}}; \wdot{{3*cos(180)}}{{3*sin(180)}}; \wdot{{3*cos(120)}}{{3*sin(120)}};
    \bdot{0}{.7};\wdot{0}{-.7}; 
    \node [above right] at ({3*cos(60)},{3*sin(60)}) {1};
    \node [right] at ({3*cos(0)},{3*sin(0)}) {2};
    \node [below right] at ({3*cos(-60)},{3*sin(-60)}) {3};
    \node [below left] at ({3*cos(-120)},{3*sin(-120)}) {4};
    \node [left] at ({3*cos(180)},{3*sin(180)}) {5};
    \node [above left] at ({3*cos(120)},{3*sin(120)}) { 6};
    \normalsize
\end{tikzpicture}} 
    & 
    {\hspace{-0.6cm} $\rightarrow$ 
    \hspace{-0.4cm}} 
    &
    {\begin{tikzpicture}[scale = \constructionscale,baseline={(0,0)}]
    \tiny
    \draw[line width=\constructioncircle,opacity = 0.3] (0,0) circle (3);
    \draw[line width=\constructionline, blue] ({3*cos(60)},{3*sin(60)}) -- ({2*cos(60)},{2*sin(60)})--(0,.7)--(-1.1,-.3);
    \draw[line width=\constructionline, blue] ({-3*cos(60)},{-3*sin(60)}) -- ({-2*cos(60)},{-2*sin(60)})--(0,-.7)--(1.1,.3);
    \draw[line width=\constructionline, blue] (0,.7)--(0,-.7);
    \draw[line width=\constructionline, blue] (1.1,.3) --({3*cos(-60)},{3*sin(-60)});
    \draw[line width=\constructionline, blue] (-1.1,-.3) -- ({-3*cos(-60)},{-3*sin(-60)});
    
    \bdot{{3*cos(60)}}{{3*sin(60)}}; \bdot{{3*cos(0)}}{{3*sin(0)}}; \wdot{{3*cos(-60)}}{{3*sin(-60)}}; \wdot{{3*cos(-120)}}{{3*sin(-120)}}; \wdot{{3*cos(180)}}{{3*sin(180)}}; \bdot{{3*cos(120)}}{{3*sin(120)}};\wdot{{2*cos(60)}}{{2*sin(60)}};\bdot{0}{.7}; \wdot{-1.1}{-.3};\bdot{{2*cos(-120)}}{{2*sin(-120)}};\wdot{0}{-.7}; \bdot{1.1}{.3};
    \node [above right] at ({3*cos(60)},{3*sin(60)}) {1};
    \node [right] at ({3*cos(0)},{3*sin(0)}) {2};
    \node [below right] at ({3*cos(-60)},{3*sin(-60)}) {3};
    \node [below left] at ({3*cos(-120)},{3*sin(-120)}) {4};
    \node [left] at ({3*cos(180)},{3*sin(180)}) {5};
    \node [above left] at ({3*cos(120)},{3*sin(120)}) { 6};
    \normalsize
\end{tikzpicture}
    }
    &
    {\hspace{-0.6cm} $\rightarrow$
    \hspace{-0.4cm}}
    &
    {\begin{tikzpicture}[scale = \constructionscale,baseline={(0,0)}]
    \tiny
    \draw[line width=\constructioncircle,opacity = 0.3] (0,0) circle (3);
    \draw[line width=\constructionline, blue] ({3*cos(60)},{3*sin(60)}) -- ({2*cos(60)},{2*sin(60)})--(0,.7)--(-1.1,-.3);
    \draw[line width=1.5*\constructionline, red](-1.1,-.3)--({-2*cos(60)},{-2*sin(60)});
    \draw[line width=\constructionline, blue] ({-3*cos(60)},{-3*sin(60)}) -- ({-2*cos(60)},{-2*sin(60)})--(0,-.7)--(1.1,.3);
    \draw[line width=1.5*\constructionline, red](1.1,.3)--({2*cos(60)},{2*sin(60)});
    \draw[line width=\constructionline, blue] (0,.7)--(0,-.7);
    \draw[line width=\constructionline, blue] (1.1,.3) --({3*cos(-60)},{3*sin(-60)});
    \draw[line width=\constructionline, blue] (-1.1,-.3) -- ({-3*cos(-60)},{-3*sin(-60)});
    
    \bdot{{3*cos(60)}}{{3*sin(60)}}; \bdot{{3*cos(0)}}{{3*sin(0)}}; \wdot{{3*cos(-60)}}{{3*sin(-60)}}; \wdot{{3*cos(-120)}}{{3*sin(-120)}}; \wdot{{3*cos(180)}}{{3*sin(180)}}; \bdot{{3*cos(120)}}{{3*sin(120)}};\wdot{{2*cos(60)}}{{2*sin(60)}};\bdot{0}{.7}; \wdot{-1.1}{-.3};\bdot{{2*cos(-120)}}{{2*sin(-120)}};\wdot{0}{-.7}; \bdot{1.1}{.3};
    \node [above right] at ({3*cos(60)},{3*sin(60)}) {1};
    \node [right] at ({3*cos(0)},{3*sin(0)}) {2};
    \node [below right] at ({3*cos(-60)},{3*sin(-60)}) {3};
    \node [below left] at ({3*cos(-120)},{3*sin(-120)}) {4};
    \node [left] at ({3*cos(180)},{3*sin(180)}) {5};
    \node [above left] at ({3*cos(120)},{3*sin(120)}) { 6};
    \normalsize
\end{tikzpicture}}
    &
    {\hspace{-0.6cm} $\rightarrow$ 
    \hspace{-0.4cm}}
    &
    {\begin{tikzpicture}[scale = \constructionscale,baseline={(0,0)}]
    \tiny
    \draw[line width=\constructioncircle,opacity = 0.3] (0,0) circle (3);
    \draw[line width=\constructionline, blue] ({3*cos(60)},{3*sin(60)}) -- ({2*cos(60)},{2*sin(60)})--(0,.7)--(-1.1,-.3);
    \draw[line width=\constructionline, blue](-1.1,-.3)--({-2*cos(60)},{-2*sin(60)});
    \draw[line width=\constructionline, blue] ({-3*cos(60)},{-3*sin(60)}) -- ({-2*cos(60)},{-2*sin(60)})--(0,-.7)--(1.1,.3);
    \draw[line width=\constructionline, blue](1.1,.3)--({2*cos(60)},{2*sin(60)});
    \draw[line width=\constructionline, blue] (0,.7)--(0,-.7);
    \draw[line width=\constructionline, blue] (1.1,.3) --({2*cos(-30)},{2*sin(-30)});
    \draw[line width=\constructionline, blue] (-1.1,-.3) -- ({-2*cos(-30)},{-2*sin(-30)});
    \draw[line width=1.5*\constructionline, red] ({3*cos(0)},{3*sin(0)}) -- ({2*cos(-30)},{2*sin(-30)});
    \draw[line width=\constructionline, blue]({2*cos(-30)},{2*sin(-30)})-- ({3*cos(-60)},{3*sin(-60)});
    \draw[line width=1.5*\constructionline, red] ({-3*cos(0)},{-3*sin(0)}) -- ({-2*cos(-30)},{-2*sin(-30)});
    \draw[line width=\constructionline, blue]({-2*cos(-30)},{-2*sin(-30)})-- ({-3*cos(-60)},{-3*sin(-60)});
    
    \bdot{{3*cos(60)}}{{3*sin(60)}}; \bdot{{3*cos(0)}}{{3*sin(0)}}; \bdot{{3*cos(-60)}}{{3*sin(-60)}}; \wdot{{3*cos(-120)}}{{3*sin(-120)}}; \wdot{{3*cos(180)}}{{3*sin(180)}}; \wdot{{3*cos(120)}}{{3*sin(120)}};\wdot{{2*cos(60)}}{{2*sin(60)}};\bdot{0}{.7}; \wdot{-1.1}{-.3};\bdot{{2*cos(-120)}}{{2*sin(-120)}};\wdot{{2*cos(-30)}}{{2*sin(-30)}};\bdot{{-2*cos(-30)}}{{-2*sin(-30)}};\wdot{0}{-.7}; \bdot{1.1}{.3};
    \node [above right] at ({3*cos(60)},{3*sin(60)}) {1};
    \node [right] at ({3*cos(0)},{3*sin(0)}) {2};
    \node [below right] at ({3*cos(-60)},{3*sin(-60)}) {3};
    \node [below left] at ({3*cos(-120)},{3*sin(-120)}) {4};
    \node [left] at ({3*cos(180)},{3*sin(180)}) {5};
    \node [above left] at ({3*cos(120)},{3*sin(120)}) { 6};
    \normalsize
\end{tikzpicture}}
\end{tabular}
}
\end{center}
    \caption{The process of constructing a symmetric graph}
\end{figure}
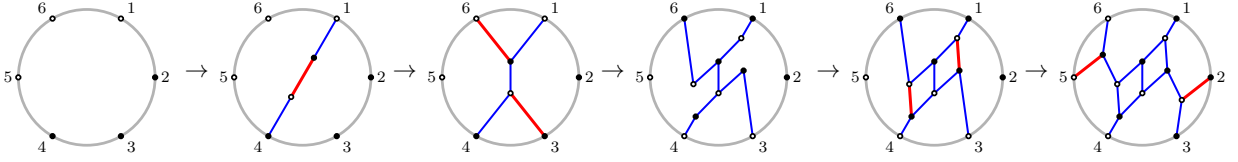

\subsection{Can \tsym graphs be reduced?}
\label{subsection: tsym-reduced}

Recall that reduced graphs are those with the minimal possible number of faces. So we wish to study the minimal number of faces for a \tsym graph representing $f$.

\begin{lemma}
    \label{symm-bounding-number-of-faces-lemma}
    Suppose $f \in \sbn$ and $G$ is a \tsym graph parameterizing $\celllgn{f}$. Then 
    $$\# \text{faces} (G) \geq \begin{cases}
        2\symmdim(f), &\exists i, f(i) = i+n;\\
        2\symmdim(f) + 1, & \text{otherwise}.
    \end{cases}.$$
\end{lemma}

\begin{proof}
    Since $G$ parametrizes $\celllgn{f}$, different choices of positive face weights (with the restriction of symmetric pairs having the product of $1$) correspond to all possible points in $\celllgn{f}$. This provides a continuous surjection $\RR^{\lfloor{\frac{\# \text{faces}}{2}}\rfloor} \to \celllgn{f} \simeq \RR^{\symmdim(f)}$. Therefore $\lfloor{\frac{\# \text{faces} (G)}{2}}\rfloor \geq \symmdim(f)$, so $\# \text{faces} (G) \geq 2 \symmdim(f)$. This proves the statement in the first case.

    For the second case, suppose $f(i) \neq i+n$ for all $i$. By \Cref{d-vs-symmdim}, we have
    $$d(f) = 2 \symmdim(f) -  \frac{\# \{i \in [2n] \st f(i) = i+n\}}{2} =  2 \symmdim(f).$$
    But then since the reduced graph representing $\cell{f}$ has the minimal possible number of faces, $d(f) + 1$ (see \Cref{reduced-number-of-faces-lemma}), we get
    $$\# \text{faces} (G) \geq d(f) + 1 = 2\symmdim(f) + 1.$$
    This proves the inequality in the second case as well.
    
\end{proof}

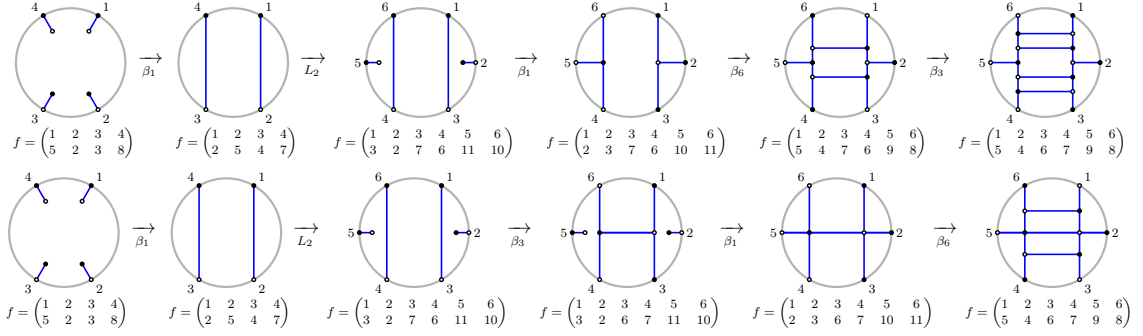
\begin{figure}
    \centering
    \begin{center}
\scalebox{0.6}{
\begin{tabular}{ccccccccccc}

    {\newcommand{\xx}{2}
\begin{tikzpicture}[scale = \constructionscale,baseline={(0,0)}]
    \tiny

    \node [above right] at ({3*cos(60)},{3*sin(60)}) {1};
    \node [below right] at ({3*cos(-60)},{3*sin(-60)}) {2};
    \node [below left] at ({3*cos(-120)},{3*sin(-120)}) {3};
    \node [above left] at ({3*cos(120)},{3*sin(120)}) { 4};

    \draw[line width=\constructioncircle,opacity = 0.3] (0,0) circle (3);
    \normalsize

    \symmlineblue{{3*cos(60)}}{{3*sin(60)}}{{\xx*cos(60)}}{{\xx*sin(60)}}{\constructionline};
    \symmlineblue{{3*cos(60)}}{{-3*sin(60)}}{{\xx*cos(60)}}{{-\xx*sin(60)}}{\constructionline};

    \symmwdot{{\xx*cos(60)}}{{\xx*sin(60)}};
    \symmbdot{{\xx*cos(60)}}{{-\xx*sin(60)}};

    \symmbdotsix{1};
    \symmbdotsix{6};
    
\end{tikzpicture}
    }
    &
    {\hspace{-0.6cm} $\xrightarrow[\symbridge{1}]{}$ \hspace{-0.4cm}}
    &
    {\newcommand{\xx}{2.3}
\begin{tikzpicture}[scale = \constructionscale,baseline={(0,0)}]
    \tiny

    \node [above right] at ({3*cos(60)},{3*sin(60)}) {1};
    \node [below right] at ({3*cos(-60)},{3*sin(-60)}) {2};
    \node [below left] at ({3*cos(-120)},{3*sin(-120)}) {3};
    \node [above left] at ({3*cos(120)},{3*sin(120)}) { 4};

    \draw[line width=\constructioncircle,opacity = 0.3] (0,0) circle (3);
    \normalsize

    \symmlineblue{{3*cos(60)}}{{3*sin(60)}}{{3*cos(60)}}{{-3*sin(60)}}{\constructionline};

    \symmbdotsix{1};
    \symmbdotsix{6};
    
\end{tikzpicture}
    }
    &
    {\hspace{-0.6cm} $\xrightarrow[\lollipop{2}]{}$ \hspace{-0.4cm}}
    &
    {\newcommand{\xx}{2.3}
\begin{tikzpicture}[scale = \constructionscale,baseline={(0,0)}]
    \tiny
    \draw[line width=\constructioncircle,opacity = 0.3] (0,0) circle (3);
    \circlenodessix;
    \normalsize

    \symmlineblue{{3*cos(60)}}{{3*sin(60)}}{{3*cos(60)}}{{-3*sin(60)}}{\constructionline};
    \symmlineblue{\xx}{0}{3}{{0}}{\constructionline};

    \symmbdot{\xx}{0};

    \symmbdotsix{1};
    \symmwdotsix{2};
    \symmwdotsix{3};
    
\end{tikzpicture}
    }
    &
    {\hspace{-0.6cm} $\xrightarrow[\symbridge{1}]{}$ \hspace{-0.4cm}}
    &
    {\begin{tikzpicture}[scale = \constructionscale,baseline={(0,0)}]
    \tiny
    \draw[line width=\constructioncircle,opacity = 0.3] (0,0) circle (3);
    \circlenodessix;
    \normalsize

    \symmlineblue{{3*cos(60)}}{{3*sin(60)}}{{3*cos(60)}}{{-3*sin(60)}}{\constructionline};
    \symmlineblue{{3*cos(60)}}{0}{3}{{0}}{\constructionline};

    \symmbdot{{-3*cos(60)}}{0};

    \symmbdotsix{1};
    \symmbdotsix{2};
    \symmbdotsix{3};
    
\end{tikzpicture}
    }
    &
    {\hspace{-0.6cm} $\xrightarrow[\symbridge{6}]{}$ \hspace{-0.4cm}}
    &
    {\newcommand{\yy}{0.8}
\begin{tikzpicture}[scale = \constructionscale,baseline={(0,0)}]
    \tiny
    \draw[line width=\constructioncircle,opacity = 0.3] (0,0) circle (3);
    \circlenodessix;
    \normalsize

    \symmlineblue{{3*cos(60)}}{{3*sin(60)}}{{3*cos(60)}}{{-3*sin(60)}}{\constructionline};
    \symmlineblue{{3*cos(60)}}{\yy}{{-3*cos(60)}}{\yy}{\constructionline};
    \symmlineblue{{3*cos(60)}}{0}{3}{{0}}{\constructionline};

    \symmbdot{{-3*cos(60)}}{0};
    \symmbdot{{3*cos(60)}}{\yy};
    \symmbdot{{3*cos(60)}}{-\yy};

    \symmwdotsix{1};
    \symmbdotsix{2};
    \symmwdotsix{3};
    
\end{tikzpicture}
    }
    &
    {\hspace{-0.6cm} $\xrightarrow[\symbridge{3}]{}$ \hspace{-0.4cm}}
    &
    {\newcommand{\yy}{0.8}
\begin{tikzpicture}[scale = \constructionscale,baseline={(0,0)}]
    \tiny
    \draw[line width=\constructioncircle,opacity = 0.3] (0,0) circle (3);
    \circlenodessix;
    \normalsize

    \symmlineblue{{3*cos(60)}}{{3*sin(60)}}{{3*cos(60)}}{{-3*sin(60)}}{\constructionline};
    \symmlineblue{{3*cos(60)}}{\yy}{{-3*cos(60)}}{\yy}{\constructionline};
    \symmlineblue{{3*cos(60)}}{{2*\yy}}{{-3*cos(60)}}{{2*\yy}}{\constructionline};
    \symmlineblue{{3*cos(60)}}{0}{3}{{0}}{\constructionline};

    \symmbdot{{-3*cos(60)}}{0};
    \symmbdot{{3*cos(60)}}{\yy};
    \symmbdot{{3*cos(60)}}{-\yy};
    \symmbdot{{-3*cos(60)}}{{2*\yy}};
    \symmbdot{{-3*cos(60)}}{{-2*\yy}};

    \symmbdotsix{1};
    \symmbdotsix{2};
    \symmbdotsix{3};
    
\end{tikzpicture}}\\ 
    {\tiny{$f = \begin{pmatrix}
        1 & 2 & 3 & 4\\
        5 & 2 & 3 & 8
    \end{pmatrix}$}} &&
    {\tiny{$f = \begin{pmatrix}
        1 & 2 & 3 & 4\\
        2 & 5 & 4 & 7
    \end{pmatrix}$}} &&
    {\tiny{$f = \begin{pmatrix}
        1 & 2 & 3 & 4 & 5 & 6\\
        3 & 2 & 7 & 6 & 11 & 10
    \end{pmatrix}$}} &&
    {\tiny{$f = \begin{pmatrix}
        1 & 2 & 3 & 4 & 5 & 6\\
        2 & 3 & 7 & 6 & 10 & 11
    \end{pmatrix}$}} &&
    {\tiny{$f = \begin{pmatrix}
        1 & 2 & 3 & 4 & 5 & 6\\
        5 & 4 & 7 & 6 & 9 & 8
    \end{pmatrix}$}} && 
    {\tiny{$f = \begin{pmatrix}
        1 & 2 & 3 & 4 & 5 & 6\\
        5 & 4 & 6 & 7 & 9 & 8
    \end{pmatrix}$}} 
    \normalsize
\end{tabular}
}

\end{center}

    \begin{center}
\scalebox{0.6}{
\begin{tabular}{ccccccccccc}

    {\newcommand{\xx}{2}
\begin{tikzpicture}[scale = \constructionscale,baseline={(0,0)}]
    \tiny

    \node [above right] at ({3*cos(60)},{3*sin(60)}) {1};
    \node [below right] at ({3*cos(-60)},{3*sin(-60)}) {2};
    \node [below left] at ({3*cos(-120)},{3*sin(-120)}) {3};
    \node [above left] at ({3*cos(120)},{3*sin(120)}) { 4};

    \draw[line width=\constructioncircle,opacity = 0.3] (0,0) circle (3);
    \normalsize

    \symmlineblue{{3*cos(60)}}{{3*sin(60)}}{{\xx*cos(60)}}{{\xx*sin(60)}}{\constructionline};
    \symmlineblue{{3*cos(60)}}{{-3*sin(60)}}{{\xx*cos(60)}}{{-\xx*sin(60)}}{\constructionline};

    \symmwdot{{\xx*cos(60)}}{{\xx*sin(60)}};
    \symmbdot{{\xx*cos(60)}}{{-\xx*sin(60)}};

    \symmbdotsix{1};
    \symmbdotsix{6};
    
\end{tikzpicture}
    }
    &
    {\hspace{-0.6cm} $\xrightarrow[\symbridge{1}]{}$ \hspace{-0.4cm}}
    &
    {\newcommand{\xx}{2.3}
\begin{tikzpicture}[scale = \constructionscale,baseline={(0,0)}]
    \tiny

    \node [above right] at ({3*cos(60)},{3*sin(60)}) {1};
    \node [below right] at ({3*cos(-60)},{3*sin(-60)}) {2};
    \node [below left] at ({3*cos(-120)},{3*sin(-120)}) {3};
    \node [above left] at ({3*cos(120)},{3*sin(120)}) { 4};

    \draw[line width=\constructioncircle,opacity = 0.3] (0,0) circle (3);
    \normalsize

    \symmlineblue{{3*cos(60)}}{{3*sin(60)}}{{3*cos(60)}}{{-3*sin(60)}}{\constructionline};

    \symmbdotsix{1};
    \symmbdotsix{6};
    
\end{tikzpicture}
    }
    &
    {\hspace{-0.6cm} $\xrightarrow[\lollipop{2}]{}$ \hspace{-0.4cm}}
    &
    {\newcommand{\xx}{2.3}
\begin{tikzpicture}[scale = \constructionscale,baseline={(0,0)}]
    \tiny
    \draw[line width=\constructioncircle,opacity = 0.3] (0,0) circle (3);
    \circlenodessix;
    \normalsize

    \symmlineblue{{3*cos(60)}}{{3*sin(60)}}{{3*cos(60)}}{{-3*sin(60)}}{\constructionline};
    \symmlineblue{3}{0}{\xx}{{0}}{\constructionline};
    
    \symmbdot{\xx}{0}

    \symmbdotsix{1};
    \symmwdotsix{2};
    \symmwdotsix{3};
    
\end{tikzpicture}
    }
    &
    {\hspace{-0.6cm} $\xrightarrow[\symbridge{3}]{}$ \hspace{-0.4cm}}
    &
    {\newcommand{\xx}{2.3}
\begin{tikzpicture}[scale = \constructionscale,baseline={(0,0)}]
    \tiny
    \draw[line width=\constructioncircle,opacity = 0.3] (0,0) circle (3);
    \circlenodessix;
    \normalsize

    \symmlineblue{{3*cos(60)}}{{3*sin(60)}}{{3*cos(60)}}{{-3*sin(60)}}{\constructionline};
    \symmlineblue{1.5}{0}{-1.5}{{0}}{\constructionline};
    \symmlineblue{3}{0}{\xx}{{0}}{\constructionline};

    \symmbdot{{-3*cos(60)}}{0};
    \symmbdot{\xx}{0}

    \symmbdotsix{1};
    \symmwdotsix{2};
    \symmbdotsix{3};
    
\end{tikzpicture}
    }
    &
    {\hspace{-0.6cm} $\xrightarrow[\symbridge{1}]{}$ \hspace{-0.4cm}}
    &
    {\begin{tikzpicture}[scale = \constructionscale,baseline={(0,0)}]
    \tiny
    \draw[line width=\constructioncircle,opacity = 0.3] (0,0) circle (3);
    \circlenodessix;
    \normalsize

    \symmlineblue{{3*cos(60)}}{{3*sin(60)}}{{3*cos(60)}}{{-3*sin(60)}}{\constructionline};
    \symmlineblue{3}{0}{-3}{{0}}{\constructionline};

    \symmbdot{{-3*cos(60)}}{0};

    \symmbdotsix{1};
    \symmbdotsix{2};
    \symmbdotsix{3};
    
\end{tikzpicture}
    }
    &
    {\hspace{-0.6cm} $\xrightarrow[\symbridge{6}]{}$ \hspace{-0.4cm}}
    &
    {\newcommand{\yy}{1.2}
\begin{tikzpicture}[scale = \constructionscale,baseline={(0,0)}]
    \tiny
    \draw[line width=\constructioncircle,opacity = 0.3] (0,0) circle (3);
    \circlenodessix;
    \normalsize

    \symmlineblue{{3*cos(60)}}{{3*sin(60)}}{{3*cos(60)}}{{-3*sin(60)}}{\constructionline};
    \symmlineblue{{3*cos(60)}}{\yy}{{-3*cos(60)}}{\yy}{\constructionline};
    \symmlineblue{3}{0}{-3}{{0}}{\constructionline};

    \symmbdot{{-3*cos(60)}}{0};
    \symmbdot{{3*cos(60)}}{\yy};
    \symmbdot{{3*cos(60)}}{-\yy};

    \symmwdotsix{1};
    \symmbdotsix{2};
    \symmwdotsix{3};
    
\end{tikzpicture}}\\ 
    {\tiny{$f = \begin{pmatrix}
        1 & 2 & 3 & 4\\
        5 & 2 & 3 & 8
    \end{pmatrix}$}} &&
    {\tiny{$f = \begin{pmatrix}
        1 & 2 & 3 & 4\\
        2 & 5 & 4 & 7
    \end{pmatrix}$}} &&
    {\tiny{$f = \begin{pmatrix}
        1 & 2 & 3 & 4 & 5 & 6\\
        3 & 2 & 7 & 6 & 11 & 10
    \end{pmatrix}$}} &&
    {\tiny{$f = \begin{pmatrix}
        1 & 2 & 3 & 4 & 5 & 6\\
        3 & 2 & 6 & 7 & 11 & 10
    \end{pmatrix}$}} &&
    {\tiny{$f = \begin{pmatrix}
        1 & 2 & 3 & 4 & 5 & 6\\
        2 & 3 & 6 & 7 & 10 & 11
    \end{pmatrix}$}} && 
    {\tiny{$f = \begin{pmatrix}
        1 & 2 & 3 & 4 & 5 & 6\\
        5 & 4 & 6 & 7 & 9 & 8
    \end{pmatrix}$}} 
    \normalsize
\end{tabular}
}

\end{center}

    \caption{An example of inductively constructing a \tsym \tred graph corresponding to the permutation $f = \begin{pmatrix}
        1 & 2 & 3 & 4 & 5 & 6\\
        5 & 4 & 6 & 7 & 9 & 8
    \end{pmatrix}$ in two ways (following the bridge construction).}
    \label{fig:bridge-construction-odd-even-faces}
\end{figure}

\begin{lemma}
    \label{lemma:middle-face}
    Suppose $f \in \sbn$. Then there exists a \tsym graph $G$ parameterizing $\celllgn{f}$ (obtained by the bridge construction) such that
    $$\# \text{faces} (G) = \begin{cases}
        2\symmdim(f_G), &\exists i, f_G(i) = i+n;\\
        2\symmdim(f_G) + 1, & \text{otherwise}.
    \end{cases}$$
\end{lemma}
\begin{remark}
    Note that not every \tsym graph obtained through the bridge construction will satisfy these conditions as demonstrated in \Cref{fig:bridge-construction-odd-even-faces}. Specifically even though both graphs are obtained through the bridge construction, only one of the graphs (the one on the bottom) satisfies the condition of \Cref{lemma:middle-face}.
\end{remark}

\begin{figure}
    \centering
    \begin{center}

\begin{tabularx}{\textwidth}{XXXX}

    {\newcommand{\xx}{1.8}
\begin{tikzpicture}[scale = \graphscale,baseline={(0,0)}]
    \tiny
    \draw[line width=\graphcircle,opacity = 0.3] (0,0) circle (3);
    \node [above] at ({3*cos(90)},{3*sin(90)}) (1) {1};
    \node [below] at ({3*cos(-90)},{3*sin(-90)}) (2) {2};
    
    \draw[line width=\graphline,blue] (0,3) -- (0,-3);

    \bdot{0}{3};
    \wdot{0}{-3};
    \normalsize
\end{tikzpicture}}
    &
    {
\begin{tikzpicture}[scale = \graphscale, baseline={(0,0)}]
        \draw[line width=\graphcircle,opacity=0.3] (0,0) circle (3);

        \draw[line width=\graphline,blue] (2.12132034356,2.12132034356) -- (0,1) -- (0,-1) -- (-2.12132034356,-2.12132034356);
        \draw[line width=\graphline,blue] (0,-1) -- (2.12132034356,-2.12132034356);
        \draw[line width=\graphline,blue] (-2.12132034356,2.12132034356) -- (0,1);

        \bdot{2.12132034356}{2.12132034356}; \wdot{2.12132034356}{-2.12132034356};\wdot{-2.12132034356}{-2.12132034356};\bdot{-2.12132034356}{2.12132034356};
        \bdot{0}{-1}; \wdot{0}{1};
        

        \scriptsize
        \node [above right] at (2.12132034356,2.12132034356) {1};
        \node [below right] at (2.12132034356,-2.12132034356) {2};
        \node [below left] at (-2.12132034356,-2.12132034356) {3};
        \node [above left] at (-2.12132034356,2.12132034356) {4};
        
\end{tikzpicture}

}
    &
    {\newcommand{\xx}{1.8}
\begin{tikzpicture}[scale = \graphscale,baseline={(0,0)}]
    \tiny
    \draw[line width=\graphcircle,opacity = 0.3] (0,0) circle (3);
    \node [above right] at ({3*cos(60)},{3*sin(60)}) (1) {1};
    \node [right] at ({3*cos(0)},{3*sin(0)}) (2) {2};
    \node [below right] at ({3*cos(-60)},{3*sin(-60)}) (3) {3};
    \node [below left] at ({3*cos(-120)},{3*sin(-120)}) (4) {4};
    \node [left] at ({3*cos(180)},{3*sin(180)}) (5) {5};
    \node [above left] at ({3*cos(120)},{3*sin(120)}) (6) {6};

    \draw[line width=\graphline,blue] ({\xx*cos(60)},{\xx*sin(60)}) -- ({\xx*cos(0)},{\xx*sin(0)}) -- ({\xx*cos(-60)},{\xx*sin(-60)}) -- ({\xx*cos(-120)},{\xx*sin(-120)}) -- ({\xx*cos(180)},{\xx*sin(180)}) -- ({\xx*cos(120)},{\xx*sin(120)}) -- ({\xx*cos(60)},{\xx*sin(60)});
    \draw[line width=\graphline,blue] ({3*cos(0)},{3*sin(0)}) -- ({\xx*cos(0)},{\xx*sin(0)});
    \draw[line width=\graphline,blue] ({3*cos(180)},{3*sin(180)}) -- ({\xx*cos(180)},{\xx*sin(180)});
    \draw[line width=\graphline,blue] ({3*cos(60)},{3*sin(60)}) -- ({\xx*cos(60)},{\xx*sin(60)});
    \draw[line width=\graphline,blue] ({3*cos(-60)},{3*sin(-60)}) -- ({3*cos(120)},{3*sin(120)}) ;
    \draw[line width=\graphline,blue] ({3*cos(-120)},{3*sin(-120)}) -- ({\xx*cos(-120)},{\xx*sin(-120)});

    \wdot{{3*cos(60)}}{{3*sin(60)}}; 
    \bdot{{3*cos(0)}}{{3*sin(0)}}; 
    \wdot{{3*cos(-60)}}{{3*sin(-60)}}; 
    \bdot{{3*cos(-120)}}{{3*sin(-120)}}; 
    \wdot{{3*cos(180)}}{{3*sin(180)}}; 
    \bdot{{3*cos(120)}}{{3*sin(120)}};

    \bdot{{\xx*cos(60)}}{{\xx*sin(60)}}; \wdot{{\xx*cos(0)}}{{\xx*sin(0)}}; \bdot{{\xx*cos(-60)}}{{\xx*sin(-60)}}; \wdot{{\xx*cos(-120)}}{{\xx*sin(-120)}}; \bdot{{\xx*cos(180)}}{{\xx*sin(180)}}; \wdot{{\xx*cos(120)}}{{\xx*sin(120)}};
    \normalsize
\end{tikzpicture}}
    &
    {


\newcommand{\ang}{45/2}

\newcommand{\cwdot}[1]{\wdot{{3*cos(45*(3-#1) - 45/2)}}{{3*sin(45*(3-#1) - 45/2)}}}
\newcommand{\cbdot}[1]{\bdot{{3*cos(45*(3-#1) - 45/2)}}{{3*sin(45*(3-#1) - 45/2)}}}
\newcommand{\cpointx}[1]{3*cos(45*(3-#1) - 45/2)}
\newcommand{\cpointy}[1]{3*sin(45*(3-#1) - 45/2)}

\pgfmathsetmacro{\xx}{0.55}

\newcommand{\wdotsym}[2]{\wdot{#1}{#2}; \bdot{-#1}{-#2}}
\newcommand{\bdotsym}[2]{\bdot{#1}{#2}; \wdot{-#1}{-#2}}

\newcommand{\newpropxxpointx}[6]{%
    {(#3 * #1 + #4 * #5)}
}
\newcommand{\newpropxxpointy}[6]{%
    {(#3 * #2 + #4 * #6)}%
}

\newcommand{\propx}[3]{
    \newpropxxpointx{\cpointx{#1}}{\cpointy{#1}}{#2}{#3}{\xx}{0}
}
\newcommand{\propxminus}[3]{
    -\newpropxxpointx{\cpointx{#1}}{\cpointy{#1}}{#2}{#3}{\xx}{0}
}
\newcommand{\propy}[3]{
    \newpropxxpointy{\cpointx{#1}}{\cpointy{#1}}{#2}{#3}{\xx}{0}
}

\newcommand{\wproppsym}[3]{
\wdotsym{\propx{#1}{#2}{#3}}{\propy{#1}{#2}{#3}}
}

\newcommand{\bproppsym}[3]{
\bdotsym{\newpropxxpointx{\cpointx{#1}}{\cpointy{#1}}{#2}{#3}{\xx}{0}}{\newpropxxpointy{\cpointx{#1}}{\cpointy{#1}}{#2}{#3}{\xx}{0}};
}

\newcommand{\linesym}[4]{\draw[line width=\graphline,blue] (#1,#2) -- (#3,#4); \draw[line width=\graphline,blue] (-#1,-#2) -- (-#3,-#4);}

\newcommand{\ranx}{\newpropxxpointx{\cpointx{1}}{\cpointy{1}}{1/4}{3/4}{\xx}{0}
}
\newcommand{\rany}{\newpropxxpointy{\cpointx{1}}{\cpointy{1}}{1/4}{3/4}{\xx}{0}
}

\pgfmathsetmacro{\ranxx}{0.375 * \xx + 0.125 * \cpointx{1} + 0.5 * \cpointx{2}}
\pgfmathsetmacro{\ranyy}{0.125 * \cpointy{1} + 0.5 * \cpointy{2}}

\begin{tikzpicture}[scale = \graphscale, baseline={(0,0)}]
        \draw[line width=\graphcircle,opacity=0.3] (0,0) circle (3);

        \linesym{\xx}{0}{{\cpointx{1}}}{{\cpointy{1}}};
        \linesym{\xx}{0}{{\cpointx{4}}}{{\cpointy{4}}}
        \linesym{\xx}{0}{{\cpointx{3}}}{{\cpointy{3}}};
        \linesym{\ranx}{\rany}{{\cpointx{2}}}{{\cpointy{2}}};

        \linesym{\ranxx}{\ranyy}{\propx{3}{1/2}{1/2}}{\propy{3}{1/2}{1/2}};

        \draw[line width=\graphline,blue] (\xx,0) -- (-\xx, 0);
        \draw[line width=\graphline,blue] (\propx{1}{1/2}{1/2},\propy{1}{1/2}{1/2}) -- (-\propx{1}{1/2}{1/2},\propy{1}{1/2}{1/2});
        \draw[line width=\graphline,blue] (\propx{1}{3/4}{1/4},\propy{1}{3/4}{1/4}) -- (-\propx{1}{3/4}{1/4},\propy{1}{3/4}{1/4});
        \draw[line width=\graphline,blue] (\propx{4}{1/2}{1/2},\propy{4}{1/2}{1/2}) -- (-\propx{4}{1/2}{1/2},\propy{4}{1/2}{1/2});
        \draw[line width=\graphline,blue] (\propx{4}{3/4}{1/4},\propy{4}{3/4}{1/4}) -- (-\propx{4}{3/4}{1/4},\propy{4}{3/4}{1/4});
        
        \cbdot{1};
        \cwdot{2};
        \cbdot{3};
        \cbdot{4};
        \cwdot{5};
        \cbdot{6};
        \cwdot{7};
        \cwdot{8};

        \bdotsym{\xx}{0};

        \bproppsym{1}{1/2}{1/2};
        \wproppsym{1}{1/4}{3/4};
        \wproppsym{1}{3/4}{1/4};
        
        \bproppsym{4}{1/2}{1/2};
        \wproppsym{4}{1/4}{3/4};
        \wproppsym{4}{3/4}{1/4};

        \wproppsym{3}{1/2}{1/2};

        \bdotsym{\ranxx}{\ranyy};


        \small
        \node [above right] at ({\cpointx{1}}, {\cpointy{1}}) {1};
        \node [above right] at ({\cpointx{2}}, {\cpointy{2}}) {2};
        \node [below right] at ({\cpointx{3}}, {\cpointy{3}}) {3};
        \node [below right] at ({\cpointx{4}}, {\cpointy{4}}) {4};
        \node [below left] at ({\cpointx{5}}, {\cpointy{5}}) {5};
        \node [below left] at ({\cpointx{6}}, {\cpointy{6}}) {6};
        \node [above left] at ({\cpointx{7}}, {\cpointy{7}}) {7};
        \node [above left] at ({\cpointx{8}}, {\cpointy{8}}) {8};
\end{tikzpicture}

}
    \vspace{0.2cm} \\ 
    {$n=1$}&{$n=2$} & {$n = 3$} & {$n=4$}\vspace{0.2cm}\\
    {\tiny $f = \begin{pmatrix}
        1 & 2 \\
        2 & 3
    \end{pmatrix}$} & 
    {\tiny $f = \begin{pmatrix}
        1 & 2 & 3 & 4\\
        3 & 4 & 6 & 5
    \end{pmatrix}$} & 
    {\tiny $f = \begin{pmatrix}
        1 & 2 & 3 & 4 & 5 & 6\\
        4 & 5 & 7 & 6 & 9 & 8
    \end{pmatrix}$} & 
    {\tiny$f = \begingroup
        \setlength{\arraycolsep}{2pt}  
        \begin{pmatrix}
          1 & 2 & 3 & 4 & 5 & 6 & 7 & 8\\
          5 & 6 & 8 & 7 & 10 & 9 & 12 & 11
        \end{pmatrix}
    \endgroup$
}
    \normalsize
\end{tabularx}
\end{center}
\normalsize
    \caption{Examples of \tred graph with maximal possible number of faces for $n \in \{1,2,3,4\}$.
    \label{sym-reduced-figure}}
\end{figure}
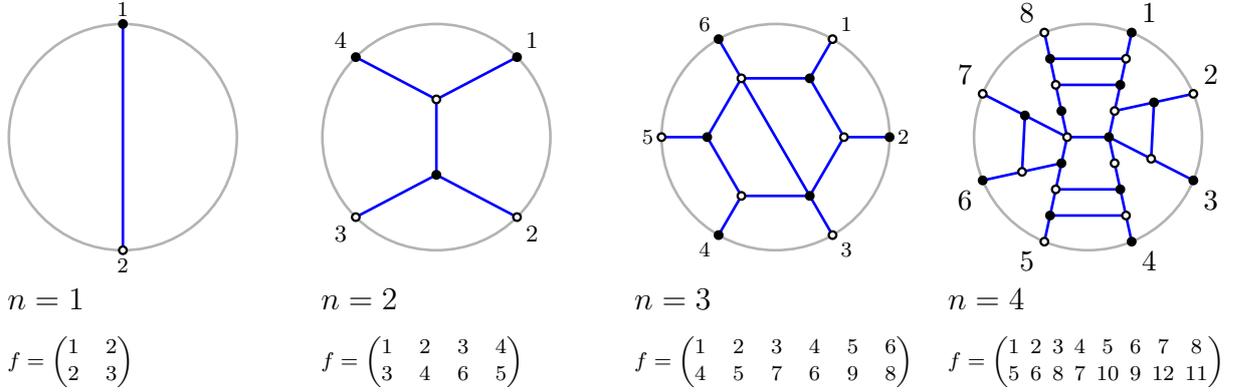

\begin{proof}[Proof of \Cref{lemma:middle-face}]

    We will prove the statement by induction on $\symmdim(f)$. 

    First, if $\symmdim(f) = 0$, then $f$ falls into the first case and we can represent it with a \tsym graph with lollipops. This graph satisfy $1 = 2 \symmdim(f) + 1$ as in case (1).

    Now, suppose $\symmdim(f) \geq 1$. Note that adding and removing fixed points and lollipops does not change the conditions (1) or (2) on $f$ and does not change the parity of faces in the graph. Thus we can assume $f$ has no fixed points. Then by \Cref{going-down-bruhat}, there exists $i$ and a \tsym permutation $f' < f$ such that $f = \symcross{i}{f'}$. If $f$ does not fall into cases (a) and (b) of \Cref{lemma:number-of-i-in} (considered below), then $\# \{i \in [2n]\st f'(i) = i+n\}$ and $\# \{i \in [2n]\st f(i) = i+n\}$ are of the same sign, so $f$ and $f'$ fall into the same case. But then if $G'$ is a \tsym graph corresponding to $f'$ that satisfies the conditions of the proposition (that exists by induction hypothesis), then $G=\symbridge{i}(G')$ is a \tsym graph that corresponds to $f$ and satisfies the desired condition (since $\# \text{faces} (G) \leq \# \text{faces} (G') + 2$). 

    We now consider cases (a) and (b) of \Cref{lemma:number-of-i-in}.
    \begin{enumerate}[(a)]
        \item Suppose $f(i) = i+n$, $f(i+1) = i+1+n$ and $f(j) \neq j+n$ for all $j \neq i, i+1$. That means $f'(i) = i+n+1,f'(i+1) = i+n$ and $f'(j) = f(j)$ for $j \neq i, i+1$. Since $f'$ does not only have fixed points, there exists $j \in [2n]$ such that $fs_j \in \Bnn$ and $fs_j < j$, i.e. $j \leq f'(j) < f'(j+1) \leq j+1 + 2n$ (in particular $i \neq 1$). But since the values of $f(j)$ and $f'(j)$ differ at most by $1$ for each $j$, we also get $j \leq f'(j) < f'(j+1) \leq j+1 + 2n$, which means that $fs_j < f$. Then by \Cref{hultman-cross-equals-sym-cross}, we get that there exists a \tsym permutation $f'' < f$ such that $f = \symcross{j}{f''}$. And since $j \neq i$, we get that $f(j) \neq j+n$ or $f(j+1) \neq j+1 +n$, so $f$ and $f''$ do not fall into cases (a) or (b). Therefore if $G''$ is a \tsym graph corresponding to $f''$ satisfying the conditions, then we can take $G = \symbridge{j}(G'')$.
        
        \item $n = 1$ and $f(1) = 2, f(2) = 1$. Then we can just take $G$ to be the graph with one bridge in the middle, between the two vertices as in \Cref{sym-reduced-figure}.
    \end{enumerate}

\end{proof}

After \Cref{symm-bounding-number-of-faces-lemma} and \Cref{lemma:middle-face}, the following definition is natural
\begin{definition}
    \label{tred-definition}
    Say that a \tsym graph $G$ parameterizing a cell $\celllgn{f}$ is \term{\tred} if it satisfies
    $$\# \text{faces} (G) = \begin{cases}
        2\symmdim(f), &\exists i, f(i) = i+n;\\
        2\symmdim(f) + 1, & \text{otherwise}.
    \end{cases}$$
    In other words, a \tred graph is the graph parameterizing $\celllgn{f}$ with the minimal possible number of faces.
\end{definition}

The next statement follows directly from \Cref{lemma:middle-face}.
\begin{corollary}
    \label{cor:tred-graph-exists}
    For every $f \in \sbn$, the cell $\celllgn{f}$ can be parametrized with a \tred graph $G$.
\end{corollary}

\begin{proposition}
    \label{prop:reduced-and-tsym-equals-tred}
    If a \tsym graph is reduced, then it is also \tred.
\end{proposition}
\begin{proof}
    Suppose a graph $G$ is \tsym and reduced and corresponds to a permutation $f \in \sbn$.  First, replicating the proof from Theorem 4.2 in \cite{karpman-su-combinatorics}, we get that every $X \in \celllgn{f}$ can be represented with a symmetric weighting on $G$. That is, $G$ parametrizes $\celllgn{f}$. And since $G$ has the minimal number of faces out of all graphs parameterizing $\cell{f}$, it also has the minimal number of faces out of \tsym graphs parameterizing $\celllgn{f}$. Therefore, $G$ is \tred.
\end{proof}

\begin{proposition}
    \label{prop:condition-symm-and-reduced}
    Suppose $f \in \sbn$. Then there exists a \tsym reduced graph parameterizing $\celllgn{f}$ if and only if 
    $$\# \{i\in [2n]\st f(i) = i+n\} \in \{0,2\}.$$
\end{proposition}
\begin{remark}
    Note that the points $i$ such that $f(i) \neq i + n$ come in pairs (if $f(i) = j+n$, then $f(j) = i + n$). Therefore, $\# \{i\in [2n]\st f(i) = i+n\}$ is always even.
\end{remark}
\begin{proof}[Proof of \Cref{prop:condition-symm-and-reduced}]

    By \Cref{cor:tred-graph-exists}, there exists a \tsym graph $G$ parameterizing $\celllgn{f}$. By \Cref{prop:reduced-and-tsym-equals-tred}, $G$ is reduced if and only if $\celllgn{f}$ can be parametrized by a graph that is \tsym and reduced. Therefore, we wish to show that $G$ is reduced if and only if $$\# \{i\in [2n]\st f(i) = i+n\} = 0 \text{ or } 2.$$

    By \Cref{tred-definition}, we have
    \begin{align}
        \label{eq:faces-1-proof}
        \# \text{faces} (G) = \begin{cases}
        2\symmdim(f), &\exists i, f(i) = i+n;\\
        2\symmdim(f) + 1, & \text{otherwise}.
    \end{cases}
    \end{align}
    
    On the other hand, applying \Cref{d-vs-symmdim} $G$ is reduced if and only if 
    \begin{align}
        \label{eq:faces-2-proof}
        \# \text{faces} (G) = d(f) + 1 = 2 \symmdim(f) - \frac{\# \{i \in [2n] \st f(i) = i+n\}}{2} + 1.
    \end{align}
    Combining \Cref{eq:faces-1-proof} and \Cref{eq:faces-2-proof}, we obtain that $G$ is reduced if and only if 
    $$\# \{i \in [2n] \st f(i) = i+n\} = \begin{cases}
        2, &\exists i, f(i) = i+n;\\
        0, & \text{otherwise}
    \end{cases}.$$
    This completes the proof.

\end{proof}

\begin{lemma}
    \label{lemma:number-of-i-in}
    Suppose $f'<f$ are two \tsym permutations such that $f = \symcross{i}{f'}$. Then 
    $$\# \{i \in [2n]\st f'(i) = i+n\} \leq \# \{i \in [2n]\st f(i) = i+n\}.$$
    Moreover, $\# \{i \in [2n]\st f'(i) = i+n\}$ and $\# \{i \in [2n]\st f(i) = i+n\}$ have the same sign (i.e. are either both zero or both positive) unless one of the following holds
    \begin{enumerate}[(a)]
        \item $f(i) = i+n$, $f(i+1) = i+1+n$ and $f(j) \neq j+n$ for all $j \neq i, i+1$. 
        \item $n = 1$ and $f(1) = 2, f(2) = 1$.
    \end{enumerate}
    In both of these cases, we get 
    $$\# \{i \in [2n]\st f'(i) = i+n\} = 0 \text{ and } \# \{i \in [2n]\st f(i) = i+n\} = 2.$$
\end{lemma}
\begin{proof}
    The case $n=1$ is obvious, so we assume $n \geq 2$. Notice that $\# \{j \in [2n]\st f'(j) = j+n\} \neq \# \{j \in [2n]\st f(j) = j+n\}$ only if there exists $j$ such that $f(j) = j+n$ and $f'(j) \neq j+n$ or vice versa. Since $f = \symcross{i}{f'}$, we get $j=i$ or $j=i+1$. Without loss of generality, we may assume that $j = i$. We now consider cases.

    \begin{enumerate}
        \item Suppose $f'(i) = i+n$. Let $f'(i+1) = \ell  + n$, which implies $f'(\ell) = i + n + 1$. Since $s_i f' > f'$, we have $\ell \neq i+1$. Then we have $f' = s_{i+n} f s_i$, i.e.
        $$f(j) = \begin{cases}
            i + n + 1, & j = i;\\
            i + n, & j = i + 1\;\\
            \ell + n, & j = \ell;\\
            f'(j),& \text{otherwise}.
        \end{cases}$$
        In this case, we can clearly see that  $\# \{j \in [2n]\st f'(j) = j+n\} = \# \{j \in [2n]\st f(j) = j+n\}$.

        \item Suppose $f(i) = i+n$. Then we have two subcases: 
        \begin{enumerate}[(i)]
            \item $f(i+1) = \ell + n \neq i + 1 + n$. The logic for this case is similar to the above.
            
            \item $f(i+1) = i + n + 1$. In this case, we get $\# \{i \in [2n]\st f'(i) = i+n\}= \# \{i \in [2n]\st f(i) = i+n\} - 2$, which leads to either them having the same sign or case (a).
        \end{enumerate}
    \end{enumerate}
\end{proof}

\begin{corollary}
    Suppose $f \in \sbn$ such that $\celllgn{f}$ can be parametrized with a graph $G$ that is both \tsym and reduced. Then
    $$\codim \celllgn{f} = \symmell(f) \geq n-1.$$

    Conversely, for each $n$, there exists $f\in \sbn$ such that $\symmell(f) = n-1$ and $\celllgn{f}$ can be parametrized with a graph that is both \tsym and reduced.
\end{corollary}
\begin{proof}

    For the first part of the statement, suppose there exists a graph representing $f$ that is both \tsym and reduced. Then, by \Cref{prop:condition-symm-and-reduced}, we get 
    $$\# \{i\in [2n]\st f(i) = i+n\} = 0 \text{ or } 2.$$
    Now, applying \Cref{eq:symm-alignments-i-in}, we obtain
    \begin{align*}
        \symmell(f) = \# \{\text{\tpairals\}} \geq \\
        \# \{\text{\tcenterals\}}= \frac{2n - \#\{i \st f(i) = i+n\}}{2} \geq n-1.
    \end{align*}
    This proves the first part of the statement. 

    For the second part, it is enough to show that for any given $n$, there exists $f \in \sbn$ such that 
    $$\symmell(f) = n-1 \text{ and } \# \{i\in [2n]\st f(i) = i+n\} = 0 \text{ or } 2.$$
    For that define $f$ as follows 
    $$f(i) = \begin{cases}
        i + n, & i = 1 \text{ or } 2;\\
        2s + n & i = 2s-1, s \neq 1;\\
        2s + n - 1 & i = 2s, s \neq 1.
    \end{cases}$$
    Then we can see that $f \in \sbn$ and $\# \{i\in [2n]\st f(i) = i+n\} = 2$. Moreover, we can see that the only alignments in $f$ are the ones formed by $2s-1$ and $2s$ for $s \neq 1$ (and all of these alignments are \tsym, i.e. symmetric to themselves). Therefore for this $f$, we get 
    $$\symmell(f) = \#\{\text{symmetric pairs of alignments}(f)\} = n - 1.$$
    This concludes the proof.
    
\end{proof}

\begin{remark}
    In particular, we get that for $n \geq 2$, the top cell of $\lgnn$ cannot be parametrized with a graph that is both \tsym and reduced.
\end{remark}

\subsection{Moves between \tsym graphs} 

Recall that every two reduced plabic graphs, corresponding to the same permutation, can be connected through a series of local moves (see \Cref{plabic-graphs-subsection}). We introduce similar moves for \tsym plabic graphs and conjecture that \tred graphs corresponding to the same \tsym permutation can be connected though \tsym moves (see \Cref{conjecture: tsym-moves}).

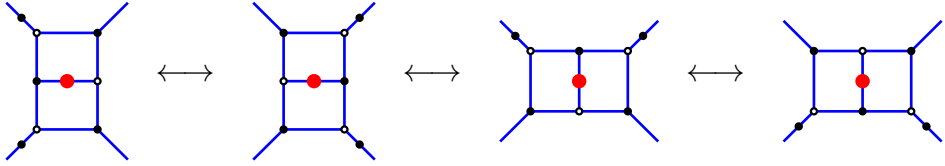
\begin{figure}
    \centering
    \newcommand{\movesline}{1}
\newcommand{\movesscale}{0.6}
\newcommand{\lo}{1.6}
\newcommand{\sh}{1}
\newcommand{\dd}{1}
\begin{tabular}{ccccccc}
   {    \begin{tikzpicture}[scale = 0.4,baseline={(0,0)}]
        \draw[line width=\movesline, blue] (\sh, \lo) -- (\sh, -\lo);
        \draw[line width=\movesline, blue] (\sh, \lo) -- (-\sh, \lo);
        \draw[line width=\movesline, blue] (-\sh,-\lo) -- (-\sh, \lo);
        \draw[line width=\movesline, blue] (-\sh, -\lo) -- (\sh, -\lo);
        \draw[line width=\movesline, blue] (-\sh, 0) -- (\sh, 0);
        \draw[line width=\movesline, blue] (\sh, \lo) -- (\sh + \dd, \lo + \dd);
        \draw[line width=\movesline, blue] (-\sh, \lo) -- (-\sh - \dd, \lo + \dd);
        \draw[line width=\movesline, blue] (\sh, -\lo) -- (\sh + \dd,-\lo - \dd);
        \draw[line width=\movesline, blue] (-\sh, -\lo) -- (-\sh - \dd,-\lo - \dd);
        
        \bdot{\sh}{\lo}; \bdot{\sh}{-\lo}, \wdot{-\sh}{\lo}; \wdot{-\sh}{-\lo}; \wdot{\sh}{0}; \bdot{-\sh}{0};

        \bdot{-\sh -0.5*\dd}{\lo + 0.5*\dd};
        \bdot{-\sh -0.5*\dd}{-\lo - 0.5*\dd};
        \cntr;
    \end{tikzpicture}}  & $\longleftrightarrow$ & {        \begin{tikzpicture}[scale = 0.4,baseline={(0,0)}]
        \draw[line width=\movesline, blue] (\sh, \lo) -- (\sh, -\lo);
        \draw[line width=\movesline, blue] (\sh, \lo) -- (-\sh, \lo);
        \draw[line width=\movesline, blue] (-\sh,-\lo) -- (-\sh, \lo);
        \draw[line width=\movesline, blue] (-\sh, -\lo) -- (\sh, -\lo);
        \draw[line width=\movesline, blue] (-\sh, 0) -- (\sh, 0);
        \draw[line width=\movesline, blue] (\sh, \lo) -- (\sh + \dd, \lo + \dd);
        \draw[line width=\movesline, blue] (-\sh, \lo) -- (-\sh - \dd, \lo + \dd);
        \draw[line width=\movesline, blue] (\sh, -\lo) -- (\sh + \dd,-\lo - \dd);
        \draw[line width=\movesline, blue] (-\sh, -\lo) -- (-\sh - \dd,-\lo - \dd);
        
        \wdot{\sh}{\lo}; \wdot{\sh}{-\lo}, \bdot{-\sh}{\lo}; \bdot{-\sh}{-\lo}; \bdot{\sh}{0}; \wdot{-\sh}{0};

        \bdot{\sh +0.5*\dd}{\lo + 0.5*\dd};
        \bdot{\sh +0.5*\dd}{-\lo - 0.5*\dd};
        \cntr;
    \end{tikzpicture}} & $\longleftrightarrow$ & {        \begin{tikzpicture}[scale = 0.4,baseline={(0,0)}]
        \draw[line width=\movesline, blue] (\lo,\sh) -- (-\lo, \sh);
        \draw[line width=\movesline, blue] (\lo,\sh) -- (\lo,-\sh);
        \draw[line width=\movesline, blue] (-\lo,-\sh) -- (\lo,-\sh);
        \draw[line width=\movesline, blue] (-\lo, -\sh) -- (-\lo, \sh);
        \draw[line width=\movesline, blue] (0, -\sh) -- (0,\sh);
        \draw[line width=\movesline, blue] (\lo,\sh) -- (\lo + \dd, \sh + \dd);
        \draw[line width=\movesline, blue] (\lo,-\sh) -- (\lo + \dd, -\sh - \dd);
        \draw[line width=\movesline, blue] (-\lo, \sh) -- (-\lo - \dd, \sh + \dd);
        \draw[line width=\movesline, blue] ( -\lo,-\sh) -- (-\lo - \dd, -\sh - \dd);
        
        \wdot{\lo}{\sh}; \wdot{-\lo}{\sh}, \bdot{\lo}{-\sh}; \bdot{-\lo}{-\sh}; \bdot{0}{\sh}; \wdot{0}{-\sh};

        \bdot{\lo + 0.5*\dd}{\sh +0.5*\dd};
        \bdot{-\lo - 0.5*\dd}{\sh +0.5*\dd};
        \cntr;
    \end{tikzpicture}}  & $\longleftrightarrow$ & {        \begin{tikzpicture}[scale = 0.4,baseline={(0,0)}]
        \draw[line width=\movesline, blue] (\lo,\sh) -- (-\lo, \sh);
        \draw[line width=\movesline, blue] (\lo,\sh) -- (\lo,-\sh);
        \draw[line width=\movesline, blue] (-\lo,-\sh) -- (\lo,-\sh);
        \draw[line width=\movesline, blue] (-\lo, -\sh) -- (-\lo, \sh);
        \draw[line width=\movesline, blue] (0, -\sh) -- (0,\sh);
        \draw[line width=\movesline, blue] (\lo,\sh) -- (\lo + \dd, \sh + \dd);
        \draw[line width=\movesline, blue] (\lo,-\sh) -- (\lo + \dd, -\sh - \dd);
        \draw[line width=\movesline, blue] (-\lo, \sh) -- (-\lo - \dd, \sh + \dd);
        \draw[line width=\movesline, blue] ( -\lo,-\sh) -- (-\lo - \dd, -\sh - \dd);
        
        \bdot{\lo}{\sh}; \bdot{-\lo}{\sh}, \wdot{\lo}{-\sh}; \wdot{-\lo}{-\sh}; \wdot{0}{\sh}; \bdot{0}{-\sh};

        \bdot{\lo + 0.5*\dd}{-\sh -0.5*\dd};
        \bdot{-\lo - 0.5*\dd}{-\sh -0.5*\dd};
        \cntr;
    \end{tikzpicture}}
\end{tabular}
    \label{double-square-move-figure}
    \caption{Double square move}
\end{figure}

\begin{figure}
    \centering
    \newcommand{\movesline}{1}
\newcommand{\movesscale}{0.6}
\newcommand{\lo}{1.6}
\newcommand{\dd}{1}
\begin{tabular}{ccc}
   {    \begin{tikzpicture}[scale = 0.4,baseline={(0,0)}]
        \draw[line width=\movesline, blue] (2.2*\lo, 0) -- (-2.2*\lo, 0);
        \draw[line width=\movesline, blue] (1.5*\lo, 0) -- (2.2*\lo, 0.8);
        \draw[line width=\movesline, blue] (1.5*\lo, 0) -- (2.2*\lo, -0.8);
        \draw[line width=\movesline, blue] (-1.5*\lo, 0) -- (-2.2*\lo, 0.8);
        \draw[line width=\movesline, blue] (-1.5*\lo, 0) -- (-2.2*\lo, -0.8);
        \bdot{0.6*\lo}{0}; \wdot{1.5*\lo}{0}; \wdot{-0.6*\lo}{0}; \bdot{-1.5*\lo}{0}; 
        \cntr;
    \end{tikzpicture}}  & $\longleftrightarrow$ & {       
    \begin{tikzpicture}[scale = 0.4,baseline={(0,0)}]
        \draw[line width=\movesline, blue] (1.8*\lo, 0) -- (-1.8*\lo, 0);
        \draw[line width=\movesline, blue] (1.1*\lo, 0) -- (1.8*\lo, 0.8);
        \draw[line width=\movesline, blue] (1.1*\lo, 0) -- (1.8*\lo, -0.8);
        \draw[line width=\movesline, blue] (-1.1*\lo, 0) -- (-1.8*\lo, 0.8);
        \draw[line width=\movesline, blue] (-1.1*\lo, 0) -- (-1.8*\lo, -0.8);
        \wdot{1.1*\lo}{0}; \bdot{-1.1*\lo}{0};
        \cntr;
    \end{tikzpicture}
    } 
\end{tabular}
    \label{double-m2-move-figure}
    \caption{Double (M2) move}
\end{figure}
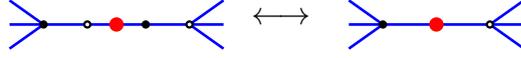

First, if regular move does not involve an edge passing through the center of the rotation (i.e. the edge that is symmetric to itself), we may apply this move with the move rotationally symmetric to it at the same time. However, if a move includes an edge passing through the center, we might not be able to do both this move and the move symmetric to it at the same time. Because of this we define two new moves, demonstrated in \Cref{double-square-move-figure} and \Cref{double-m2-move-figure}. We say that \textbf{a \tsym move} is either a combination of a pair of regular moves symmetric to each other or one of the two new moves.

One can see that each \tsym move can be represented as a sequence of regular moves. Therefore, we get the following.

\begin{lemma}
    Suppose $G$ and $G'$ are two \tsym plabic graphs related by a sequence of \tsym moves. Then $f_G = f_{G'}$. 
\end{lemma}

We expect that it will extend to the following.

\begin{conjecture}
    \label{conjecture: tsym-moves}
    Suppose $G$ and $G'$ are two \tred graphs such that $f_G = f_{G'}$. Then $G$ and $G'$ are related though a sequence of \tsym moves. 
\end{conjecture}

\section{Cell structure} \label{cell-structure-section}


In this section, we will describe the cell structure of the TNN Lagrangian Grassmannian. We already established that the non-empty cells of 
$\lgnn$ correspond to \tsym permutations. Additionally, it was shown in \Cref{cell-topology} that each \tsym permutation is associated with a non-empty cell of 
$\lgnn$ of a dimension $\symmdim(f)$. Our next focus is to investigate the relationships between these cells. Specifically, we aim to determine the conditions under which two \tsym bounded affine permutations 
$f$ and $g$ satisfy 
$\celllgn{f} \subset \overline{\celllgn{g}}$. We claim that this relationship can be characterized as follows.
\begin{theorem}
    \label{theorem-cell-structure}
    Suppose $f$ and $g$ are \tsym bounded affine permutations.

    Then if $f \sleq g$, then $\celllgn{f} \subset \overline{\celllgn{g}}$. Otherwise, we have $\celllgn{f} \cap \overline{\celllgn{g}} = \varnothing$.
\end{theorem}
\begin{figure}
    \centering
    \input{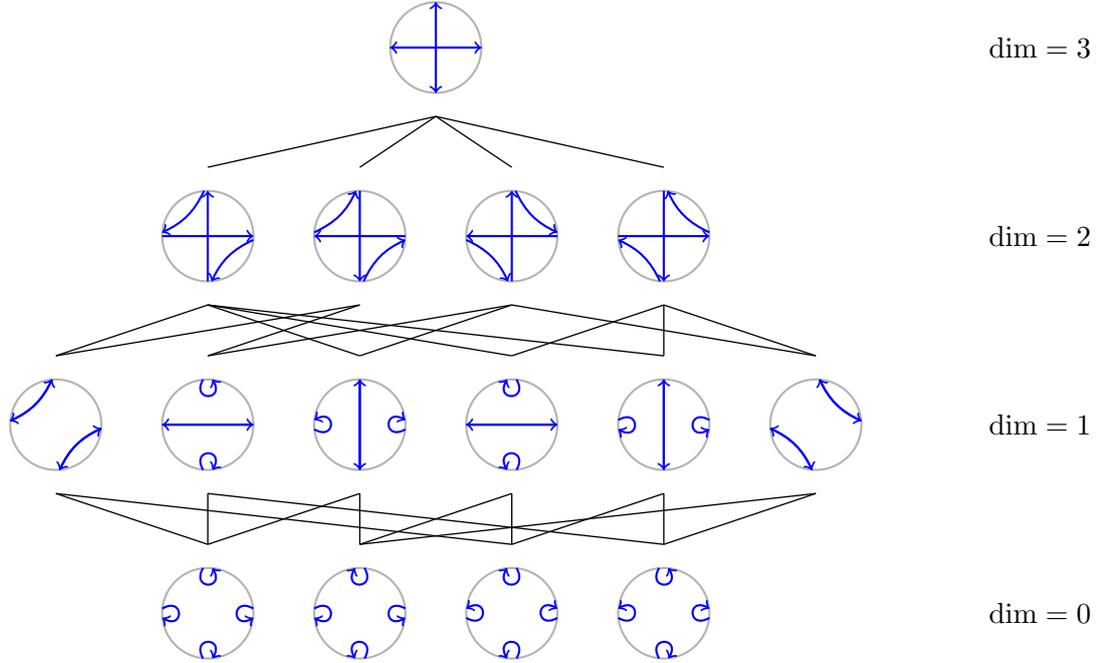}
    \caption{The Hasse diagram of the poset of \tsym permutations for $n=2$, ordered by inclusion of the closures of the corresponding cells. Here, we use permutations instead of affine permutations by taking everything  modulo $2n$ (fixed points are pictured as counterclockwise or clockwise loops depending on whether $f(i) = i$ or $f(i) = i+2n$). We list the dimensions of the corresponding cells on the right.}
    \label{cell-structure-figure}
\end{figure}

\begin{example}
    This theorem is illustrated in \Cref{cell-structure-figure}, which shows the cell closure poset of \( \lgrnn(2,4) \). Specifically, the diagram depicts the Hasse diagram of the poset of \tsym permutations for \( n = 2 \) with order relation given by $\sleq$.
\end{example}

The proof of \Cref{theorem-cell-structure} is based on the following topological lemma.\footnote{This fact was stated and proved by Pavel Galashin. We are grateful for his permission to include it here.} Here, $D^n$ is the closed ball of dimension $n$ and $S^{n-1}$ is its boundary.

\begin{lemma}
    \label{lemma: main-topology-lemma}
    Let $\mainop: D^n \to D^n$ be a continuous involution such that $\mainop(S^{n-1}) = S^{n-1}$. Suppose $x_0 \in S^{n-1}$ is a fixed point of $\mainop$. Then there exists a sequence $x_j \in \interior(D^n)$ of fixed points of $\mainop$ such that $x_j \to x_0$.
\end{lemma}
\begin{proof} 
    It is enough to show that for every open neighborhood $U$ of $x_0$, there exists an $\mainop$-fixed point in $U \cap \interior(D^n)$. Fix such $U$. We show the statement in two steps.

    First, we take $U_n = U$ and find open balls $U_{n-1}, \ldots, U_0$ such that for each $i \in \{0, \ldots, n\}$, we have $x_0 \in U_i$ and for each $i \in \{1, \ldots, n\}$,
    $$U_{i-1} \subset \mainop^{-1}(U_i) \cap U_i.$$
    (We can choose such a ball $U_{i-1}$ inductively since $\mainop^{-1}(U_i) \cap U_i$ is an open neighborhood of $x_0$.) Equivalently, we get that for every $i \in \{0, \ldots, n-1\}$, 
    \begin{align}
        \label{eq:U-i-def}
        U_i \subset U_{i+1}, \quad \mainop(U_i) \subset U_{i+1}.
    \end{align}

    Now, we will show inductively that for every $i \in \{0, \ldots, n-1\}$, there exists a continuous map $g_i: S^{i} \to D^n$ such that 
    $$\text{for all } x \in S^i, g_i(-x) = \mainop(g_i(x)).$$
    and
    $$g_i(S^i) \subset U_{i+1} \cap \interior(D^n).$$
    First, for $n = 0$, $S^0 = \{1,-1\}$, so we can define $g_0(1) = y$ and $g_0(-1) = \mainop(y)$, where $y \in U_0 \cap \interior(D^n) \subset U_1 \cap \interior(D^n)$, so by \Cref{eq:U-i-def}, we have $\mainop(y) \in U_1 \cap \interior(D^n)$.

    Now, if $g_i$ is already constructed on $S^i$, notice that it is null-homotopic in $U_{i+1} \cap \interior(D^n)$. We can use the homotopy to extend $g_i$ to a map $g_{i+1}$ on the upper hemisphere $S^{i+1}_{\text{upper}}$. Now, on the lower semisphere $S^{i+1}_{\text{lower}}$, we define $g_{i+1}$ by 
    $$g_{i+1}(-x) = \mainop(g_{i+1}(x)) \in \mainop(U_{i+1} \cap \interior(D^n)) \subset U_{i+2} \cap \interior(D^n) \text{ for } x \in S^{i+1}_{\text{upper}}.$$
    Therefore $g_i$ satisfies the conditions above. 

    Next, let $g = g_{n-1}: S^{n-1} \to U \cap \interior(D^n)$. If $\mainop$ has no fixed points in $U \cap \interior(D^n)$, we may define $h: S^{n-1} \to S^{n-1}$ as follows
    $$h(x) = \frac{g(x) - \mainop(g(x))}{\|g(x) - \mainop(g(x))\|}.$$
    Then we can notice that $h$ satisfies $h(-x) = -h(x)$, and therefore has an odd degree. But on the other hand, since $g$ is null-homotopic on $S^n$, 
    so is $h$, so it must have an even degree. This implies that $\mainop$ must have fixed points in $U \cap \interior(D^n)$.
\end{proof}

We are now ready to prove the main theorem of the section.
\begin{proof}[Proof of \Cref{theorem-cell-structure} ]

    First, suppose $f \leq g$. By \Cref{postnikov-theorem}, this means that 
    \begin{align}
        \label{eq: cell-inside-cell}
        \cell{f} \subset \overline{\cell{g}}.
    \end{align}
    Notice that since $g \in \sbn$, $\mainop$ restricts to an involution $\mainop: \overline{\cell{g}} \to \overline{\cell{g}}$. Moreover, by Theorem 1.1 in \cite{galashin-karp-lam-regularity-theorem}, $\overline{\cell{g}}$ is homeomorphic to a closed ball $D^{d(g)}$. Therefore, we may apply \Cref{lemma: main-topology-lemma}. Since $\mainop$-fixed points are exactly points in $\lgnn$, we obtain that 
    \begin{align}
        \label{eq: closure-intersection}
        \overline{\cell{g}} \cap \lgnn \subset \overline{\cell{g} \cap \lgnn}.
    \end{align}
    On the other hand, combining  \Cref{eq: cell-inside-cell} and \Cref{eq: closure-intersection}, we get that
    $$ \celllgn{f} = \cell{f} \cap \lgnn \subset \overline{\cell{g}} \cap \lgnn \subset \overline{\cell{g} \cap \lgnn} = \overline{\celllgn{g}}.$$
    This proves the first part of the theorem.

    For the second part, suppose $$\celllgn{f} \cap \overline{\celllgn{g}} \neq \varnothing.$$
    Since $\celllgn{f} \subset \cell{f}$ and $\celllgn{g} \subset \cell{g}$, we also get
    $$\cell{f} \cap \overline{\cell{g}} \neq \varnothing.$$
    From the cell structure of the TNN Grassmannian (see \Cref{postnikov-theorem}), we know that 
    $$\overline{\cell{g}} = \bigsqcup_{h \leq g} \cell{h}.$$
    And since all cells are disjoint, we get that $f$ must satisfy $f \leq g$ (so that $\cell{f}$ is one of the cells in the union). But since both $f$ and $g$ are \tsym, this means that $f \sleq g$.

    For the other part of the theorem, suppose 
    
\end{proof}

\subsection{Totally Nonnegative Lagrangian Grassmannian is a ball}

In this subsection, we prove the following.

\begin{theorem}
    \label{thm: lgnn-is-a-ball}
    For every $n$, $\lgnn$ is homeomorphic to a closed ball of dimension ${n+1} \choose 2$.
\end{theorem}

The proof of this statement follows the approach from \cite{galashin-karp-lam-ball}. First, we define the left cyclic shift
$$S: \RR^{2n} \to \RR^{2n}: (x_1, \ldots, x_{2n}) \mapsto (x_2, \ldots, x_{2n}, (-1)^{n-1} x_1)$$
and the right cyclic shift
$$S^T: \RR^{2n} \to \RR^{2n}: (x_1, \ldots, x_{2n}) \mapsto ((-1)^{n-1} x_n, x_1, x_2, \ldots, x_{2n-1}).$$
\begin{lemma}[\cite{galashin-karp-lam-ball}, Lemma 3.5]
    \label{lemma-3.5-ball}
    For all $X \in \grnn$ and $t>0$, we have
    $$X \cdot \exp(t (S + S^T)) \in \grnn.$$
\end{lemma}
We can easily extend the result to $\lgnn$.

\begin{lemma}
    \label{lemma: lgnn-exp}
    Suppose $X \in \lgnn$ and $t>0$. Then 
    $$X \cdot \exp(t (S + S^T)) \in \lgnp.$$
\end{lemma}
\begin{proof}
    We are following the proof in \cite{pashas} (Lemma 7.2). 

    First, it follows from \Cref{lemma-3.5-ball}, it is enough to show that $X \cdot \exp(t (S + S^T)) \in \lgr$ for all $X \in \lgr$. Since $\lgr = \{X \in \gr\st X R X^T = 0\}$, it is enough to show that $\exp(t (S + S^T))$ belongs to the symplectic group $\Spn$ corresponding to $R$. And since $\exp$ restricts to a map from the symplectic Lie algebra $\Spalgn$ to $\Spn$, and $\Spalgn$ is closed under multiplication by a scalar, it is enough to show that $S + S^T \in \Spalgn$. That is, we wish to show that $(S + S^T) R = - R (S + S^T)$. We can compute that 
    $$S + S^T = \begin{pmatrix}
    \begin{matrix}
        0 & 1 &  & \\
        1 & 0 & 1 & \\
        & 1 & 0 & &\\
        &&& \ddots &&\\
        & &&   & 0 & 1\\
        & &&   & 1 & 0
    \end{matrix}
    & \rvline &  
    \begin{matrix}
        &&&&& (-1)^{n-1}\\
        &&&&&\\
        &&&&&\\
        &&&&&\\
        &&&&&\\
        1&&&&&\\
    \end{matrix}
    \\
    \hline
    \begin{matrix}
        &&&&& 1\\
        &&&&&\\
        &&&&&\\
        &&&&&\\
        &&&&&\\
        (-1)^{n-1}&&&&&\\
    \end{matrix} & \rvline &
    \begin{matrix}
        0 & 1 &  && \\
        1 & 0 & 1 & &\\
        & 1 & 0 & &&\\
        &&& \ddots &&\\
        & &&   & 0 & 1\\
        & &&   & 1 & 0
    \end{matrix}
\end{pmatrix} 
$$

$$R = \begin{pmatrix}
    & \rvline &  
    \begin{matrix}
        -1 &  &   \\
         & 1 &   \\
        && \ddots &\\
        & & &  (-1)^{n}
    \end{matrix}
    \\
    \hline
    \begin{matrix}
        1 &  &   \\
         & -1 &   \\
        && \ddots &\\
        & & &  (-1)^{n-1}
    \end{matrix} & \rvline &
    
\end{pmatrix},$$

So we get
$$(S + S^T)R = \begin{pmatrix}
    \begin{matrix}
        &&& 1\\
        &&&\\
        &&&\\
        1&&&\\
    \end{matrix}
    & \rvline &
    \begin{matrix}
        0 & 1 &   \\
        -1 & 0 & -1  \\
        & 1 & 0 \\
        &&& \ddots\\
    \end{matrix}
    \\
    \hline
    \begin{matrix}
        0 & 1 &   \\
        -1 & 0 & -1  \\
        & 1 & 0 \\
        &&& \ddots\\
    \end{matrix} & \rvline &
    \begin{matrix}
        &&& (-1)^n\\
        &&&\\
        &&&\\
        (-1)^n&&&\\
    \end{matrix}
\end{pmatrix} $$

Notice that the product is symmetric, so $\big( (S+S^T) R \big)^T = (S+S^T) R$. On the other hand, we have $\big( S+S^T) R \big)^T = R^T (S + S^T)^T = -R (S + S^T)$. This shows that $(S + S^T) R = - R(S + S^T)$, which completes the proof.
\end{proof}

\begin{proof}[Proof of \Cref{thm: lgnn-is-a-ball}]
    Using \Cref{lemma: lgnn-exp}, the proof of the main theorem translates directly from the proof of \cite[Theorem~1.1]{galashin-karp-lam-ball}. 
\end{proof}


\section*{Acknowledgements}
The author would like to thank Pavel Galashin for the help in finding this problem and for guidance on every step and Thomas Martinez for many helpful discussions.

\bibliographystyle{alpha}
\bibliography{biblio}

\end{document}